\documentclass[10pt, twoside, open=right, pointlessnumbers, a4paper]{article}

\usepackage{amsmath}
\usepackage{titlesec}
\usepackage{wrapfig}
\usepackage{graphicx}
\usepackage{color}
\usepackage{xcolor, framed}
\usepackage{soul}
\usepackage[latin1]{inputenc}

\usepackage{multicol}
\usepackage{needspace}

\makeatletter
\DeclareRobustCommand{\em}{%
	\@nomath\em \if b\expandafter\@car\f@series\@nil
	\normalfont \else \slshape \fi}

\usepackage{tocloft}
\usepackage{amsmath}
\usepackage{amssymb}
\usepackage{amsfonts}
\usepackage{upgreek}
\usepackage{amsthm}
\usepackage{mathtools}

\usepackage{mathrsfs}
\usepackage{pifont}
\DeclareMathAlphabet{\mathpzc}{OT1}{pzc}{m}{it}
\usepackage{bbm}
\usepackage{bbding}
\usepackage{setspace}
\usepackage{mathdots}
\usepackage{relsize}
\usepackage{tocloft}

\numberwithin{equation}{section}
\usepackage{mathtools}
\mathtoolsset{showonlyrefs}
\numberwithin{equation}{section}
\usepackage{enumitem} 

\newenvironment{myenumerate}{\begin{enumerate}[topsep=2pt,parsep=2pt,partopsep=2pt,itemsep=0pt,label={\normalfont(\alph*)}]\itemsep0pt}{\end{enumerate}}

\newtheoremstyle{style1}
  {13pt}
  {13pt}
  {}
  {}
  {\normalfont\bfseries}
  {.}
  {.5em}
  {}
\theoremstyle{style1}

\newtheorem{definition}[equation]{Definition}
\newtheorem{example}[equation]{Example}
\newtheorem{remark}[equation]{Remark}
\newtheorem{preremarks}[equation]{Remarks}

\newtheoremstyle{style2}
  {13pt}
  {13pt}
  {\slshape}
  {}
  {\normalfont\bfseries}
  {.}
  {.5em}
  {}

\theoremstyle{style2}

\newtheorem{lemma}[equation]{Lemma}
\newtheorem{theorem}[equation]{Theorem}
\newtheorem{proposition}[equation]{Proposition}
\newtheorem{corollary}[equation]{Corollary}

\usepackage{hhline}
\usepackage{colortbl}
\usepackage{tabularx}
\usepackage{longtable}
\usepackage{pgf}
\usepackage{tikz}
\usepackage{tikz-cd}
\usepackage{fancybox}
\usetikzlibrary{arrows,shapes,trees,matrix}
\usetikzlibrary{decorations.markings}
\usetikzlibrary{knots}
\usepackage{pictexwd}
\usepackage[hidelinks]{hyperref}

\newcommand{\myforall}{\quad \text{for all}\quad }

\newcommand{\id}{\operatorname{id}}

\newcommand{\Hom}{\operatorname{Hom}}
\newcommand{\cat}[1]{\mathcal{#1}}
\newcommand{\catf}[1]{{\normalfont\textbf{#1}}}

\newcommand{\FinGrpd}{\catf{FinGrpd}}

\newcommand{\Vect}{\catf{Vect}}
\newcommand{\FinVect}{\catf{FinVect}}
\newcommand{\RepGrpd}[1]{\catf{VecBun}_{#1}\catf{Grpd}}
\newcommand{\TRepGrpd}{\catf{2VecBunGrpd}}

\newcommand{\Mod}{\catf{Mod}}

\newcommand{\TwoVect}{\catf{2Vect}}

\newcommand{\Par}{\operatorname{Par}}
\newcommand{\dir}[2]{to [out=#1, in=#2]}

\newcommand{\Map}{\operatorname{Map}}

\newcommand{\Aut}{\operatorname{Aut}}
\newcommand{\End}{\operatorname{End}}
\newcommand{\td}{\text{\normalfont d}}

\let\Phi\undefined\DeclareMathSymbol{\Phi}{\mathalpha}{letters}{"08}
\let\Psi\undefined\DeclareMathSymbol{\Psi}{\mathalpha}{letters}{"09}
\let\Sigma\undefined\DeclareMathSymbol{\Sigma}{\mathalpha}{letters}{"06}
\let\Xi\undefined\DeclareMathSymbol{\Xi}{\mathalpha}{letters}{"04}
\let\Pi\undefined\DeclareMathSymbol{\Pi}{\mathalpha}{letters}{"05}
\let\Gamma\undefined\DeclareMathSymbol{\Gamma}{\mathalpha}{letters}{"00}
\let\Omega\undefined\DeclareMathSymbol{\Omega}{\mathalpha}{letters}{"0A}
\let\Lambda\undefined\DeclareMathSymbol{\Lambda}{\mathalpha}{letters}{"03}

\let\to=\longrightarrow
\let\mapsto=\longmapsto
\newcommand{\VecBun}{\catf{VecBun}}

\newcommand{\TVecBun}{\catf{2VecBun}}

\usepackage{geometry}

\definecolor{mygray}{rgb}{0.81,0.81,0.81}

\usepackage{epstopdf}

\begin{document}\thispagestyle{empty}
\setlength{\parskip}{0pt}
\setlength{\abovedisplayskip}{7pt}
\setlength{\belowdisplayskip}{7pt}
\setlength{\abovedisplayshortskip}{0pt}
\setlength{\belowdisplayshortskip}{0pt}

\begin{flushright}
\textsf{ZMP-HH/17-28}\\
\textsf{Hamburger Beiträge zur Mathematik Nr. 706}
\end{flushright}

\vspace*{1cm}

\begin{center}
\Large \textbf{A Parallel Section Functor for 2-Vector Bundles} \normalsize \\[4ex] Christoph Schweigert and Lukas Woike
\end{center}

\begin{center}
\emph{Fachbereich Mathematik,   Universität Hamburg}\\
\emph{Bereich Algebra und Zahlentheorie}\\
\emph{Bundesstra{\ss}e 55,   D -- 20 146  Hamburg}
\end{center}

\begin{abstract}
\noindent \textbf{Abstract.} 
 We associate to a 2-vector bundle over an essentially finite groupoid a 2-vector space
of parallel sections, or, in representation theoretic terms, of higher invariants, which can be described as homotopy fixed points. Our main result is the extension of this assignment to a symmetric monoidal 2-functor $\operatorname{Par} : \mathbf{2VecBunGrpd} \to \mathbf{2Vect}$. It is defined on the symmetric monoidal bicategory $\mathbf{2VecBunGrpd}$ whose morphisms arise from spans of groupoids in such a way that the functor $\operatorname{Par}$ provides pull-push maps between 2-vector spaces of parallel sections of 2-vector bundles.
The direct motivation for our construction comes from the orbifoldization of
extended equivariant topological field theories. 
\end{abstract}

\vspace*{1cm}
\tableofcontents

\section{Introduction and summary}
A representation of a group $G$ on, say, a complex vector space $V$ can be seen as a functor $\star // G \to \Vect$ from the groupoid $\star //G$ with one object $\star$ and automorphism group $G$ to the category $\Vect$ of complex vector spaces sending $\star$ to $V$. 
It is an obvious generalization to replace $\star //G$ by a groupoid $\Gamma$ and call any functor $\varrho: \Gamma \to \Vect$ a representation of $\Gamma$. The limit of the functor $\varrho$ yields the invariants of the representation. 

A functor $\varrho: \Gamma \to \Vect$, i.e.\ a representation of $\Gamma$, is a purely algebraic object. It can also be seen a (flat) vector bundle over the groupoid $\Gamma$. This profitable point of view is for instance emphasized in \cite{willterongerbesgrpds}. It allows us to think of the algebraic notion of invariants of a representation in a geometric way, namely in terms of parallel sections. We will take this convenient geometric point of view throughout our paper. 

If we denote by $\VecBun(\Gamma)$ the category of finite-dimensional vector bundles over a groupoid $\Gamma$, then taking parallel sections yields a functor
\begin{align}
	\Par_\Gamma : \VecBun(\Gamma) \to \FinVect,\label{eqnparovergroupoid}
\end{align} namely the limit functor on the functor category $\VecBun(\Gamma)$. 

There is a higher analogue of a vector bundle over a groupoid, namely a 2-vector bundle over a groupoid, i.e.\ a 2-functor from a given groupoid (seen as a bicategory) to the bicategory $\TwoVect$ of 2-vector spaces, see \cite{baez2group} and \cite{kirrilovg04} for related notions. To a 2-vector bundle $\varrho: \Gamma \to \TwoVect$ over a groupoid $\Gamma$ we associate the category of parallel sections (Definition~\ref{defprallelsectionsof2vectorbundle}) and prove that this category is naturally a 2-vector space if $\Gamma$ is essentially finite (Proposition~\ref{satzactinterparsec}). Hence, as a categorification of \eqref{eqnparovergroupoid} we obtain a 2-functor
\begin{align} \Par_\Gamma  : \TVecBun(\Gamma) \to \TwoVect \label{eqnparovergroupoid2} \end{align} from the bicategory of 2-vector bundles over a fixed groupoid $\Gamma$ to the category of 2-vector spaces. 

The parallel section functors \eqref{eqnparovergroupoid} and \eqref{eqnparovergroupoid2} are mathematically not very challenging and \emph{not} the main concern of this paper. Instead, we are interested in a variant of parallel section functors meeting the requirements determined by our motivation, namely the orbifoldization of equivariant topological field theories: An orbifold construction for (non-extended) equivariant topological field theories can be formulated \cite{schweigertwoikeofk} by means of a parallel section functor
\begin{align}
	\Par : \RepGrpd{} \to \FinVect       \label{eqnfirstparallelsectionfunctor}
\end{align}
for vector bundles over varying groupoids whose construction is given in \cite{trova}. Here, $\RepGrpd{}$ is the symmetric monoidal bicategory from \cite[Section~3.2]{schweigertwoikeofk} whose objects are vector bundles over essentially finite groupoids and whose morphisms come from spans of groupoids and intertwiners. The key point about this functor is that it provides pull-push maps between the vector spaces of parallel sections of vector bundles over different groupoids which are related by a span of groupoids. Having in mind that our parallel section functor is tailored to the application in equivariant topological field theory also explains the importance of spans of groupoids in this paper: Equivariant topological field theories assign quantities to bordisms equipped with principal fiber bundles, and the application of the bundle stack to bordisms, which can be seen as cospans in manifolds, yields exactly spans of groupoids.  Hence, the biased reader may think of all groupoids in this paper as groupoids of principal fiber bundles over some manifold. 

Applications in mathematical physics and representation theory involve an orbifold construction for \emph{extended} equivariant topological field theory (in our terminology \emph{extended} means that the field theory is defined on manifolds up to codimension two).  
In order to generalize the orbifoldization procedure given in \cite{schweigertwoikeofk} to the case of extended field theories in \cite{extofk}, we need a higher analogue of the parallel section functor \eqref{eqnfirstparallelsectionfunctor} or, in other words, the extension of the 2-functor \eqref{eqnparovergroupoid2} to a symmetric monoidal 2-functor defined on a symmetric monoidal bicategory $\TRepGrpd$ of 2-vector bundles over varying groupoids.  
The construction of this symmetric monoidal 2-functor
\begin{align}
	\Par : \TRepGrpd \to \TwoVect,\label{eqntheparallelsectionfunctor}
\end{align} is the main result of this paper (Theorem~\ref{thmextparsecfunctor}). 
The relation between the different parallel section functors is summarized in the diagram
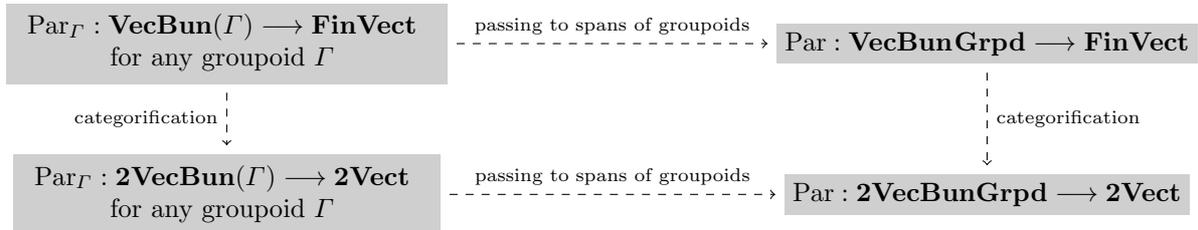
\footnotesize
\begin{center}
	\begin{tikzpicture}[scale=1.9, implies/.style={double,double equal sign distance,-implies},
		dot/.style={shape=circle,fill=black,minimum size=2pt,
			inner sep=0pt,outer sep=2pt},]
		\node (A1) at (0,1) {\colorbox{mygray}{$\begin{array}{c}\Par_\Gamma : \VecBun(\Gamma) \to \FinVect \\ \text{for any groupoid}\ \Gamma\end{array}$ }};
		\node (A2) at (5,1) {\colorbox{mygray}{$\Par : \RepGrpd{} \to \FinVect$}};
		\node (B1) at (0,0) {\colorbox{mygray}{$\begin{array}{c}\Par_\Gamma :  \TVecBun(\Gamma) \to \TwoVect \\ \text{for any groupoid}\ \Gamma\end{array}$ }};
		\node (B2) at (5,0) {\colorbox{mygray}{$\Par : \TRepGrpd \to \TwoVect$}};
		\path[->,font=\scriptsize,dashed]
		(A1) edge node[above]{passing to spans of groupoids} (A2)
		(A1) edge node[left]{categorification} (B1)
		(B1) edge node[above]{passing to spans of groupoids} (B2)
		(A2) edge node[right]{categorification} (B2);
	\end{tikzpicture},
\end{center} \normalsize in which the object of main interest is sitting in the right lower corner. Note that the upper half of the diagram was already discussed in \cite{schweigertwoikeofk}.

Concretely, this paper is organized as follows: In Section~\ref{section2vectorbundlespar} we first recall ordinary vector bundles over groupoids and their parallel sections including the pull-push operations used in \cite{schweigertwoikeofk}. Afterwards, we discuss the higher analogues of these notions, i.e.\ 2-vector bundles over groupoids and their parallel sections and hence the left lower corner of the above diagram. 

Section~\ref{secpullpushext} is devoted to the introduction of pullback and pushforward maps on two different categorical levels needed for the construction of the parallel section functor \eqref{eqntheparallelsectionfunctor}. The discussion of pullback and pushforward 2-morphisms in Section~\ref{pullpush2mor} leads to a higher version of the equivariant Beck-Chevalley condition (Proposition~\ref{mscextversionvonsatz440ausbscwoikeext}), which is of independent interest. 

In Section~\ref{secparallelsectionfunctor} we construct the parallel section functor \eqref{eqntheparallelsectionfunctor}, i.e.\ the right lower corner of the above diagram. To this end, we first have to introduce the domain symmetric monoidal bicategory $\TRepGrpd$ in Section~\ref{secdefdefTRepGrpd}. The objects are 2-vector bundles over essentially finite groupoids, 1-morphisms arise from spans of essentially finite groupoids and intertwiners and 2-morphisms from spans of spans of essentially finite groupoids and higher intertwiners. 
Section~\ref{secparsectionfunctorconstruc} contains the formulation and proof of the main result (Theorem~\ref{thmextparsecfunctor}). Finally, we show how to recover the categorical parallel section functor \eqref{eqnfirstparallelsectionfunctor} from the bicategorical one \eqref{eqntheparallelsectionfunctor} introduced in this paper (Proposition~\ref{satzresofparfunctor}).

\subsection*{Acknowledgements}
We would like to thank Lukas M\" uller for his constant interest in this project and his numerous suggestions on a draft version of this article.
We are also grateful to Ehud Meir and Alexis Virelizier for helpful discussions.

CS is partially supported by the Collaborative Research Centre 676 ``Particles,
Strings and the Early Universe -- the Structure of Matter and Space-Time"
and by the RTG 1670 ``Mathematics inspired by String theory and Quantum
Field Theory''.
LW is supported by the RTG 1670 ``Mathematics inspired by String theory and Quantum
Field Theory''.

\subsection*{Conventions}
All vector spaces or higher analogues thereof encountered in this article will be over the field of complex numbers. Therefore we suppress the field in the notation and write $\Vect$ instead of $\Vect_\mathbb{C}$. Still all constructions would also work over a field of characteristic zero.

\section{2-Vector bundles and their parallel sections\label{section2vectorbundlespar}}
After recalling the notion of a vector bundle over a groupoid, we fix the definition 
of 2-vector bundle used in this text and define the category of parallel sections of a 2-vector bundle. 

\subsection{A brief reminder on vector bundles over groupoids\label{secrecallpullpush}}
Let us review some of the notions and constructions used or introduced in \cite{schweigertwoikeofk}
while implementing also some mild generalizations: A vector bundle over a groupoid $\Omega$ 
with values in a 2-vector space $\cat{V}$ is a functor $\xi: \Omega \to \cat{V}$ 
(in \cite{schweigertwoikeofk} this 2-vector space was always chosen to be the category of vector spaces). 
For the definition of a 2-vector space and the bicategory $\TwoVect$ of 2-vector spaces we refer to \cite{mortonvec} or Example~\ref{extwovectandalg2} below. 

By $\VecBun(\Gamma,\cat{V})$ we denote the category of $\cat{V}$-valued vector bundles over $\Gamma$. 
In case $\Gamma$ is essentially finite, this category naturally carries the structure of a 2-vector space. 

If $\xi: \Omega \to \cat{V}$ is a vector bundle and $\Phi : \Gamma \to \Omega$ a functor between groupoids, 
then we can form the pullback $\Phi^* \xi := \xi \circ \Phi$ of $\xi$ to $\Gamma$. In fact, $\Phi$ gives rise to a pullback functor
\begin{align}
	\Phi^* : \VecBun(\Omega,\cat{V}) \to \VecBun(\Gamma,\cat{V}).\label{eqnpullbackfunctordef}
\end{align} More concisely,
\begin{align}
	\VecBun(?,\cat{V}) : \FinGrpd^\text{\normalfont opp} \to \TwoVect
\end{align} naturally extends to a 2-functor defined on the bicategory of essentially finite groupoids, 
functors and natural transformations. Concretely, it sends a groupoid $\Gamma$ to the 2-vector space $\VecBun(\Gamma,\cat{V})$, 
a functor $\Phi : \Gamma \to \Omega$ to the pullback functor $\Phi^*$ and a natural transformation 
$\eta : \Phi \Rightarrow \Phi'$ of functors $\Phi,\Phi' : \Gamma \to \Omega$ to the obvious 
natural transformation $? (\eta) : \Phi^* \Rightarrow {\Phi'}^*$ whose component 
\begin{align} \xi(\eta) : \Phi^* \xi \to {\Phi'}^* \xi \end{align} 
for $\xi$ in $\VecBun(?,\cat{V})$ consists of the maps $\xi(\eta_x) : \xi(\Phi(x)) \to \xi(\Phi'(x))$ for all $x\in \Gamma$. 

The space $\Par \xi$ of parallel sections of a $\cat{V}$-valued vector bundle $\xi$ over $\Omega$ 
is defined as the limit of the functor $\xi$, 
\begin{align} \Par \xi := \lim \xi,\end{align}
see \cite[Section~3.1]{schweigertwoikeofk}. 
This construction yields a functor
\begin{align}
	\Par_\Omega : \VecBun(\Omega,\cat{V}) \to \cat{V} \end{align} for each essentially finite groupoid $\Omega$. These functors constitute a 1-morphism 
\begin{align}
	\Par : \VecBun(?,\cat{V}) \to \cat{V}
\end{align} in the bicategory of 2-functors $\FinGrpd^\text{\normalfont opp} \to \TwoVect$, where we use $\cat{V}$ to denote the constant 2-functor with value $\cat{V}$. 

By the following standard fact limits can be pulled back:

\begin{lemma}\label{msclemmalimesbildungaeqkat}
Let $\cat{C}$ be a complete category and $X: \cat{J} \to \cat{C}$ a functor from a small category $\cat{J}$ to $\cat{C}$.
Then any functor $\Phi : \cat{I} \to \cat{J}$ of small categories induces a morphism
\begin{align}   \lim X \to \lim \Phi^* X . \end{align} If $\Phi$ is an equivalence, then this morphism is an isomorphism. The dual statement for colimits is true if $\cat{C}$ is cocomplete.
\end{lemma}

For a functor $\Phi : \Gamma \to \Omega$ between groupoids we obtain a natural map
\begin{align} \Phi^* : \Par \xi \to \Par \Phi^* \xi,\end{align} the \emph{pullback map}. By abuse of notation it is denoted by the same symbol as the pullback functor \eqref{eqnpullbackfunctordef}, but should not be confused with the latter.

In case that $\Gamma$ and $\Omega$ are essentially finite, we introduced in \cite[Section~3.4]{schweigertwoikeofk} also a \emph{pushforward map}
\begin{align} \Phi_* : \Par \Phi^* \xi \to \Par  \xi\end{align}
by integration over the homotopy fiber $\Phi^{-1}[y]$ over $y\in \Omega$. Recall that for a given $y\in \Omega$, an object $(x,g) \in \Phi^{-1}[y]$ in the homotopy fiber over $y$ is an object $x\in \Gamma$ together with a morphism $g: \Phi(x) \to y$ in $\Gamma$.
Since $\Par \Phi^* \xi$ is the limit of $\Phi^* \xi$, it comes equipped with maps $\pi_x : \Par \Phi^* \xi \to \xi(\Phi(x))$ for each $x\in \Gamma$, which we can use to form the concatenation
\begin{align}
	\nu_{x,g} : \Par \Phi^* \xi \stackrel{\pi_x}{\to} \xi(\Phi(x)) \xrightarrow{\xi(g)} \xi(y).
\end{align} Recall that a morphism $(x,g) \to (x',g')$ in the groupoid $\Phi^{-1}[y]$ is a morphism $h: x \to x'$ such that $g' \Phi(h) = g$. Now an easy computation shows that the morphism $\nu_{x,g}$ only depends on the isomorphism class of $(x,g)$ in $\Phi^{-1}[y]$. This allows us to define
\begin{align} \label{eqnoverhomotopyfibers} \int_{\Phi^{-1}[y]} \nu_{x,g}\,\td (x,g) := \sum_{[x,g] \in \pi_0(\Phi^{-1}(y))} \frac{\nu_{x,g}}{|\!  \Aut(x,g)|} : \Par \Phi^* \xi \to \xi(y).  \end{align} 
The morphisms $\nu_{x,g}$ can be added and multiplied by scalars since $\Hom_\cat{V}(\Phi^* \xi,\xi(y))$ is a complex vector space.
Formula~\eqref{eqnoverhomotopyfibers} provides us with an instance of an \emph{integral with respect to groupoid cardinality}, i.e.\ a sum over the values of an invariant function on an essentially finite groupoid, here $\Phi^{-1}[y]\ni (x,g) \mapsto \nu_{x,g}$, taking values in a complex vector space, here $\Hom_\cat{V}(\Phi^* \xi,\xi(y))$, weighted by the cardinalities of the automorphism groups in our groupoid. For more background on groupoid cardinality we refer to \cite{baezgroup}. The integral with respect to groupoid cardinality was also an essential concept for the construction of the parallel section functor in \cite{schweigertwoikeofk} and is also fully recalled there. 

An easy computation shows that for any morphism $a: y \to y'$
\begin{align} \xi(a) \int_{\Phi^{-1}[y]} \nu_{x,g}\,\td (x,g) = \int_{\Phi^{-1}[y']} \nu_{x',g'}\,\td (x',g').\end{align} This implies that the maps \eqref{eqnoverhomotopyfibers} induce a natural map
\begin{align} \Phi_* :\Par  \Phi^* \xi \to \Par \xi,\end{align} the so-called \emph{pushforward map}.

The most important properties of pullback and the pushforward map include the composition laws and the equivariant Beck-Chevalley property. The proofs are obvious generalizations of those in \cite{schweigertwoikeofk}.

\begin{proposition}\label{satzcompositionslawspullpush}
Let $\cat{V}$ be a 2-vector space $\Phi : \Gamma \to \Omega$ and $\Psi : \Omega \to \Lambda$ be functors between essentially finite groupoids.
\begin{enumerate}
	\item For any $\cat{V}$-valued vector bundle $\xi$ over $\Lambda$ the composition law $(\Psi \circ \Phi)^* = \Phi^* \circ \Psi^*$ for the pullback maps holds.\label{satzcompositionslawspullpusha}
	\item For any $\cat{V}$-valued vector bundle $\xi$ over $\Gamma$ the composition law $(\Psi \circ \Phi)_* = \Psi_*\circ \Phi_*$ for the pushforward maps holds.\label{satzcompositionslawspullpushb}
\end{enumerate}
\end{proposition}

\begin{proposition}[Equivariant Beck-Chevalley condition]\label{satzparsecbeckchevalley}
For the homotopy pullback 
\begin{center}
	\begin{tikzpicture}[scale=2, implies/.style={double,double equal sign distance,-implies},
		dot/.style={shape=circle,fill=black,minimum size=2pt,
			inner sep=0pt,outer sep=2pt},]
		\node (A1) at (0,1) {$\Gamma\times_\Omega \Lambda$};
		\node (A2) at (1,1) {$\Gamma$};
		\node (B1) at (0,0) {$\Lambda$};
		\node (B2) at (1,0) {$\Omega$};
		\path[->,font=\scriptsize]
		(A1) edge node[above]{$\pi_\Gamma$} (A2)
		(A1) edge node[left]{$\pi_\Lambda$} (B1)
		(A2) edge node[right]{$\Phi$} (B2)
		(B1) edge node[below]{$\Psi$} (B2);
		\draw (A2) edge[implies] node[above] {\scriptsize$\eta\ $} (B1);
	\end{tikzpicture}
\end{center}
of a cospan $\Lambda \stackrel{\Psi}{\to} \Omega \stackrel{\Phi}{\longleftarrow} \Gamma$ of essentially finite groupoids and any $\cat{V}$-valued vector bundle $\xi$ over $\Omega$ the pentagon relating different pull-push combinations
\begin{center}
	\begin{tikzpicture}[scale=3]
		\node (A1) at (0,2) {$\Par \Phi^*\xi$};
		\node (A2) at (2,2) {$\Par \xi$};
		\node (C) at (0,1) {$\Par \pi_\Gamma^* \Phi^* \xi$};
		\node (B1) at (1,0) {$\Par \pi_\Lambda^* \Psi^* \xi$};
		\node (B2) at (2,0) {$\Par \Psi^* \xi$};
		\path[->,font=\scriptsize]
		(A1) edge node[above]{$\Phi_*$} (A2)
		(A1) edge node[left]{$\pi_\Gamma^*$} (C)
		(C) edge node[left]{$\xi(\eta)_*$} (B1)
		(B1) edge node[above]{${\pi_\Lambda}_*$} (B2)
		(A2) edge node[right]{$\Psi^*$} (B2);
	\end{tikzpicture}
\end{center} commutes. 
\end{proposition}

\subsection{2-Vector bundles}
The goal of this article is a bicategorical generalization of the parallel section functor. Hence, we need a higher analogue of a groupoid representation (or a vector bundle over a groupoid). This will be the notion of a 2-vector bundle below. For related notions and generalizations see \cite{baez2group} and \cite{kirrilovg04}. The definition requires the notion of a bicategory and 2-functors, see \cite{leinster}, and of a symmetric monoidal bicategory, see \cite{schommerpries}. 

\begin{definition}[2-Vector bundle]\label{defsymmonbicatof2vectorbundles}
	A \emph{2-vector bundle $\varrho$ over a groupoid $\Gamma$ with values in a symmetric monoidal bicategory $\cat{S}$} is a representation of $\Gamma$ on $\cat{S}$, i.e.\ a 2-functor $\varrho : \Gamma \to \cat{S}$, where $\Gamma$ is seen as a bicategory without non-trivial 2-morphisms. (There are no monoidality requirements on $\varrho$.) By $\TVecBun(\Gamma,\cat{S})$ we denote the symmetric monoidal bicategory of $\cat{S}$-valued 2-vector bundles over $\Gamma$. 
\end{definition}

\remark\label{bmkdescr2rep}
\begin{enumerate}
	
	\item We use the term `2-vector bundle' although we do not assume any (higher) linear structure on the target $\cat{S}$. 
	
	\item 	Let us partly unpack the definition of a 2-vector bundle $\varrho : \Gamma \to \cat{S}$:
	\begin{itemize}
		\item To $x\in \Gamma$ the 2-vector bundle $\varrho$ assigns an object $\varrho(x)$ in $\cat{S}$, also called the \emph{fiber of $\varrho$ over $x$}.
		\item To a morphism $g: x \to y$ in $\Gamma$ the 2-vector bundle assigns a 1-morphism $\varrho(g) : \varrho(x) \to \varrho(y)$, which is geometric terms can be thought of as a parallel transport operator.
		\item The data comprises natural isomorphisms
		\begin{align} \eta_x : \varrho(\id_x) &\cong \id_{\varrho(x)},\\ \alpha_{gh} : \varrho(g) \circ \varrho(h) &\cong \varrho(gh) \end{align} for composable morphisms in $\Gamma$. These natural isomorphisms are subject to obvious coherence conditions. 
	\end{itemize} \label{bmkdescr2rep1}
	
	\item Let us describe the bicategory $\TVecBun(\Gamma,\cat{S})$ in more detail: \begin{enumerate}
		\item[(0)] Objects are 2-vector bundles over $\Gamma$.
		\item[(1)] 1-morphisms are 2-vector bundles morphisms or, equivalently, intertwiners. An intertwiner $\phi : \varrho \to \xi$ of 2-vector bundles over $\Gamma$ consists of 1-morphisms $\phi_x : \varrho(x) \to \xi(x)$ for each $x\in \Gamma$ and natural morphisms 
		\begin{align}
			\xi(g) \circ \phi_x  \stackrel{\theta_g} {\to}  \phi_y \circ \varrho(g) \myforall g: x \to y 
		\end{align} subject to obvious coherence conditions. These coherence conditions entail in particular that all $\theta_g$ are 2-isomorphisms. For this it is crucial that $\Gamma$ is a groupoid. 
		\item[(2)] A 2-morphism $\eta : \phi \to \psi$ between 1-morphisms $(\phi,\theta)$ and $(\psi,\kappa)$ between the 2-vector bundles $\varrho$ and $\xi$ consists of 2-morphisms $\eta_x : \phi_x \to \psi_x$ such that for all $g: x \to y$ the square
		\begin{center}
			\begin{tikzpicture}[scale=2.5]
				\node (A1) at (0,1) {$\phi_y \circ \varrho(g)$};
				\node (A2) at (1,1) {$\xi(g)\circ \phi_x$};
				\node (B1) at (0,0) {$\psi_y \circ \varrho(g)$};
				\node (B2) at (1,0) {$\xi(g) \circ \psi_x$};
				\path[->,font=\scriptsize]
				(A1) edge node[left]{$\eta_y$} (B1);
				\path[->,font=\scriptsize]
				(A1) edge node[above]{$\theta_g$} (A2)
				(B1) edge node[above]{$\kappa_g$} (B2)
				(A2) edge node[right]{$\eta_x$} (B2);
			\end{tikzpicture}
		\end{center} commutes. Here $\eta_x : \xi(g) \circ \phi_x \to \xi(g) \circ \psi_x$ is the 2-morphism induced by $\eta_x$ and the identity on $\varrho(x)$, but we suppress the identity morphism in the notation for readability. 
		
	\end{enumerate}
	The tensor product in $\TVecBun(\Gamma,\cat{S})$ is the tensor product in $\cat{S}$ applied object-wise. The monoidal unit $\mathbb{I}_\Gamma$ in $\TVecBun(\Gamma,\cat{S})$ assigns to each $x\in \Gamma$ is the monoidal unit $\mathbb{I}$ in $\cat{S}$ and to all morphisms in $\Gamma$ the identity 1-morphism. 
	\label{bmkdescr2rep2}
	
\end{enumerate}
\endremark

\example[The symmetric monoidal bicategory $\TwoVect$]\label{extwovectandalg2}
Let us review the main example of a symmetric monoidal bicategories relevant in this text, namely the symmetric monoidal bicategory $\TwoVect$ of \emph{2-vector spaces (of Kapranov-Voevodsky type\footnote{There are other types of 2-vector spaces, but throughout this text we will always mean 2-vector spaces of Kapranov-Voevodsky type when talking about 2-vector spaces. In particular, we always work over the complex field when talking about 2-vector spaces.})}, see \cite{mortonvec}:
\begin{enumerate}
	\item[(0)] Objects are 2-vector spaces, i.e.\ $\mathbb{C}$-linear additive semisimple categories with biproducts and finitely many simple objects up to isomorphism. 
	
	\item[(1)] 1-Morphisms are $\mathbb{C}$-linear functors, which are also called \emph{2-linear maps}. 
	
	\item[(2)] 2-Morphisms are natural transformations of $\mathbb{C}$-linear functors.
	
\end{enumerate} 
The tensor product is the Deligne product, the tensor unit is the category $\FinVect$ of finite-dimensional complex vector spaces. 
For any 2-vector space $\cat{V}$ we can choose a \emph{basis}, i.e.\ a family of representatives for the finitely many isomorphism classes of simple objects. Having chosen a basis $\mathscr{B}$ of a 2-vector space $\cat{V}$ we can write any object $X$ in $\cat{V}$ as a biproduct
\begin{align}
	X \cong \bigoplus_{B\in \mathscr{B}} V_B * B,
\end{align} where the $V_B$ are finite-dimensional complex vector spaces and where $V_B*B$ is the $\dim V_B$-fold biproduct of $B$ with itself.
The $\mathbb{C}$-linearity of a functor between 2-vector spaces is equivalent to the preservation of biproducts.
Consequently, any 2-linear map $\cat{V} \to \cat{W}$ is determined by its values on this basis, which allows us the describe 2-linear maps in terms of matrices with vector space valued entries. 
Moreover note that 2-linear maps $\cat{V} \to \cat{W}$ are precisely the exact functors. Indeed, exactness implies preservations of biproducts. The converse holds since all short exact sequences in $\cat{V}$ split by semisimplicity. 

Up to $\mathbb{C}$-linear equivalences a 2-vector spaces is determined by the cardinality of its basis, which we will also refer to as \emph{dimension}. For instance, any $n$-dimensional 2-vector space is equivalent to the category $\mathbb{C}[\mathbb{Z}_n]\text{-}\Mod$ of finite-dimensional complex modules over the group algebra of the cyclic group $\mathbb{Z}_n$.

Since the symmetric monoidal bicategory $\TwoVect$ will be the most important one for us in the sequel, we agree on the notation $\TVecBun(\Gamma):= \TVecBun(\Gamma,\TwoVect)$, i.e.\ 2-vector bundles with unspecified target category always have to be understood as $\TwoVect$-valued 2-vector bundles. 

For later use we recall that the symmetric monoidal category obtained by restriction of $\TwoVect$ to the endomorphisms of the monoidal unit is the category $\FinVect$ of finite-dimensional complex vector spaces. 

\endexample

\subsection{Parallel sections of 2-vector bundles\label{parsecsecfor2vecbund}}
Parallel sections of a given vector bundle $\varrho$ (or, equivalently, invariants of the representation $\varrho$) can be obtained by taking the morphisms of the trivial vector bundle to $\varrho$. This principle can be directly generalized to 2-vector bundles.

\begin{definition}[Parallel sections of a 2-vector bundle]\label{defprallelsectionsof2vectorbundle}
	Let $\cat{S}$ be a symmetric monoidal bicategory. 
	The \emph{category of parallel sections} of an $\cat{S}$-valued 2-vector bundle $\varrho$ over a groupoid $\Gamma$ is the category
	\begin{align} \Par \varrho := \Hom_{\TVecBun(\Gamma,\cat{S})}(\mathbb{I}_\Gamma,\varrho).\end{align}
\end{definition}

\remark
\begin{enumerate}
	\item A parallel section $s \in \Par \varrho$ gives us a 1-morphism $s(x): \mathbb{I} \to \varrho(x)$ in $\cat{S}$ for each $x\in \Gamma$ and coherent isomorphisms $s(y) \cong \varrho(g) \circ s(x)$ for all $g: x \to y$ in $\Gamma$. Instead of $\varrho(x) \circ s(x)$ we will often write $\varrho(x) s(x)$ or even $g.s(x)$ if the vector bundle is clear from the context.  For $\cat{S}=\TwoVect$ the monoidal unit $\mathbb{I}$ is given by the category $\FinVect$ of finite-dimensional complex vector spaces. Note that 1-morphisms $\FinVect \to \varrho(x)$ can be identified with the value on $\mathbb{C}$ and hence with an object in the fiber $\varrho(x)$. 
	
	\item The parallel sections of a 2-vector bundle are thus `parallel up to isomorphism', where the isomorphism is part of the data. Hence, being parallel is no longer a property, but structure. In other contexts the parallel sections considered here would be called homotopy fixed points, see e.g. \cite{hsvgroupbicat}.
	
\end{enumerate}\endremark

\noindent For a $\TwoVect$-valued 2-vector bundle we would like to find conditions under which the category of parallel sections is naturally a 2-vector space again. 
It is easy to see that the category of parallel sections inherits all the needed structure and properties from the 2-vector bundle except for finite semisimplicity. 
In order to look at this last missing point more closely, we use techniques and results from \cite{kirrilovorbII}.

First of all we note that it suffices to study 2-vector bundles over connected groupoids, i.e.\ we can concentrate on 2-vector bundles $\varrho: \star //G \to \TwoVect$ for a finite group $G$. In this case we obtain a 2-vector space $\cat{V}:= \varrho(\star)$, and any $g\in G$ yields a 2-linear equivalence $\varrho(g) : \cat{V} \to \cat{V}$. These equivalences fulfill the properties of a representation only up to isomorphism as discussed in Remark~\ref{bmkdescr2rep}, \ref{bmkdescr2rep1}, but still we will refer to this as a representation of $G$ on $\cat{V}$. We denote the action of $g\in G$ on some object $X\in \cat{V}$ by $g. X$ and the evaluation of the coherence isomorphisms on $X$ by
\begin{align} \beta_{g,h}^X : g.h.X \to (gh).X.\end{align}
According to Definition~\ref{defprallelsectionsof2vectorbundle}, $\Par \varrho$ is the category of pairs $(X,\phi=(\phi_g)_{g\in G})$, where $X$ is in $\cat{V}$ and $\phi$ is a family of coherent isomorphisms $\phi_g : g.X := \varrho(g)(X) \to X$. 

Let $(X_s)_{s\in\mathscr{S}}$ be a basis of $\cat{V}$. 
Since any $g\in G$ acts as an equivalence, it maps simple objects to simple objects. Hence, when forgetting about the coherence data, $g\in G$ just acts as a permutation of the basis. Consequently, we obtain an action of $G$ on $\mathscr{S}$. We denote the corresponding action groupoid by $\mathscr{S}//G$ and the set of orbits by $\mathscr{S}/G$. 

For a given orbit $\mathcal{O} \in \mathscr{S}/G$ and $s\in \mathcal{O}$ 
there is an isomorphism $\xi_g^s : g. X_s \to X_{g.s}$. It is unique up to multiplication by an element in $\mathbb{C}^\times=\mathbb{C}\setminus \{0\}$ since $X_{g.s}$ is simple, and we fix such an isomorphism. Now for $g,h \in G$ we find $\xi_g^{h.s} \circ g. \xi_h^s \circ \left(\beta_{gh}^{X_s}\right)^{-1} = \alpha_{gh}^s \xi^s_{gh}$ for some $\alpha_{gh}^s \in \mathbb{C}^\times$ since $X_{(gh).s}=X_{g.h.s}$ is simple. 
The scalars $\alpha_{gh}^s$ form a cocycle $\alpha_\mathcal{O}\in Z^2(G;\Map(\mathcal{O},\mathbb{C}^\times))$ with coefficients in the Abelian group of functions $\mathcal{O} \to \mathbb{C}^\times$, and the class $[\alpha_\mathcal{O}] \in H^2(G;\Map(\mathcal{O},\mathbb{C}^\times))$ does not depend on the chosen isomorphisms $\xi_g^s$. This cocycle (or its class) is used to define the \emph{twisted group algebra} $\mathfrak{A}_{\alpha_\mathcal{O}}(G,\mathcal{O})$, a semisimple finite-dimensional complex algebra. 
Using the category of modules over these twisted group algebras we can state a version of the following result of \cite{kirrilovorbII}:

\begin{proposition}[Kir01, Theorem~3.5]  \label{satzthm35kir}
Let $\varrho: \star //G \to \TwoVect$ be a 2-vector bundle, i.e.\ a representation of the group $G$ on a 2-vector space $\cat{V}:= \varrho(\star)$. 
There is an equivalence 
\begin{align}
	\Par \varrho \cong	\bigoplus_{\mathcal{O} \in \mathscr{S}/G} \mathfrak{A}_{\alpha_\mathcal{O}}(G,\mathcal{O})\text{-}\Mod
\end{align} of Abelian categories. Hence if $G$ is finite, $\Par \varrho$ is semisimple with finitely many simple objects and hence a 2-vector space. 
\end{proposition}

Since 
\begin{align}
	\Map (\mathscr{S},\mathbb{C}^\times) \cong \Map \left( \coprod_{\mathcal{O} \in \mathscr{S} / G} \mathcal{O},\mathbb{C}^\times  \right) \cong \prod_{\mathcal{O} \in \mathscr{S} / G} \Map (\mathcal{O},\mathbb{C}^\times)
\end{align} we can combine the cocycles $\alpha_\mathcal{O}$ with coefficients in $\Map (\mathcal{O},\mathbb{C}^\times)$ into a cocycle $\alpha$ with coefficients in $\Map (\mathscr{S},\mathbb{C}^\times)$. In the sequel we will rather use the cocycle $\alpha$ and write $\mathfrak{A}_{\alpha}(G,\mathcal{O})\text{-}\Mod$ instead of $\mathfrak{A}_{\alpha_\mathcal{O}}(G,\mathcal{O})\text{-}\Mod$ (the dependence on the orbit is still present in the notation, so there is no risk of confusion).


Summarizing and extending to non-connected groupoids, we conclude that for an essentially finite groupoid $\Gamma$ taking parallel sections of a $\TwoVect$-valued 2-vector bundle $\varrho$ over $\Gamma$ produces a 2-vector space $\Par \varrho$, which is entirely determined  by the action groupoid $\mathscr{S} //\Gamma$ and a gerbe on $\Gamma$, i.e.\ a class $H^2 (\Gamma; \Map(\mathscr{S},\mathbb{C}^\times))$ with coefficients in the Abelian group of functions $\mathscr{S} \to \mathbb{C}^\times$. Obviously, taking parallel sections is also 2-functorial:

\begin{proposition}\label{satzactinterparsec}
Taking parallel sections of 2-vector bundles over an essentially finite groupoid $\Gamma$ naturally extends to a 2-functor
\begin{align}
	\Par_\Gamma: \TVecBun(\Gamma) \to \TwoVect.
\end{align} The image of a 1-morphism $\lambda$ or a 2-morphism $\eta$ will be denoted by $\lambda_*$ or $\eta_*$, respectively. 
\end{proposition}

\noindent We should emphasize that this is just the parallel section functor for 2-vector bundles over \emph{one fixed groupoid}. It is \emph{not} the parallel section functor we intend to construct in this article as our main result.

\section{Pullback and pushforward\label{secpullpushext}}
The construction of the parallel section functor in Section~\ref{secparallelsectionfunctor} relies on certain pullback and pushforward maps that we will introduce in this section on two different categorical levels.

\subsection{Pullback and pushforward 1-morphisms}
Just like ordinary bundles, 2-vector bundles have an obvious notion of pullback: For any functor $\Phi : \Gamma \to \Omega$ between groupoids we obtain a pullback functor 
\begin{align}
	\Phi ^* : \TVecBun(\Omega,\cat{S}) \to \TVecBun(\Gamma,\cat{S}) 
\end{align} by precomposition. 

Additionally, for any 2-vector bundle $\varrho$ over $\Omega$ we get a pullback 1-morphism in $\cat{S}$
\begin{align}
	\Phi^* : \Par \varrho \to \Par \Phi^* \varrho, \quad s \mapsto \Phi^* s
\end{align} denoted by the same symbol and given by
\begin{align}
	(\Phi^* s)(x) := s(\Phi(x)) \myforall x\in \Gamma
\end{align} together with the isomorphisms
\begin{align}
	(\Phi^* s)(y)=s(\Phi(y))\cong\Phi(g). s(\Phi(x)) = (\Phi^* \varrho)(g) (\Phi^* s)(x)  \myforall g: x \to y \quad \text{in $\Gamma$}.
\end{align}
It is now easy to prove the following statements: 

\begin{proposition}[Pullback 1-morphism]\label{satzpulloperationsofk}
Let $\Phi : \Gamma \to \Omega$ be a functor between groupoids.
\begin{enumerate}

	\item Contravariance: If $\Psi: \Omega \to \Lambda$ is another functor, then we have $(\Psi \circ \Phi)^* =\Phi^* \circ \Psi^*$ for both the functors as well as the pullback 1-morphisms induced by them.\label{satzpulloperationsofkb}
	
	\item If $\eta : \Phi \to \Phi'$ is natural transformation of functors $\Gamma \to \Omega$, then for any 2-vector bundle $\varrho$ over $\Omega$ we obtain a 1-isomorphism $\varrho(\eta): \Phi^*\varrho \to {\Phi'}^*\varrho$. \label{satzpulloperationsofkc}
	
	\item Naturality: The pull maps are natural in the sense that for any 1-morphism $\lambda : \varrho \to \xi$ of 2-vector bundles over $\Omega$ the square 
	\begin{center}
		\begin{tikzpicture}[scale=2, implies/.style={double,double equal sign distance,-implies},
			dot/.style={shape=circle,fill=black,minimum size=2pt,
				inner sep=0pt,outer sep=2pt},]
			\node (A1) at (0,1) {$\Par \varrho$};
			\node (A2) at (1,1) {$\Par \Phi^* \varrho$};
			\node (B1) at (0,0) {$\Par \xi$};
			\node (B2) at (1,0) {$\Par \Phi^* \xi$};
			\path[->,font=\scriptsize]
			(A1) edge node[above]{$\Phi^*$} (A2)
			(A1) edge node[left]{$\lambda_*$} (B1)
			(B1) edge node[above]{$\Phi^*$} (B2)
			(A2) edge node[right]{$(\Phi^* \lambda)_*$} (B2);
		\end{tikzpicture}
	\end{center}
	commutes strictly.\label{satzpulloperationsofkd}
	The vertical arrows are the images of 2-vector bundle morphisms under the functor from Proposition~\ref{satzactinterparsec}.
	
\end{enumerate}

\end{proposition}

\noindent We define the pushforward 1-morphisms via the limit of diagrams with shape of a homotopy fiber and values in spaces of 1-morphisms:

\begin{definition}[Pushforward 1-morphism]\label{satztwolinpushmap}
	Let $\cat{S}$ be a symmetric monoidal bicategory with complete categories of 1-morphisms between any two objects. 
	Let $\Phi : \Gamma \to \Omega$ be a functor between groupoids. For an $\cat{S}$-valued 2-vector bundle $\varrho$ over $\Omega$ and $s\in \Par \Phi^* \varrho$ we define the parallel section $\Phi_* s \in \Par  \varrho$ by
	\begin{align}
		(\Phi_* s)(y) := \lim_{(x,g)\in \Phi^{-1}[y]} g.s(x)  \myforall y\in \Omega .
	\end{align} The limit is taken in $\Hom_\cat{S}(\mathbb{I},\varrho(y))$ and has the shape of the homotopy fiber $\Phi^{-1}[y]$ of $\Phi$ over $y\in\Omega$. 
	We call the resulting 2-linear map
	\begin{align} \Phi_* : \Par \Phi^* \varrho \to \Par \varrho, \quad s \mapsto \Phi_* s \end{align} \emph{pushforward 1-morphism}. 
\end{definition}
\remark
\begin{enumerate} \item An easy computation shows that $\Phi_*$ actually takes values in parallel sections.

	\item In the special case $\cat{S}=\TwoVect$ the functor
	\begin{align} \Phi^{-1}[y] \to \Hom_\TwoVect(\FinVect,\varrho(y))\cong\varrho(y), \quad (x,g) \mapsto g.s(x) \end{align} is a $\varrho(y)$-valued vector bundle over the homotopy fiber $\Phi^{-1}[y]$ of $\Phi$ over $y$. Its limit, by definition, coincides with $(\Phi_* s)(y)$ and is the space of parallel sections of this vector bundle as recalled in Section~\ref{secrecallpullpush}.  
	
\end{enumerate}
\endremark

\noindent We now recall an observation already needed in \cite[Lemma~4.12]{schweigertwoikeofk}:

\begin{lemma}\label{lemmahtppullbackcomposition}
For composable functors $\Phi : \Gamma \to \Omega$ and $\Psi : \Omega \to \Lambda$ between groupoids, there is a canonical equivalence
\begin{align}
	(\Psi \circ \Phi)^{-1}[z] \cong \Psi^{-1}[z] \times_\Omega \Gamma.
\end{align}
\end{lemma}

\begin{proposition}\label{satzmsceigpushmapsext}
Let $\Phi : \Gamma \to \Omega$ and $\Psi : \Omega \to \Lambda$ between be composable functors between groupoids.
\begin{enumerate}
	\item Covariance: The pushforward 1-morphisms obey the composition law $(\Psi\circ \Phi)_* \cong \Psi_* \circ \Phi_*$ by a canonical isomorphism.\label{satzmsceigpushmapsexta}
	
	\item Naturality: The pushforward 1-morphisms are natural in the sense that for any 1-morphism $\lambda : \varrho \to \xi$ of 2-vector bundles over $\Omega$ the square 
	\begin{center}
		\begin{tikzpicture}[scale=2, implies/.style={double,double equal sign distance,-implies},
			dot/.style={shape=circle,fill=black,minimum size=2pt,
				inner sep=0pt,outer sep=2pt},]
			\node (A1) at (0,1) {$\Par \Phi^*\varrho$};
			\node (A2) at (1,1) {$\Par \varrho$};
			\node (B1) at (0,0) {$\Par \Phi^* \xi$};
			\node (B2) at (1,0) {$\Par \xi$};
			\path[->,font=\scriptsize]
			(A1) edge node[above]{$\Phi_*$} (A2)
			(A1) edge node[left]{$(\Phi^* \lambda)_*$} (B1)
			(B1) edge node[above]{$\Phi_*$} (B2)
			(A2) edge node[right]{$\lambda_*$} (B2);
			\draw (A2) edge[implies] node[above] {\scriptsize$\cong\ $} (B1);
		\end{tikzpicture}
	\end{center}
	commutes up to a canonical natural isomorphism arising from the coherence isomorphism that $\lambda$ comes equipped with.\label{satzmsceigpushmapsextb}

	\item The naturality isomorphisms and the composition of pushforward 1-morphisms are compatible in the sense that for a 1-morphism $\lambda : \varrho \to \xi$ of 2-vector bundles over $\Lambda$ we have the equality of 2-isomorphisms \label{satzmsceigpushmapsextc}
	\small
	
	\begin{align}\begin{array}{c}
			\begin{tikzpicture}[scale=1.5,     implies/.style={double,double equal sign distance,-implies},
				dot/.style={shape=circle,fill=black,minimum size=2pt,
					inner sep=0pt,outer sep=2pt},]
				\node (A1) at (0,1) {$\Par \Phi^* \Psi^* \varrho$};
				\node (A2) at (2,1) {$ \Par \Psi^*\varrho$};
				\node (A3) at (4,1) {$\Par \varrho$};
				\node (B1) at (0,0) {$\Par \Phi^* \Psi^* \xi$};
				\node (B2) at (2,0) {$\Par \Psi^* \xi$};
				\node (B3) at (4,0) {$\Par \xi$};
				\node (E) at (4.7,0.5) {$=$};
				\node (C1) at (6,1) {$\Par \Phi^* \Psi^* \varrho$};
				\node (C2) at (8,1) {$ \Par \varrho$};
				\node (D1) at (6,0) {$\Par \Phi^* \Psi^* \xi$};
				\node (D2) at (8,0) {$ \Par \xi$};
				\path[->,font=\scriptsize]
				(A1) edge node[left]{$(\Phi^* \Psi^* \lambda)_*$} (B1)
				(C1) edge node[above]{$(\Psi\circ\Phi)_*$} (C2)
				(D1) edge node[above]{$(\Psi\circ\Phi)_*$} (D2)
				(A2) edge node[left]{$(\Psi^* \lambda)_*$} (B2)
				(C1) edge node[left]{$(\Phi^* \Psi^* \lambda)_*$} (D1)
				(A3) edge node[left]{$\lambda_*$} (B3)
				(C2) edge node[left]{$\lambda_*$} (D2)
				(A1) edge node[above]{$\Phi_*$} (A2)
				(A2) edge node[above]{$\Psi_*$} (A3)
				(B1) edge node[above]{$\Phi_*$} (B2)
				(B2) edge node[above]{$\Psi_*$} (B3);
				\draw (A2) edge[implies] node[above] {\scriptsize$\cong $} (B1);
				\draw (A3) edge[implies] node[above] {\scriptsize$\cong $} (B2);
				\draw (C2) edge[implies] node[above] {\scriptsize$\cong $} (D1);
				\draw (3,2) edge[implies] node[above] {\scriptsize$\cong $} (1,1.5);
				\draw (3,-0.5) edge[implies] node[above] {\scriptsize$\cong $} (1,-1);
				\draw[line width=0.5pt,->]
				(A1)   \dir{90}{90} node[above] {\scriptsize$(\Psi\circ \Phi)_*$} (A3) 
				;
				\draw[line width=0.5pt,->]
				(B1)   \dir{270}{270} node[below] {\scriptsize$(\Psi\circ \Phi)_*$} (B3) 
				;
			\end{tikzpicture} .
	\end{array}\end{align}
	\normalsize
\end{enumerate}
\end{proposition}

\begin{proof}
For a parallel section $s$ of a 2-vector bundle $\varrho$ over $\Lambda$ and $z\in \Lambda$ we find by definition
\begin{align}
	(\Psi_* \Phi_* s)(z) &\cong \lim_{(y,h) \in \Psi^{-1}[z]} \lim_{ (x,g) \in \Phi^{-1}[y]} (h\Psi(g)).s(x)\label{eqndoublelimit},\\
	((\Psi \circ \Phi)_*s)(z) &= \lim_{(x,k) \in (\Psi \circ \Phi)^{-1}[z]} k.s(x).
\end{align} The double limit \eqref{eqndoublelimit} can be seen as a limit over the homotopy pullback $\Psi^{-1}[z] \times_\Omega \Gamma$, which by Lemma~\ref{lemmahtppullbackcomposition} is canonically equivalent to $(\Psi \circ \Phi)^{-1}[z] \cong \Psi^{-1}[z] \times_\Omega \Gamma$. This equivalence yields the needed isomorphism $(\Psi\circ \Phi)_* \cong \Psi_* \circ \Phi_*$ by Lemma~\ref{msclemmalimesbildungaeqkat}. The remaining assertions can be directly verified.    
\end{proof}

\subsection{Pullback and pushforward 2-morphisms and the equivariant Beck-Chevalley condition\label{pullpush2mor}}
So far, we have established pullback and pushforward 1-mor\-phisms. In the next step, we will provide pull and push 2-morphisms between different pull-push combinations.

We consider a weakly commuting square 
\begin{center}
	\begin{tikzpicture}[scale=2, implies/.style={double,double equal sign distance,-implies},
		dot/.style={shape=circle,fill=black,minimum size=2pt,
			inner sep=0pt,outer sep=2pt},]
		\node (A1) at (0,1) {$\Pi$};
		\node (A2) at (1,1) {$\Gamma$};
		\node (B1) at (0,0) {$\Lambda$};
		\node (B2) at (1,0) {$\Omega$};
		\path[->,font=\scriptsize]
		(A1) edge node[above]{$P$} (A2)
		(A1) edge node[left]{$Q$} (B1)
		(A2) edge node[right]{$\Phi$} (B2)
		(B1) edge node[below]{$\Psi$} (B2);
		\draw (A2) edge[implies] node[above] {\scriptsize$\eta\ $} (B1);
	\end{tikzpicture}
\end{center} of essentially finite groupoids. 
By the definition of the homotopy fiber $Q^{-1}[y]$ we obtain a weakly commutative square
\begin{center} \begin{tikzpicture}[scale=1.8, implies/.style={double,double equal sign distance,-implies},
		dot/.style={shape=circle,fill=black,minimum size=2pt,
			inner sep=0pt,outer sep=2pt},]
		\node (A1) at (0,1) {$\Pi$};
		\node (A2) at (1,1) {$\Gamma$};
		\node (B1) at (0,0) {$\Lambda$};
		\node (B2) at (1,0) {$\Omega$};
		\node (D) at (-1,0) {$\star$};
		\node (E) at (-1,1) {$Q^{-1} [y]$};
		\path[->,font=\scriptsize]
		(A1) edge node[left]{$Q$} (B1);
		\path[->,font=\scriptsize]
		(A1) edge node[above]{$P$} (A2)
		(E) edge node[left]{$$} (D)
		(D) edge node[above]{$y$} (B1)
		(E) edge node[right]{$$} (A1)
		(B1) edge node[above]{$\Psi$} (B2)
		(A2) edge node[right]{$\Phi$} (B2);
		\draw (A2) edge[implies] node[above] {\scriptsize$\eta\ $} (B1);
	\end{tikzpicture}
\end{center} (the natural isomorphism being part of the square for the homotopy fiber is suppressed in the notation) and by the universal property of the homotopy fiber $\Phi^{-1}[\Psi(y)]$ a functor $F: Q^{-1}[y] \to \Phi^{-1}[\Psi(y)]$, which is explicitly given by
\begin{align} F (z,g) := (P(z),\Psi(g) \eta_z) \myforall(z,g) \in  Q^{-1}[y]. \label{derivationpullpush1}   \end{align}
Now for any 2-vector bundle $\varrho$ over $\Omega$ taking values in 2-vector spaces, $s\in  \Par \Phi^* \varrho$ and $y\in \Lambda$ we define the vector bundle with values in the 2-vector space $\varrho(\Psi(y))$
\begin{align}
	\xi : \Phi^{-1}[\Psi(y)] \to \varrho(\Psi(y)), \quad (x,h) \mapsto \varrho(h)s(x)=h.s(x).
\end{align} 
We should emphasize that $\xi$ is an ordinary vector bundle, so we obtain a pullback map
\begin{align}
	F^* : \Par \xi \to \Par F^* \xi, \label{derivationpullpush2}  
\end{align} and a pushforward map
\begin{align}
	F_* : \Par F^* \xi \to \Par \xi 
	\label{derivationpullpush3}  
\end{align} by the constructions recalled in Section~\ref{secrecallpullpush}.
Of course, the auxiliary objects $\xi$ and $F$ depend on all the functors involved and on $s$ and $y$, which is suppressed in the notation. Next observe
\begin{align}
	\Par \xi = \lim_{(x,h) \in \Phi^{-1}[\Psi(y)]} \varrho(h)s(x) = (\Psi^* \Phi_* s)(y)
\end{align} and
\begin{align}
	\Par F^* \xi = \lim_{(z,g) \in Q^{-1}[y]} \varrho(\Psi(g)\eta_z) s(P(z)) = (Q_* \varrho(\eta) P^*  s)(y).
\end{align} Restoring the previously suppressed dependence on $s$ and $y$ we obtain maps
\begin{align}
	(\eta^*)_s^y := F^* : (\Psi^* \Phi_* s)(y) \to (Q_* \varrho(\eta)_* P^*  s)(y) \label{eqndefineetapull}
\end{align} and
\begin{align}
	(\eta_*)_s^y := F_* : (Q_* \varrho(\eta)_* P^*  s)(y) \to (\Psi^* \Phi_* s)(y)\label{eqndefineetapush} .
\end{align} If we let $s$ and $y$ run over all parallel sections of $\Phi^* \varrho$ and all objects of $\Lambda$, respectively, they combine into the following natural transformations:

\begin{proposition}\label{satz440transf}
Consider a weakly commuting square 
\begin{center}
	\begin{tikzpicture}[scale=2, implies/.style={double,double equal sign distance,-implies},
		dot/.style={shape=circle,fill=black,minimum size=2pt,
			inner sep=0pt,outer sep=2pt},]
		\node (A1) at (0,1) {$\Pi$};
		\node (A2) at (1,1) {$\Gamma$};
		\node (B1) at (0,0) {$\Lambda$};
		\node (B2) at (1,0) {$\Omega$};
		\path[->,font=\scriptsize]
		(A1) edge node[above]{$P$} (A2)
		(A1) edge node[left]{$Q$} (B1)
		(A2) edge node[right]{$\Phi$} (B2)
		(B1) edge node[below]{$\Psi$} (B2);
		\draw (A2) edge[implies] node[above] {\scriptsize$\eta\ $} (B1);
	\end{tikzpicture}
\end{center} of essentially finite groupoids and a 2-vector bundle $\varrho$ over $\Omega$. Then we have two natural transformations \begin{align}\begin{array}{c}
		\begin{tikzpicture}[scale=1, implies/.style={double,double equal sign distance,-implies},
			dot/.style={shape=circle,fill=black,minimum size=2pt,
				inner sep=0pt,outer sep=2pt},]
			\draw[line width=0.5pt,->]
			(0,2) node {$\Par \Phi^* \varrho$}
			(4,2) node {$\Par \Psi^* \varrho$}
			(2,4.2) node {$\Psi^* \Phi_*$}
			(2,-0.2) node {$Q_* \varrho(\eta)_* P^* $}
			(0,2.6)   \dir{90}{90} (4,2.6) ;
			\draw[line width=0.5pt,->]
			(0,1.4)   \dir{270}{270} (4,1.5) 
			;
			\draw (1.5,3.2) edge[implies] node[left] {$\eta^*\ $} (1.5,0.8);
			\draw (2.5,0.8) edge[implies] node[right] {$\eta_*\ $} (2.5,3.2);
		\end{tikzpicture} 
	\end{array}
\end{align} of 2-linear functors by the construction just given. We call $\eta^*$ the pull map and $\eta_*$ the push map.  
\end{proposition}

\noindent We will need the special case in which the square in Proposition~\ref{satz440transf} is a homotopy pullback square. For the investigation of this special case we need the following easy Lemma:

\begin{lemma}\label{lemmahomotopyfibersandequivalences}
Let $\Phi : \Gamma \to \Omega$ be an equivalence of groupoids, then $\Phi^{-1}[y]$ is equivalent to the groupoid consisting of one object $(x,g)$, where $x\in \Gamma$ and $g: \Phi(x) \to y$, and trivial automorphism group.
\end{lemma}

\noindent Now if the square in Proposition~\ref{satz440transf} is a homotopy pullback square, then the functor $F$ from \eqref{derivationpullpush1} is an equivalence by the fiberwise characterization of homotopy pullbacks in \cite[5.2]{cpshtp}. 
Using Lemma~\ref{lemmahomotopyfibersandequivalences} we can deduce that in this case $F^*$ from \eqref{derivationpullpush2} is inverse to $F_*$ from \eqref{derivationpullpush3}. Interpreting this in terms of $\eta^*$ and $\eta_*$, see \eqref{eqndefineetapull} and \eqref{eqndefineetapush}, we obtain a generalization of the equivariant Beck-Chevalley condition in \cite{schweigertwoikeofk} to 2-vector bundles:

\begin{proposition}[Equivariant Beck-Chevalley condition for 2-vector bundles]\label{mscextversionvonsatz440ausbscwoikeext}
For the homotopy pullback 
\begin{center}
	\begin{tikzpicture}[scale=2, implies/.style={double,double equal sign distance,-implies},
		dot/.style={shape=circle,fill=black,minimum size=2pt,
			inner sep=0pt,outer sep=2pt},]
		\node (A1) at (0,1) {$\Gamma\times_\Omega \Lambda$};
		\node (A2) at (1,1) {$\Gamma$};
		\node (B1) at (0,0) {$\Lambda$};
		\node (B2) at (1,0) {$\Omega$};
		\path[->,font=\scriptsize]
		(A1) edge node[above]{$\pi_\Gamma$} (A2)
		(A1) edge node[left]{$\pi_\Lambda$} (B1)
		(A2) edge node[right]{$\Phi$} (B2)
		(B1) edge node[below]{$\Psi$} (B2);
		\draw (A2) edge[implies] node[above] {\scriptsize$\eta\ $} (B1);
	\end{tikzpicture}
\end{center}
of a cospan $\Lambda \stackrel{\Psi}{\to} \Omega \stackrel{\Phi}{\longleftarrow} \Gamma$ of essentially finite groupoids and any 2-vector bundle $\varrho$ over $\Omega$ the pentagon relating different pull-push combinations
\begin{center}
	\begin{tikzpicture}[scale=2.5, implies/.style={double,double equal sign distance,-implies},
		dot/.style={shape=circle,fill=black,minimum size=2pt,
			inner sep=0pt,outer sep=2pt},]
		\node (A1) at (0,2) {$\Par \Phi^*\varrho$};
		\node (A2) at (2,2) {$\Par \varrho$};
		\node (C) at (0,1) {$\Par \pi_\Gamma^* \Phi^* \varrho$};
		\node (B1) at (1,0) {$\Par \pi_\Lambda^* \Psi^* \varrho$};
		\node (B2) at (2,0) {$\Par \Psi^* \varrho$};
		\path[->,font=\scriptsize]
		(A1) edge node[above]{$\Phi_*$} (A2)
		(A1) edge node[left]{$\pi_\Gamma^*$} (C)
		(C) edge node[left]{$\varrho(\eta)_*$} (B1)
		(B1) edge node[above]{${\pi_\Lambda}_*$} (B2)
		(A2) edge node[right]{$\Psi^*$} (B2);
		\draw (A2) edge[implies] node[above] {\scriptsize$\eta^*\ $} (0.7,0.7);
	\end{tikzpicture}
\end{center} commutes up to the natural isomorphism $\eta^* : \Psi^* \Phi_* \Rightarrow {\pi_\Lambda}_* \varrho(\eta)_* \pi_\Gamma^*$ with inverse $\eta_* :  Q_* \varrho(\eta)_* P^* \Rightarrow \Psi^* \Phi_*$.
\end{proposition}

\noindent For later purposes we work out a special case of Proposition~\ref{satz440transf}. To this end, let us look at a special case of the pushforward maps from Section~\ref{secrecallpullpush}: Let $V$ be an object in a 2-vector space $\cat{V}$ and $\Gamma$ an essentially finite groupoid. Then there is a constant vector bundle $\xi_V$ assigning $V$ to all objects in $\Gamma$ and the identity on $V$ to all morphisms. Obviously, $\Par \xi_V$, i.e.\ the limit of the functor $\xi_V$, is given by the product $\prod_{\pi_0(\Gamma)} V$, but this is also a coproduct, hence
\begin{align}
	\Par \xi_V = \coprod_{\pi_0(\Gamma)} V.\end{align} Now note that $\xi_V = t^* V$, where $t: \Gamma \to \star$ is the functor to the terminal groupoid $\star$ and $V$ is the object seen as vector bundle of $\star$. The pushforward along $t$ is now a map
\begin{align}
	\int_\Gamma := t_* : \coprod_{\pi_0(\Gamma)} V \to V
\end{align} that we call integral with respect to groupoid cardinality. Recalling the definition of the pushforward (Section~\ref{secrecallpullpush}), we see that on the summand belong to $[x]\in \pi_0(\Gamma)$ it is given by the endomorphism 
\begin{align} \frac{1}{|\! \Aut(x)|} \cdot \id_V : V \to V. \end{align}

\begin{corollary}\label{korbeckchevgrpdcard}
Consider a weakly commuting square 
\begin{center}
	\begin{tikzpicture}[scale=2, implies/.style={double,double equal sign distance,-implies},
		dot/.style={shape=circle,fill=black,minimum size=2pt,
			inner sep=0pt,outer sep=2pt},]
		\node (A1) at (0,1) {$\Pi$};
		\node (A2) at (1,1) {$\Gamma$};
		\node (B1) at (0,0) {$\Lambda$};
		\node (B2) at (1,0) {$\Gamma$};
		\path[->,font=\scriptsize]
		(A1) edge node[above]{$P$} (A2)
		(A1) edge node[left]{$Q$} (B1)
		(A2) edge node[right]{$\id_\Gamma$} (B2)
		(B1) edge node[below]{$\Psi$} (B2);
		\draw (A2) edge[implies] node[above] {\scriptsize$\eta\ $} (B1);
	\end{tikzpicture}
\end{center} of essentially finite groupoids and a given 2-vector bundle $\varrho$ over $\Gamma$. Then the natural transformation
\begin{align}
	\eta_* : Q_* \varrho(\eta)_* P^* \Longrightarrow \Psi^* 
\end{align}admits the following explicit description: For $s\in \Par \varrho$ and $y\in \Lambda$ we obtain the commuting diagram
\begin{center}
	\begin{tikzpicture}[scale=2, implies/.style={double,double equal sign distance,-implies},
		dot/.style={shape=circle,fill=black,minimum size=2pt,
			inner sep=0pt,outer sep=2pt},]
		\node (A1) at (0,1) {$(Q_* \varrho(\eta)_* P^* s)(y)$};
		\node (A2) at (2,1) {${\displaystyle\coprod_{\pi_0(Q^{-1}[y])} s(\Psi(y))}$};
		\node (B1) at (1,0) {$(\Psi^* s)(y)=s(\Psi(y))$};
		\path[->,font=\scriptsize]
		(A1) edge node[above]{$\cong$} (A2)
		(A1) edge node[left]{$\eta_*$} (B1)
		(A2) edge node[right]{$\quad \int_{Q^{-1}[y]}$} (B1);
	\end{tikzpicture},
\end{center} where $\cong$ denotes a natural isomorphism to $\coprod_{\pi_0(Q^{-1}[y])} s(\Psi(y))$. This expresses $\eta_*$ as an integral with respect to groupoid cardinality.
\end{corollary}

\begin{proof}
We use that $s$ is parallel to see
\begin{align} (Q_* \varrho(\eta) P^* s)(y) = \lim_{(z,g) \in Q^{-1}[y]} \varrho(\Psi(g)\eta) s(P(z))  \cong \lim_{(z,g) \in Q^{-1}[y]}  s(\Psi(y)).\end{align} The limit of this last constant diagram is obviously given by the finite product \begin{align}\prod_{\pi_0(Q^{-1}[y])} s(\Psi(y)),\end{align} which coincides with the finite coproduct $\coprod_{\pi_0(Q^{-1}[y])} s(\Psi(y))$ since $s(\Psi(y))$ is an object in a 2-vector space. By definition and Lemma~\ref{lemmahomotopyfibersandequivalences} the component $(Q_* \varrho(\eta) P^* s)(y) \to s(\Psi(y))$ of $\eta_*$ is the pushforward map along the functor $Q^{-1}[y] \to \Psi(y)$, where $\Psi(y)$ is the discrete groupoid with one object $\Psi(y)$. This implies the claim.    
\end{proof}

\section{The parallel section functor on the symmetric monoidal bicategory  $\TRepGrpd$\label{secparallelsectionfunctor}}
In this section we formulate and prove the main result of this paper: We introduce the symmetric monoidal bicategory $\TRepGrpd$, which is built from 2-vector bundles and spans of groupoids (Section~\ref{secdefdefTRepGrpd}) and construct the parallel section functor
\begin{align} \Par : \TRepGrpd \to  \TwoVect\end{align} with values in 2-vector spaces.

\subsection{The symmetric monoidal bicategory $\TRepGrpd$\label{secdefdefTRepGrpd}}
In this subsection we introduce the symmetric monoidal bicategory $\TRepGrpd$ needed for the definition of the parallel section functor. Objects are 2-vector bundles over essentially finite groupoids. The 1-morphisms and 2-morphisms come from spans of groupoids and spans of spans of groupoids decorated with (higher) intertwiners. We show below in Proposition~\ref{satzendtrepgrpd} that $\TRepGrpd$ generalizes the symmetric monoidal category $\RepGrpd{}$, which was the domain of the parallel section functor in \cite{schweigertwoikeofk}.\footnote{In \cite{schweigertwoikeofk} we additionally specified a field, but we suppress this here and always work over the complex numbers. } We remark that, in a homotopical setting, related (higher) span categories for (higher) vector bundles have been discussed in \cite{haugseng}.

\begin{definition}[The symmetric monoidal bicategory $\TRepGrpd$]\label{defTRepGrpd}
	We define the symmetric monoidal bicategory $\TRepGrpd$ as follows:
	\begin{enumerate}
		\item Objects are 2-vector bundles over essentially finite groupoids, i.e.\ pairs $(\Gamma,\varrho)$ where $\varrho$ is a 2-vector bundle over an essentially finite groupoid $\Gamma$. 
		
		\item A 1-morphism $(\Gamma_0,\varrho_0) \to (\Gamma_1,\varrho_1)$ is a span
		\begin{align} \Gamma_0 \stackrel{r_0}{\longleftarrow} \Lambda \stackrel{r_1}{\to} \Gamma_1 \end{align} of essentially finite groupoids together with an intertwiner
		\begin{align} \lambda : r_0^* \varrho_0 \to r_1^* \varrho_1 ,\end{align} i.e.\ a 1-morphism of the 2-vector bundles $r_0^* \varrho_0$ and $r_1^* \varrho_1$ over $\Lambda$. We denote such a 1-morphism by $(\Gamma_0,\varrho_0) \stackrel{r_0}{\longleftarrow} (\Lambda,\lambda) \stackrel{r_1}{\to} (\Gamma_1,\varrho_1)$.
		
		\item A 2-morphism from $(\Gamma_0,\varrho_0) \stackrel{r_0}{\longleftarrow} (\Lambda,\lambda) \stackrel{r_1}{\to} (\Gamma_1,\varrho_1)$ to $(\Gamma_0,\varrho_0) \stackrel{r_0'}{\longleftarrow} (\Lambda,\lambda) \stackrel{r_1'}{\to} (\Gamma_1,\varrho_1)$ is an equivalence of class (as explained in Remark~\ref{bmkTRepGrpd}, \ref{bmkTRepGrpd2} below) of \begin{itemize}
			\item a span of spans, i.e.\ a diagram
			
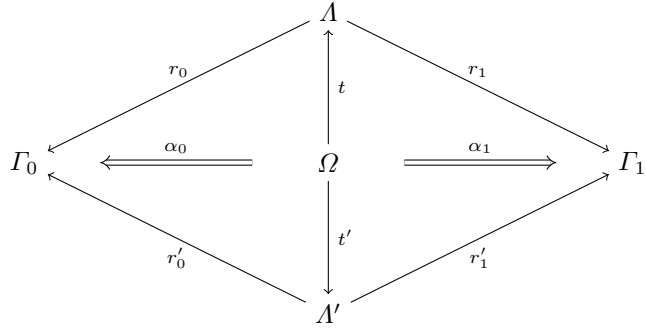
\begin{center}
				\begin{tikzpicture}[scale=2,     implies/.style={double,double equal sign distance,-implies},
					dot/.style={shape=circle,fill=black,minimum size=2pt,
						inner sep=0pt,outer sep=2pt},]
					\node (A1) at (0,0) {$\Gamma_0$};
					\node (A2) at (2,1) {$\Lambda$};
					\node (A3) at (4,0) {$\Gamma_1$};
					\node (B2) at (2,-1) {$\Lambda'$};
					\node (C) at (2,0) {$\Omega$};
					\node (B1) at (1,-0.5) {$$};
					\node (B3) at (2,-1) {$$};
					\path[->,font=\scriptsize]
					(A2) edge node[above]{$r_0$} (A1)
					(A2) edge node[above]{$r_1$} (A3)
					(B2) edge node[below]{$r_0'$} (A1)
					(B2) edge node[below]{$r_1'$} (A3)
					(C) edge node[right]{$t'$} (B2)
					(C) edge node[right]{$t$} (A2);
					\draw (1.5,0) edge[implies] node[above] {\scriptsize$\alpha_0$} (0.5,0);
					\draw (2.5,0) edge[implies] node[above] {\scriptsize$\alpha_1$} (3.5,0);
				\end{tikzpicture}
			\end{center} in essentially finite groupoids commutative up to the indicated natural isomorphisms (in this way of presentation the direction of the natural isomorphism is obtained by reading from top to bottom, e.g. $\alpha_0$ is a natural isomorphism $r_0t \Rightarrow r_0' t'$).
			
			\item together with a natural morphism
			\begin{center}
				\begin{tikzpicture}[scale=1.5,     implies/.style={double,double equal sign distance,-implies},
					dot/.style={shape=circle,fill=black,minimum size=2pt,
						inner sep=0pt,outer sep=2pt},]
					\node (A1) at (0,1) {$(r_0t)^* \varrho_0=t^* r_0^* \varrho_0$};
					\node (A2) at (4,1) {$t^* r_1^* \varrho_1 = (r_1t)^* \varrho_1$};
					\node (B1) at (0,0) {$(r_0't')^* \varrho_0 = {t'}^* {r_0'}^* \varrho_0$};
					\node (B2) at (4,0) {${t'}^* {r_1'}^* \varrho_1 = (r_1't')^* \varrho_1$};
					\path[->,font=\scriptsize]
					(A1) edge node[above]{$t^* \lambda$} (A2)
					(A1) edge node[left]{$\varrho_0({\alpha_0})$} (B1)
					(B1) edge node[above]{${t'}^* \lambda'$} (B2)
					(A2) edge node[right]{$\varrho_1({\alpha_1})$} (B2);
					\draw (A2) edge[implies] node[above] {\scriptsize$\omega$} (B1);
				\end{tikzpicture}
			\end{center}
			We will denote this 2-morphism by
			
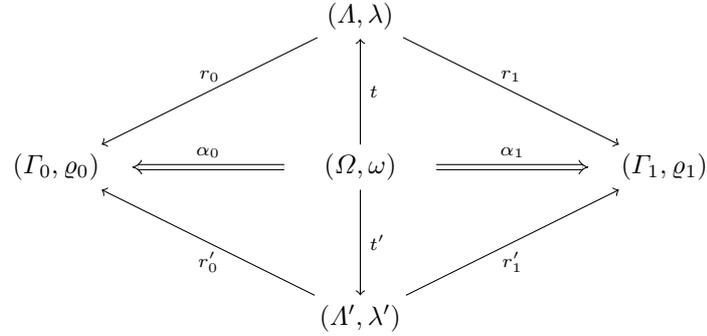
\begin{center}
				\begin{tikzpicture}[scale=2,     implies/.style={double,double equal sign distance,-implies},
					dot/.style={shape=circle,fill=black,minimum size=2pt,
						inner sep=0pt,outer sep=2pt},]
					\node (A1) at (0,0) {$(\Gamma_0,\varrho_0)$};
					\node (A2) at (2,1) {$(\Lambda,\lambda)$};
					\node (A3) at (4,0) {$(\Gamma_1,\varrho_1)$};
					\node (B2) at (2,-1) {$(\Lambda',\lambda')$};
					\node (C) at (2,0) {$(\Omega,\omega)$};
					\node (B1) at (1,-0.5) {$$};
					\node (B3) at (2,-1) {$$};
					\path[->,font=\scriptsize]
					(A2) edge node[above]{$r_0$} (A1)
					(A2) edge node[above]{$r_1$} (A3)
					(B2) edge node[below]{$r_0'$} (A1)
					(B2) edge node[below]{$r_1'$} (A3)
					(C) edge node[right]{$t'$} (B2)
					(C) edge node[right]{$t$} (A2);
					\draw (1.5,0) edge[implies] node[above] {\scriptsize$\alpha_0$} (0.5,0);
					\draw (2.5,0) edge[implies] node[above] {\scriptsize$\alpha_1$} (3.5,0);
				\end{tikzpicture}
			\end{center}
			
		\end{itemize}
		
		\item Composition of 1-morphisms: For 1-morphisms
		\begin{align}
			(\Gamma_0,\varrho_0) \stackrel{r_0}{\longleftarrow} (\Lambda,\lambda) \stackrel{r_1}{\to} (\Gamma_1,\varrho_1) 
		\end{align} and
		\begin{align} (\Gamma_1,\varrho_1) \stackrel{r_1'}{\longleftarrow} (\Lambda',\lambda') \stackrel{r_2'}{\to} (\Gamma_2,\varrho_2) \end{align} the composition $(\Lambda',\lambda') \circ (\Lambda,\lambda)$ is by homotopy pullback. More precisely, the span part of the composition is the outer span of
		
		
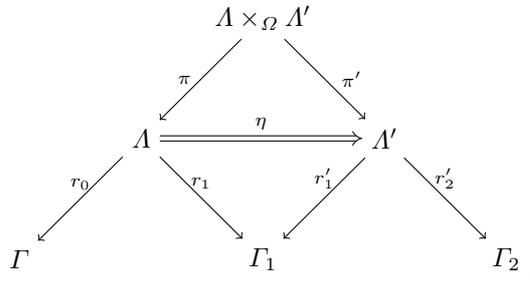
\begin{center}
			\begin{tikzpicture}[scale=0.8, implies/.style={double,double equal sign distance,-implies},
				dot/.style={shape=circle,fill=black,minimum size=2pt,
					inner sep=0pt,outer sep=2pt},]
				\node (A1) at (0,0) {$\Gamma$};
				\node (A2) at (2,2) {$\Lambda$};
				\node (A3) at (4,0) {$\Gamma_1$};
				\node (C) at (4,0.2) {$$};
				\node (A4) at (6,2) {$\Lambda'$};
				\node (A5) at (8,0) {$\Gamma_2$};
				\node (B2) at (4,4) {$\Lambda \times_\Omega \Lambda'$};
				\path[->,font=\scriptsize]
				(A2) edge node[above]{$r_0$} (A1)
				(A2) edge node[above]{$r_1$} (A3)
				(A4) edge node[above]{$r_1'$} (A3)
				(A4) edge node[above]{$r_2'$} (A5)
				(B2) edge node[left]{$\pi$} (A2)
				(B2) edge node[right]{$\ \pi'$} (A4);
				\draw (A2) edge[implies] node[above] {\scriptsize$\eta$} (A4);
			\end{tikzpicture},
		\end{center} where $\eta$ is the natural transformation the homotopy pullback comes equipped with, together with the intertwiner
		\begin{align}
			\lambda \times_{\Gamma_1} \lambda' : \pi^* r_0^* \varrho_0 \xrightarrow{\pi^* \lambda} \pi^* r_1^* \varrho_1 \xrightarrow{\varrho_1(\eta)} {\pi'}^* {r_1'}^* \varrho_1 \xrightarrow{{\pi'}^* \lambda'}  {\pi'}^* {r_2'}^* \varrho_2.
		\end{align}\label{defTRepGrpdc1}

		\item Vertical composition of 2-morphisms: The vertical composition of the 2-morphisms\label{defTRepGrpd2mor}
		\begin{center}
			\begin{tikzpicture}[scale=2,     implies/.style={double,double equal sign distance,-implies},
				dot/.style={shape=circle,fill=black,minimum size=2pt,
					inner sep=0pt,outer sep=2pt},]
				\node (A1) at (0,0) {$(\Gamma_0,\varrho_0)$};
				\node (A2) at (2,1) {$(\Lambda,\lambda)$};
				\node (A3) at (4,0) {$(\Gamma_1,\varrho_1)$};
				\node (B2) at (2,-1) {$(\Lambda',\lambda')$};
				\node (C) at (2,0) {$(\Omega,\omega)$};
				\node (B1) at (1,-0.5) {$$};
				\node (B3) at (2,-1) {$$};
				\path[->,font=\scriptsize]
				(A2) edge node[above]{$r_0$} (A1)
				(A2) edge node[above]{$r_1$} (A3)
				(B2) edge node[below]{$r_0'$} (A1)
				(B2) edge node[below]{$r_1'$} (A3)
				(C) edge node[right]{$t'$} (B2)
				(C) edge node[right]{$t$} (A2);
				\draw (1.5,0) edge[implies] node[above] {\scriptsize$\alpha_0$} (0.5,0);
				\draw (2.5,0) edge[implies] node[above] {\scriptsize$\alpha_1$} (3.5,0);
			\end{tikzpicture}
		\end{center} and 
		\begin{center}
			\begin{tikzpicture}[scale=2,     implies/.style={double,double equal sign distance,-implies},
				dot/.style={shape=circle,fill=black,minimum size=2pt,
					inner sep=0pt,outer sep=2pt},]
				\node (A1) at (0,0) {$(\Gamma_0,\varrho_0)$};
				\node (A2) at (2,1) {$(\Lambda',\lambda')$};
				\node (A3) at (4,0) {$(\Gamma_1,\varrho_1)$};
				\node (B2) at (2,-1) {$(\Lambda'',\lambda'')$};
				\node (C) at (2,0) {$(\widetilde{\Omega},\widetilde{\omega})$};
				\node (B1) at (1,-0.5) {$$};
				\node (B3) at (2,-1) {$$};
				\path[->,font=\scriptsize]
				(A2) edge node[above]{$r_0'$} (A1)
				(A2) edge node[above]{$r_1'$} (A3)
				(B2) edge node[below]{$r_0''$} (A1)
				(B2) edge node[below]{$r_1''$} (A3)
				(C) edge node[right]{$u''$} (B2)
				(C) edge node[right]{$u'$} (A2);
				\draw (1.5,0) edge[implies] node[above] {\scriptsize$\beta_0$} (0.5,0);
				\draw (2.5,0) edge[implies] node[above] {\scriptsize$\beta_1$} (3.5,0);
			\end{tikzpicture}
		\end{center} is the 2-morphism
		\begin{align}\begin{array}{c}
				\begin{tikzpicture}[scale=2,     implies/.style={double,double equal sign distance,-implies},
					dot/.style={shape=circle,fill=black,minimum size=2pt,
						inner sep=0pt,outer sep=2pt},]
					\node (A1) at (0,0) {$(\Gamma_0,\varrho_0)$};
					\node (A2) at (2,1) {$(\Lambda,\lambda)$};
					\node (A3) at (4,0) {$(\Gamma_1,\varrho_1)$};
					\node (B2) at (2,-1) {$(\Lambda'',\lambda'')$};
					\node (C) at (2,0) {$(\Omega \times_{\Lambda'} \widetilde \Omega,\omega \times_{\lambda'} \widetilde \omega)$};
					\node (B1) at (1,-0.5) {$$};
					\node (B3) at (2,-1) {$$};
					\path[->,font=\scriptsize]
					(A2) edge node[above]{$r_0$} (A1)
					(A2) edge node[above]{$r_1$} (A3)
					(B2) edge node[below]{$r_0''$} (A1)
					(B2) edge node[below]{$r_1''$} (A3)
					(C) edge node[right]{$v''$} (B2)
					(C) edge node[right]{$v$} (A2);
					\draw (1,0) edge[implies] node[above] {\scriptsize$\gamma_0$} (0.5,0);
					\draw (3,0) edge[implies] node[above] {\scriptsize$\gamma_1$} (3.5,0);
			\end{tikzpicture}\end{array}, \label{verticalcompositioneqn}
		\end{align}
		whose components are defined as follows:
		\begin{itemize}
			
			\item $\Omega \times_{\Lambda'} \widetilde \Omega$ denotes the homotopy pullback
			\begin{center}
				\begin{tikzpicture}[scale=2, implies/.style={double,double equal sign distance,-implies},
					dot/.style={shape=circle,fill=black,minimum size=2pt,
						inner sep=0pt,outer sep=2pt},]
					\node (A1) at (0,1) {$\Omega \times_{\Lambda'} \widetilde \Omega$};
					\node (A2) at (1,1) {$\Omega$};
					\node (B1) at (0,0) {$\widetilde \Omega$};
					\node (B2) at (1,0) {$\Lambda'$};
					\path[->,font=\scriptsize]
					(A1) edge node[above]{$q$} (A2)
					(A1) edge node[left]{$\widetilde q$} (B1)
					(A2) edge node[right]{$t'$} (B2)
					(B1) edge node[below]{$u'$} (B2);
					\draw (A2) edge[implies] node[above] {\scriptsize$\eta\ $} (B1);
				\end{tikzpicture},
			\end{center} which introduces $q$, $\widetilde q$ and $\eta$,
			
			\item the functors $v$ and $v''$ in \eqref{verticalcompositioneqn} are defined by $v:= tq$ and $v'' := u'' \widetilde q$,
			
			\item the natural transformations $\gamma_0$ and $\gamma_1$ in \eqref{verticalcompositioneqn} are obtained by 
			\begin{align}
				\gamma_0 : r_0v= r_0tq \stackrel{\alpha_0}{\Longrightarrow} r_0' t' q \stackrel{\eta}{\Longrightarrow} r_0'u' \widetilde q \stackrel{\beta_0}{\Longrightarrow} r_0'' u'' \widetilde q = r_0'' v''
			\end{align} and
			\begin{align}
				\gamma_1 : r_1v=r_1tq \stackrel{\alpha_1}{\Longrightarrow} r_1' t' q \stackrel{\eta}{\Longrightarrow} r_1'u' \widetilde q \stackrel{\beta_1}{\Longrightarrow} r_1'' u'' \widetilde q = r_1'' v'',
			\end{align} where identity transformations are suppressed in the notation,
			
			\item and the natural morphism $\omega \times_{\lambda'} \widetilde \omega$ is obtained as the composition
			\begin{center}
				\begin{tikzpicture}[scale=1.5,     implies/.style={double,double equal sign distance,-implies},
					dot/.style={shape=circle,fill=black,minimum size=2pt,
						inner sep=0pt,outer sep=2pt},]
					\node (A1) at (0,1) {$(r_0v)^* \varrho_0 = q^* t^* r_0^* \varrho_0$};
					\node (A2) at (4,1) {$q^* t^* r_1^*  \varrho_1 = (r_1v)^* \varrho_1$};
					\node (B1) at (0,0) {$q^* {t'}^* {r_0'}^*  \varrho_0$};
					\node (B2) at (4,0) {$q^* {t'}^* {r_1'}^* \varrho_1$};
					\node (C1) at (0,-1) {$\widetilde q^* {u'}^* {r_0'}^*   \varrho_0$};
					\node (C2) at (4,-1) {$\widetilde q^* {u'}^* {r_1'}^* \varrho_1$};
					\node (D1) at (0,-2) {$\widetilde q^* {u''}^* {r_0''}^* \varrho_0$};
					\node (D2) at (4,-2) {$\widetilde q^* {u''}^* {r_1''}^* \varrho_1$};
					\path[->,font=\scriptsize]
					(A1) edge node[above]{$q^* t^* \lambda$} (A2)
					(A1) edge node[left]{$q^*\varrho_0({\alpha_0})$} (B1)
					(B1) edge node[left]{$\varrho_0(\eta)$} (C1)
					(B2) edge node[right]{$\varrho_1(\eta)$} (C2)
					(C1) edge node[above]{$\widetilde q^* {u'}^* \lambda'$} (C2)
					(D1) edge node[above]{$\widetilde q^* {u''}^* \lambda''$} (D2)
					(B1) edge node[above]{$q^*{t'}^* \lambda'$} (B2)
					(A2) edge node[right]{$q^*\varrho_1({\alpha_1})$} (B2)
					(C1) edge node[left]{${\widetilde q} ^*\varrho_0(\beta_0)$} (D1)
					(C2) edge node[right]{${\widetilde q} ^*\varrho_1(\beta_1)$} (D2);
					\draw (A2) edge[implies] node[above] {\scriptsize$q^* \omega$} (B1);
					\draw (B2) edge[implies] node[above] {\scriptsize$\theta'$} (C1);
					\draw (C2) edge[implies] node[above] {\scriptsize$\widetilde q^* \omega'$} (D1);
				\end{tikzpicture},
			\end{center} where the middle square is decorated with the isomorphism $\theta'$ belonging to $\lambda'$ (Remark~\ref{bmkdescr2rep}, \ref{bmkdescr2rep2}), i.e.\ the evaluation of the 2-morphisms on the middle square on $(z,\widetilde z,g) \in \Omega \times_{\Lambda'} \widetilde \Omega$ is given by 
			\begin{center}
				\begin{tikzpicture}[scale=2, implies/.style={double,double equal sign distance,-implies},
					dot/.style={shape=circle,fill=black,minimum size=2pt,
						inner sep=0pt,outer sep=2pt},]
					\node (A1) at (0,1) {$\varrho_0(r_0't'(z))$};
					\node (A2) at (2,1) {$\varrho_1(r_1't'(z))$};
					\node (B1) at (0,0) {$\varrho_0(r_0'u'(\widetilde z))$};
					\node (B2) at (2,0) {$\varrho_1(r_1'u'(\widetilde z))$};
					\path[->,font=\scriptsize]
					(A1) edge node[above]{$\lambda'_{t'(z)}$} (A2)
					(A1) edge node[left]{$\varrho_0({r_0'(g)})$} (B1)
					(A2) edge node[right]{$\varrho_1({r_1'(g)})$} (B2)
					(B1) edge node[below]{$\lambda'_{u'(\widetilde z)}$} (B2);
					\draw (A2) edge[implies] node[above] {\scriptsize$\theta_g\ $} (B1);
				\end{tikzpicture}.
			\end{center}

		\end{itemize}

		\item Horizontal composition of 2-morphisms: The horizontal composition of the 2-morphisms\label{defTRepGrpd2morh}
		\begin{center}
			\begin{tikzpicture}[scale=2,     implies/.style={double,double equal sign distance,-implies},
				dot/.style={shape=circle,fill=black,minimum size=2pt,
					inner sep=0pt,outer sep=2pt},]
				\node (A1) at (0,0) {$(\Gamma_0,\varrho_0)$};
				\node (A2) at (2,1) {$(\Lambda,\lambda)$};
				\node (A3) at (4,0) {$(\Gamma_1,\varrho_1)$};
				\node (B2) at (2,-1) {$(\Lambda',\lambda')$};
				\node (C) at (2,0) {$(\Omega,\omega)$};
				\node (B1) at (1,-0.5) {$$};
				\node (B3) at (2,-1) {$$};
				\path[->,font=\scriptsize]
				(A2) edge node[above]{$r_0$} (A1)
				(A2) edge node[above]{$r_1$} (A3)
				(B2) edge node[below]{$r_0'$} (A1)
				(B2) edge node[below]{$r_1'$} (A3)
				(C) edge node[right]{$t'$} (B2)
				(C) edge node[right]{$t$} (A2);
				\draw (1.5,0) edge[implies] node[above] {\scriptsize$\alpha_0$} (0.5,0);
				\draw (2.5,0) edge[implies] node[above] {\scriptsize$\alpha_1$} (3.5,0);
			\end{tikzpicture}
		\end{center} and 
		\begin{center}
			\begin{tikzpicture}[scale=2,     implies/.style={double,double equal sign distance,-implies},
				dot/.style={shape=circle,fill=black,minimum size=2pt,
					inner sep=0pt,outer sep=2pt},]
				\node (A1) at (0,0) {$(\Gamma_1,\varrho_1)$};
				\node (A2) at (2,1) {$(\Pi,\pi)$};
				\node (A3) at (4,0) {$(\Gamma_2,\varrho_2)$};
				\node (B2) at (2,-1) {$(\Pi',\pi')$};
				\node (C) at (2,0) {$(\widetilde{\Omega},\widetilde{\omega})$};
				\node (B1) at (1,-0.5) {$$};
				\node (B3) at (2,-1) {$$};
				\path[->,font=\scriptsize]
				(A2) edge node[above]{$v_1$} (A1)
				(A2) edge node[above]{$v_2$} (A3)
				(B2) edge node[below]{$v_1'$} (A1)
				(B2) edge node[below]{$v_2'$} (A3)
				(C) edge node[right]{$u'$} (B2)
				(C) edge node[right]{$u$} (A2);
				\draw (1.5,0) edge[implies] node[above] {\scriptsize$\beta_0$} (0.5,0);
				\draw (2.5,0) edge[implies] node[above] {\scriptsize$\beta_1$} (3.5,0);
			\end{tikzpicture}
		\end{center} is
		\begin{center}
			\begin{tikzpicture}[scale=2,     implies/.style={double,double equal sign distance,-implies},
				dot/.style={shape=circle,fill=black,minimum size=2pt,
					inner sep=0pt,outer sep=2pt},]
				\node (A1) at (0,0) {$(\Gamma_0,\varrho_0)$};
				\node (A2) at (2,1) {$(\Lambda \times_{\Gamma_1} \Pi,\lambda \times_{\Gamma_1} \pi)$};
				\node (A3) at (4,0) {$(\Gamma_2,\varrho_2)$};
				\node (B2) at (2,-1) {$(\Lambda' \times_{\Gamma_1} \Pi',\lambda' \times_{\Gamma_1} \pi')$};
				\node (C) at (2,0) {$(\Omega \times_{\Gamma_1} \widetilde \Omega,\omega \times_{\Gamma_1} \widetilde \omega)$};
				\node (B1) at (1,-0.5) {$$};
				\node (B3) at (2,-1) {$$};
				\path[->,font=\scriptsize]
				(A2) edge node[above]{$r_0p$} (A1)
				(A2) edge node[above]{$v_2q$} (A3)
				(B2) edge node[below]{$r_0'p'$} (A1)
				(B2) edge node[below]{$v_2'q'$} (A3)
				(C) edge node[right]{$d$} (B2)
				(C) edge node[right]{$c$} (A2);
				\draw (1,0) edge[implies] node[above] {\scriptsize$\delta_0$} (0.5,0);
				\draw (3,0) edge[implies] node[above] {\scriptsize$\delta_1$} (3.5,0);
			\end{tikzpicture},
		\end{center}
		where
		\begin{itemize}
			
			\item $p: \Lambda \times_{\Gamma_1} \Pi\to \Lambda$, $q: \Lambda \times_{\Gamma_1} \Pi\to \Pi$, $p': \Lambda' \times_{\Gamma_1} \Pi'\to \Lambda'$, $q': \Lambda' \times_{\Gamma_1} \Pi'\to \Pi'$ are the projection functors defined on the respective homotopy pullbacks,
			
			\item $\Omega \times_{\Lambda'} \widetilde \Omega$ denotes the homotopy pullback
			\begin{center}
				\begin{tikzpicture}[scale=2, implies/.style={double,double equal sign distance,-implies},
					dot/.style={shape=circle,fill=black,minimum size=2pt,
						inner sep=0pt,outer sep=2pt},]
					\node (A1) at (0,1) {$\Omega \times_{\Gamma_1} \widetilde \Omega$};
					\node (A2) at (1,1) {$\Omega$};
					\node (B1) at (0,0) {$\widetilde \Omega$};
					\node (B2) at (1,0) {$\Gamma_1$};
					\path[->,font=\scriptsize]
					(A1) edge node[above]{$\ell$} (A2)
					(A1) edge node[left]{$\widetilde \ell$} (B1)
					(A2) edge node[right]{$r_1t$} (B2)
					(B1) edge node[below]{$v_1u$} (B2);
					\draw (A2) edge[implies] node[above] {\scriptsize$\varepsilon \ $} (B1);
				\end{tikzpicture},
			\end{center} which introduces $\ell$, $\widetilde \ell$ and $\varepsilon$; 
			note that the groupoid $\Omega \times_{\Gamma_1} \widetilde \Omega$ together with the natural isomorphism $\varepsilon'$ given as the composition
			\begin{center}
				\begin{tikzpicture}[scale=2, implies/.style={double,double equal sign distance,-implies},
					dot/.style={shape=circle,fill=black,minimum size=2pt,
						inner sep=0pt,outer sep=2pt},]
					\node (B2) at (1,0) {$\Gamma_1$};
					\node (C1) at (-1,0) {$\widetilde \Omega$};
					\node (C2) at (-1,2) {$\Omega \times_{\Gamma_1} \widetilde \Omega$};
					\node (C3) at (1,2) {$\Omega$};
					\node (D) at (1.7,1) {$$};
					\node (E) at (0,-0.7) {$$};
					\path[->,font=\scriptsize]
					(C2) edge node[above]{$\ell$} (C3)
					(C2) edge node[left]{$\widetilde\ell$} (C1)
					(C1) edge node[above]{$v_1u$} (B2)
					(C3) edge node[left]{$r_1t$} (B2);
					\draw (C3) edge[implies] node[above] {\scriptsize$\varepsilon \ $} (C1);
					\draw[line width=0.5pt,->]
					(C3)   \dir{0}{0} node[right] {\scriptsize$r_1't'$} (B2) 
					;
					\draw (D) edge[implies] node[above] {\scriptsize$\alpha_1^{-1}$} (A2);
					\draw[line width=0.5pt,->]
					(C1)   \dir{270}{270} node[below] {\scriptsize$v_1'u'$} (B2) 
					;
					\draw (B1) edge[implies] node[right] {\scriptsize$\beta_0$} (E);
				\end{tikzpicture}
			\end{center}
			is also a homotopy pullback of the cospan defined by the primed functors $r_1't'$ and $v_1' u'$, we call this homotopy pullback $(\Omega \times_{\Gamma_1} \widetilde \Omega)'$,
			
			\item the functor $c: \Omega \times_{\Gamma_1} \widetilde \Omega \to \Lambda \times_{\Gamma_1} \Pi$ is defined using the universal property of the homotopy pullback and makes the diagram
			\begin{center}
				\begin{tikzpicture}[scale=2, implies/.style={double,double equal sign distance,-implies},
					dot/.style={shape=circle,fill=black,minimum size=2pt,
						inner sep=0pt,outer sep=2pt},]
					\node (A1) at (0,1) {$\Lambda \times_{\Gamma_1} \Pi$};
					\node (A2) at (1,1) {$\Lambda$};
					\node (B1) at (0,0) {$\Pi$};
					\node (B2) at (1,0) {$\Gamma_1$};
					\node (C1) at (-1,0) {$\widetilde \Omega$};
					\node (C2) at (-1,2) {$\Omega \times_{\Gamma_1} \widetilde \Omega$};
					\node (C3) at (1,2) {$\Omega$};
					\path[->,font=\scriptsize]
					(A1) edge node[above]{$p$} (A2)
					(A1) edge node[left]{$q$} (B1)
					(A2) edge node[right]{$r_1$} (B2)
					(C2) edge node[above]{$\ell$} (C3)
					(C2) edge node[left]{$\widetilde\ell$} (C1)
					(C2) edge node[above]{$c$} (A1)
					(C1) edge node[above]{$u$} (B1)
					(C3) edge node[right]{$t$} (A2)
					(B1) edge node[below]{$v_1$} (B2);
					\draw (A2) edge[implies] node[above] {\scriptsize$\eta \ $} (B1);
				\end{tikzpicture}
			\end{center}
			commute up to $\varepsilon$ (the squares left blank are also labeled by natural isomorphisms arising from the universal property of the homotopy pullback, but we suppress them in the notation); note here that $c$ arises from the product functor $t\times u$,
			
			\item the functor $d: \Omega \times_{\Gamma_1} \widetilde \Omega \to \Lambda' \times_{\Gamma_1} \Pi'$ is defined analogously using the universal property of the homotopy pullback and makes the diagram
			\begin{center}
				\begin{tikzpicture}[scale=2, implies/.style={double,double equal sign distance,-implies},
					dot/.style={shape=circle,fill=black,minimum size=2pt,
						inner sep=0pt,outer sep=2pt},]
					\node (A1) at (0,1) {$\Lambda' \times_{\Gamma_1} \Pi'$};
					\node (A2) at (1,1) {$\Lambda'$};
					\node (B1) at (0,0) {$\Pi'$};
					\node (B2) at (1,0) {$\Gamma_1$};
					\node (C1) at (-1,0) {$\widetilde \Omega$};
					\node (C2) at (-1,2) {$\Omega \times_{\Gamma_1} \widetilde \Omega$};
					\node (C3) at (1,2) {$\Omega$};
					\node (D) at (1.7,1) {$$};
					\node (E) at (0,-0.7) {$$};
					\path[->,font=\scriptsize]
					(A1) edge node[above]{$p'$} (A2)
					(A1) edge node[left]{$q'$} (B1)
					(A2) edge node[right]{$r_1'$} (B2)
					(C2) edge node[above]{$\ell$} (C3)
					(C2) edge node[left]{$\widetilde\ell$} (C1)
					(C2) edge node[above]{$d$} (A1)
					(C1) edge node[below]{$u'$} (B1)
					(C3) edge node[right]{$t'$} (A2)
					(B1) edge node[below]{$v_1'$} (B2);
					\draw (A2) edge[implies] node[above] {\scriptsize$\eta' \ $} (B1);
					\draw[line width=0.5pt,->]
					(C3)   \dir{0}{0} node[right] {\scriptsize$r_1t$} (B2) 
					;
					\draw (D) edge[implies] node[above] {\scriptsize$\alpha_1$} (A2);
					\draw[line width=0.5pt,->]
					(C1)   \dir{270}{270} node[below] {\scriptsize$v_1u$} (B2) 
					;
					\draw (E) edge[implies] node[right] {\scriptsize$\beta_0$} (B1);
				\end{tikzpicture}
			\end{center}
			commute up to $\varepsilon$ (again, the squares left blank are also labeled by natural isomorphisms arising from the universal property of the homotopy pullback, but we suppress them in the notation); note here that $d$ arises from the product functor $t'\times u'$ and that the composition of the natural isomorphisms in the inner square is $\varepsilon'$,

			\item $\lambda \times_{\Gamma_1} \pi$ is the intertwiner
			\begin{align} \lambda \times_{\Gamma_1} \pi : (r_0p)^* \varrho_0 \cong p^* r_0^* \varrho_0 \xrightarrow{p^* \lambda} p^* r_1^* \varrho_1 \xrightarrow{\varrho_1(\eta)} q^* v_1^* \varrho_1 \xrightarrow{q^* \pi} q^* v_2^* \varrho_2
			\end{align} and $\lambda' \times_{\Gamma_1} \pi'$ is defined analogously (these definitions are recalled here for convenience, although they already follow from the definition of 1-morphisms, see Definition~\ref{defTRepGrpd}, \ref{defTRepGrpdc1}),

			\item $\delta_0: r_0pc \Rightarrow r_0'p'd$ is the natural transformation
			\begin{align}
				r_0pc = r_0t\ell \stackrel{\alpha_0}{\Longrightarrow} r_0't' \ell = r_0' p' d,
			\end{align} and, analogously, 
			$\delta_1: v_2qc \Rightarrow v_2'q'd$ is the natural transformation
			\begin{align}
				v_2qc = v_2u\widetilde \ell \stackrel{ \beta_1}{\Longrightarrow} v_2'u'\widetilde \ell = v_2' q' d,
			\end{align}
			
			\item the natural transformation $\omega\times_{\Gamma_1} \widetilde \omega$ is defined by the commutativity of the diagram
			\begin{center}
				\begin{tikzpicture}[scale=1.7,     implies/.style={double,double equal sign distance,-implies},
					dot/.style={shape=circle,fill=black,minimum size=2pt,
						inner sep=0pt,outer sep=2pt},]
					\node (A1) at (0,1) {$ (r_0pc)^* \varrho_0$};
					\node (A2) at (2,1) {$ (r_1pc)^* \varrho_1$};
					\node (A3) at (4,1) {$ (v_1qc)^* \varrho_1$};
					\node (A4) at (6,1) {$ (v_2qc)^* \varrho_2$};
					\node (B1) at (0,0) {$ (r_0'p'd)^* \varrho_0$};
					\node (B2) at (2,0) {$ (r_1'p'd)^* \varrho_1$};
					\node (B3) at (4,0) {$ (v_1'q'd)^* \varrho_1$};
					\node (B4) at (6,0) {$ (v_2'q'd)^* \varrho_2$};
					\path[->,font=\scriptsize]
					(A1) edge node[left]{$\varrho_0(\delta_0)=\varrho_0(\ell^* \alpha_0)$} (B1)
					(A2) edge node[left]{$\varrho_1(\ell^* \alpha_1)$} (B2)
					(A3) edge node[left]{$\varrho_1(\widetilde\ell^* \beta_0)$} (B3)
					(A4) edge node[right]{$\varrho_2(\widetilde \ell^* \beta_1)=\varrho_2(\delta_1)$} (B4)
					(A1) edge node[above]{$(pc)^* \lambda$} (A2)
					(A2) edge node[above]{$\varrho_1(\varepsilon)$} (A3)
					(A3) edge node[above]{$(qc)^* \pi$} (A4)
					(B1) edge node[above]{$(p'd)^* \lambda'$} (B2)
					(B2) edge node[above]{$\varrho_1(d^*\eta')$} (B3)
					(B3) edge node[above]{$(q'd)^* \pi'$} (B4);
					\draw (A2) edge[implies] node[above] {\scriptsize$\ell^* \omega \ \ $} (B1);
					\draw (A4) edge[implies] node[above] {\scriptsize$\widetilde{\ell}^*\widetilde{\omega} \ \ $} (B3);
				\end{tikzpicture}.
			\end{center}
			
		\end{itemize}

		\item The symmetric monoidal structure is inherited from the symmetric monoidal structure of $\TVecBun(\Gamma)$ for one fixed groupoid (Definition~\ref{defsymmonbicatof2vectorbundles}) and the symmetric monoidal structure on the bicategory of spans of essentially finite groupoids: For objects $X = (\Gamma,\varrho)$ and $Y = (\Omega,\xi)$ the tensor product $X \boxtimes Y$ is the bundle over $\Gamma \times \Omega$ assigning to $(x,y) \in \Gamma \times \Omega$ to 2-vector space $\varrho(x) \boxtimes \xi(y)$. 
		The tensor product of 1- and 2-morphisms is defined analogously by using the Cartesian product of groupoids and the Deligne product.
		The monoidal unit object is the trivial representation of the groupoid $\star$ with one object and trivial automorphism group on the 2-vector space $\FinVect$ of finite-dimensional complex vector spaces; we denote it again by $\star$.

	\end{enumerate}
\end{definition}

\remark	\label{bmkTRepGrpd}
\begin{enumerate}
	
	\item The composition of 1-morphisms in $\TRepGrpd$ requires a model for the homotopy pullback to be chosen. For definiteness we have chosen the one used e.g. in \cite[Definition~3.6]{schweigertwoikeofk}. Choosing a different model yields a composition naturally 2-isomorphic to the initial one. \label{bmkTRepGrpd1}
	
	\item For readability we did not spell out the equivalence relation needed to define 2-morphisms, but only worked with representatives. We justfify this by the fact that in \cite{morton1} such issues were addressed for pure span bicategories (without 2-vector bundles) and in \cite{schweigertwoikeofk} it was explicitly explained in the categorical case how these equivalence relations have to be generalized to take vector bundles into account. This generalization can be done in the bicategorical case as well following the exact same strategy one categorical level higher.  \label{bmkTRepGrpd2}
	
\end{enumerate} \endremark

The next result explains the relation between the symmetric monoidal bicategory $\TRepGrpd$ and the symmetric monoidal category $\RepGrpd{}$ from \cite[Definition~3.7]{schweigertwoikeofk}:

\begin{proposition}\label{satzendtrepgrpd}
The category $\End_{\TRepGrpd}(\star)$ of endomorphisms of the monoidal unit is canonically equivalent, as a symmetric monoidal category, to $\RepGrpd{}$.
\end{proposition}

\begin{proof}
The claim holds more or less by construction, so we only give the main arguments: 
Denote by $\tau : \star \to \FinVect$ the trivial representation of the terminal groupoid $\star$ and by $t: \Gamma \to \star$ the unique functor.
The category $\End_{\TRepGrpd}(\star)$ is symmetric monoidal. Its objects are spans $\star \stackrel{t}{\longleftarrow} \Gamma \stackrel{t}{\to} \star$ together a 1-morphism $\lambda : t^* \tau \to t^* \tau$. 
But this means that we specify for each $x \in \Gamma$ a 2-linear map $\lambda : \FinVect \to \FinVect$, i.e.\  -- by evaluation of $\lambda$ on the ground field -- a vector space $\varrho_\lambda(x)$. For a morphism $x \to y$ in $\Gamma$ we obtain a natural transformation $\lambda_x \to \lambda_y$, which is equivalent to a linear map $\varrho_\lambda(x) \to \varrho_\lambda(y)$. This shows that the objects of $\End_{\TRepGrpd}(\star)$ can be identified with vector bundles over essentially finite groupoids. 

A 1-morphism $(\Gamma_0,\lambda_0) \to (\Gamma_1,\lambda_1)$ in $\End_{\TRepGrpd}(\star)$ is a span of spans 
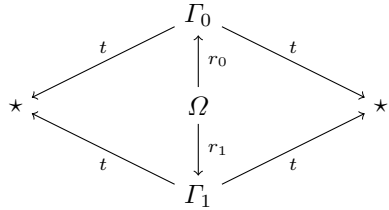
\begin{center}
	\begin{tikzpicture}[scale=1.2,     implies/.style={double,double equal sign distance,-implies},
		dot/.style={shape=circle,fill=black,minimum size=2pt,
			inner sep=0pt,outer sep=2pt},]
		\node (A1) at (0,0) {$\star$};
		\node (A2) at (2,1) {$\Gamma_0$};
		\node (A3) at (4,0) {$\star$};
		\node (B2) at (2,-1) {$\Gamma_1$};
		\node (C) at (2,0) {$\Omega$};
		\node (B1) at (1,-0.5) {$$};
		\node (B3) at (2,-1) {$$};
		
		\path[->,font=\scriptsize]
		(A2) edge node[above]{$t$} (A1)
		(A2) edge node[above]{$t$} (A3)
		(B2) edge node[below]{$t$} (A1)
		(B2) edge node[below]{$t$} (A3)
		(C) edge node[right]{$r_1$} (B2)
		(C) edge node[right]{$r_0$} (A2);
	\end{tikzpicture}
\end{center} together with a natural morphism $\omega$ between the intertwiners $t^* \lambda_0, t^* \lambda_1 : r_0^* t^* \tau \to r_1^* t^* \tau$, which is just an intertwiner from $r_0^* \varrho_{\lambda_0}$ to $r_1^* \varrho_{\lambda_1}$. All these identifications extend naturally to the composition and the symmetric monoidal structure.    
\end{proof}

\noindent In a symmetric monoidal bicategory there is a notion of a (fully) dualizable object, see \cite[Definition~2.3.5]{lurietft}: An object in a symmetric monoidal bicategory is called \emph{dualizable} if it is dualizable in the homotopy category. So informally speaking, a dualizable object has evaluation and coevaluation 1-morphisms which obey the triangle identities up to 2-isomorphism. The 2-isomorphisms are not required to be coherent. 

Exactly the same arguments which yielded duals in $\RepGrpd{}$ in \cite[Proposition~3.8]{schweigertwoikeofk} prove the following result:

\begin{proposition}\label{satztrepgrpddual}
Every object in $\TRepGrpd$ is dualizable.
\end{proposition}

\subsection{The parallel section functor\label{secparsectionfunctorconstruc}}
This subsection is devoted to the construction of the parallel section functor $\Par : \TRepGrpd \to \TwoVect$. It will send an object $(\Gamma,\varrho)$, i.e.\ a 2-vector bundle $\varrho$ over a groupoid $\Gamma$ to its space of parallel sections (Definition~\ref{defprallelsectionsof2vectorbundle}). We have already shown that this is indeed a 2-vector space (Proposition~\ref{satzactinterparsec}).
It remains to define $\Par$ on 1-morphisms and 2-morphisms in $\TRepGrpd$. This is accomplished in the following two definitions using the pullback and pushforward constructions from Section~\ref{secpullpushext}. 

\begin{definition}[Parallel section functor on 1-morphisms]\label{defparext1}
	Let $(\Gamma_0,\varrho_0) \stackrel{r_0}{\longleftarrow} (\Lambda,\lambda) \stackrel{r_1}{\to} (\Gamma_1,\varrho_1)$ be a 1-morphism in $\TRepGrpd$. Denote by $\Par(\Lambda, \lambda)$ the 2-linear map
	\begin{align}
		\Par (\Lambda,\lambda) : \Par \varrho_0 \xrightarrow{r_0^*} \Par r_0^* \varrho_0 \xrightarrow{\lambda_*} \Par r_1^* \varrho_1 \xrightarrow{{r_1}_*} \Par \varrho_1.
	\end{align} Here we use the 2-linear pull map from Proposition~\ref{satzpulloperationsofk}, the operation of intertwiners on parallel sections from Proposition~\ref{satzactinterparsec} and the 2-linear push map from Proposition~\ref{satztwolinpushmap}. 
\end{definition}

\begin{definition}[Parallel section functor on 2-morphisms]\label{defparext2}
	Let \begin{center}
		\begin{tikzpicture}[scale=2,     implies/.style={double,double equal sign distance,-implies},
			dot/.style={shape=circle,fill=black,minimum size=2pt,
				inner sep=0pt,outer sep=2pt},]
			\node (A1) at (0,0) {$(\Gamma_0,\varrho_0)$};
			\node (A2) at (2,1) {$(\Lambda,\lambda)$};
			\node (A3) at (4,0) {$(\Gamma_1,\varrho_1)$};
			\node (B2) at (2,-1) {$(\Lambda',\lambda')$};
			\node (C) at (2,0) {$(\Omega,\omega)$};
			\node (B1) at (1,-0.5) {$$};
			\node (B3) at (2,-1) {$$};
			\path[->,font=\scriptsize]
			(A2) edge node[above]{$r_0$} (A1)
			(A2) edge node[above]{$r_1$} (A3)
			(B2) edge node[below]{$r_0'$} (A1)
			(B2) edge node[below]{$r_1'$} (A3)
			(C) edge node[right]{$t'$} (B2)
			(C) edge node[right]{$t$} (A2);
			\draw (1.5,0) edge[implies] node[above] {\scriptsize$\alpha_0$} (0.5,0);
			\draw (2.5,0) edge[implies] node[above] {\scriptsize$\alpha_1$} (3.5,0);
		\end{tikzpicture}
	\end{center} be a 2-morphism in $\TRepGrpd$. Then we define the 2-morphism
	\begin{align}
		\Par (\Omega,\omega) : \Par (\Lambda,\lambda) \to \Par (\Lambda',\lambda')
	\end{align} to be
	\begin{center}
		\begin{tikzpicture}[scale=1.5,     implies/.style={double,double equal sign distance,-implies},
			dot/.style={shape=circle,fill=black,minimum size=2pt,
				inner sep=0pt,outer sep=2pt},]
			\node (A1) at (0,1) {$\Par r_0^* \varrho_0$};
			\node (A2) at (4,1) {$\Par r_1^* \varrho_1$};
			\node (B1) at (0,0) {$\Par t^* r_0^* \varrho_0$};
			\node (B2) at (4,0) {$\Par t^* r_1^* \varrho_1$};
			\node (C1) at (0,-1) {$\Par {t'}^* {r_0'}^* \varrho_0$};
			\node (C2) at (4,-1) {$\Par {t'}^* {r_1'}^* \varrho_1$};
			\node (D1) at (0,-2) {$\Par {r_0'}^* \varrho_0$};
			\node (D2) at (4,-2) {$\Par {r_1'}^* \varrho_1$};
			\node (E0) at (-3,-0.5) {$\Par \varrho_0$};
			\node (E1) at (7,-0.5) {$\Par \varrho_1$};
			\node (F0) at (-1.5,-0.8) {$$};
			\node (G0) at (-1.2,-0.1) {$$};
			\node (G1) at (5.5,-0.1) {$$};
			\node (F1) at (5.2,-0.8) {$$};
			\path[->,font=\scriptsize]
			(A1) edge node[above]{$\lambda_*$} (A2)
			(E0) edge node[above]{$r_0^*$} (A1)
			(E0) edge node[below]{${r_0'}^*$} (D1)
			(A2) edge node[above]{${r_1}_*$} (E1)
			(D2) edge node[below]{${r_1'}_*$} (E1)
			(A1) edge node[left]{$t^*$} (B1)
			(B1) edge node[left]{$\varrho_0(\alpha_0)_*$} (C1)
			(B2) edge node[right]{$\varrho_1(\alpha_1)_*$} (C2)
			(C1) edge node[above]{${t'}^* \lambda'_*$} (C2)
			(D1) edge node[below]{$\lambda_*'$} (D2)
			(B1) edge node[above]{$t^* \lambda_*$} (B2)
			(A2) edge node[right]{$t^*$} (B2)
			(C1) edge node[left]{$t'_*$} (D1)
			(C2) edge node[right]{$t'_*$} (D2);
			\draw (B2) edge[implies] node[above] {\scriptsize$\omega_*\ $} (C1);
			\draw (C2) edge[implies] node[above] {\scriptsize$\cong$} (D1);
			\draw (G0) edge[implies] node[left] {\scriptsize${\alpha_0}_*$} (F0);
			\draw (G1) edge[implies] node[right] {\scriptsize${\alpha_1^*}$} (F1);
		\end{tikzpicture},
	\end{center} where 
	
	\begin{itemize}
		
		\item the commutativity of the top square expresses the naturality of the pullback maps (Proposition~\ref{satzpulloperationsofk}, \ref{satzpulloperationsofkd}) and the commutativity of the lowest square expresses the naturality of the pushforward maps up to natural isomorphism (Proposition~\ref{satzmsceigpushmapsext}, \ref{satzmsceigpushmapsextb}),
		
		\item $\omega_*$ is the application of the functor from Proposition~\ref{satzactinterparsec} to $\omega$,
		
		\item ${\alpha_0}_* : t'_* \varrho_0(\alpha_0) (r_0 t)^* \to {r_0'}^*$ comes from the application of 
		Proposition~\ref{satz440transf} to the square
		\begin{center}
			\begin{tikzpicture}[scale=2, implies/.style={double,double equal sign distance,-implies},
				dot/.style={shape=circle,fill=black,minimum size=2pt,
					inner sep=0pt,outer sep=2pt},]
				\node (A1) at (0,1) {$\Omega$};
				\node (A2) at (1,1) {$\Gamma_0$};
				\node (B1) at (0,0) {$\Lambda'$};
				\node (B2) at (1,0) {$\Gamma_0$};
				\path[->,font=\scriptsize]
				(A1) edge node[above]{$r_0t$} (A2)
				(A1) edge node[left]{$t'$} (B1)
				(A2) edge node[right]{$\id_{\Gamma_0}$} (B2)
				(B1) edge node[below]{$r_0'$} (B2);
				\draw (A2) edge[implies] node[above] {\scriptsize$\alpha_0\ \  $} (B1);
			\end{tikzpicture},
		\end{center} 
		
		\item and $\alpha_1^* : {r_1}_* \to (r_1't')_* \varrho_1(\alpha_1) t^*$ comes similarly from the application of 
		Proposition~\ref{satz440transf} to the square
		\begin{center}
			\begin{tikzpicture}[scale=2, implies/.style={double,double equal sign distance,-implies},
				dot/.style={shape=circle,fill=black,minimum size=2pt,
					inner sep=0pt,outer sep=2pt},]
				\node (A1) at (0,1) {$\Omega$};
				\node (A2) at (1,1) {$\Lambda$};
				\node (B1) at (0,0) {$\Gamma_1$};
				\node (B2) at (1,0) {$\Gamma_1$};
				\path[->,font=\scriptsize]
				(A1) edge node[above]{$t$} (A2)
				(A1) edge node[left]{$r_1't'$} (B1)
				(A2) edge node[right]{$r_1$} (B2)
				(B1) edge node[below]{$\id_{\Gamma_1}$} (B2);
				\draw (A2) edge[implies] node[above] {\scriptsize$\alpha_1\ \  $} (B1);
			\end{tikzpicture}.
		\end{center} 
	\end{itemize}
\end{definition}

\noindent We are now ready to formulate the main result of this paper whose proof will take the rest of the paper, except for Proposition~\ref{satzresofparfunctor}.

\begin{theorem}[Parallel section functor]\label{thmextparsecfunctor}
	The assignments of Definition~\ref{defparext1} for 1-morphisms and Definition~\ref{defparext2} for 2-morphisms extend to a symmetric monoidal 2-functor
	\begin{align}
		\Par : \TRepGrpd \to \TwoVect,
	\end{align} that we call parallel section functor. 
\end{theorem}

\begin{proof}

\begin{enumerate}
	
	\item First we prove that $\Par$ is functorial on 1-morphisms. Obviously, it respects identities up to natural isomorphism. 
	For the proof of the compatibility with composition we start with two composable 1-morphisms
	\begin{align}
		(\Gamma_0,\varrho_0) \stackrel{r_0}{\longleftarrow} (\Lambda,\lambda) \stackrel{r_1}{\to} (\Gamma_1,\varrho_1) 
	\end{align} and
	\begin{align} (\Gamma_1,\varrho_1) \stackrel{r_1'}{\longleftarrow} (\Lambda',\lambda') \stackrel{r_2'}{\to} (\Gamma_2,\varrho_2) \end{align} in $\TRepGrpd$. 
	According to the definition of the composition of 1-morphisms (Definition~\ref{defTRepGrpd}, \ref{defTRepGrpdc1}) we need to form the homotopy pullback
	\begin{center}
		\begin{tikzpicture}[scale=0.8, implies/.style={double,double equal sign distance,-implies},
			dot/.style={shape=circle,fill=black,minimum size=2pt,
				inner sep=0pt,outer sep=2pt},]
			\node (A1) at (0,0) {$\Gamma_0$};
			\node (A2) at (2,2) {$\Lambda$};
			\node (A3) at (4,0) {$\Gamma_1$};
			\node (C) at (4,0.2) {$$};
			\node (A4) at (6,2) {$\Lambda'$};
			\node (A5) at (8,0) {$\Gamma_2$};
			\node (B2) at (4,4) {$\Lambda \times_\Omega \Lambda'$};
			\path[->,font=\scriptsize]
			(A2) edge node[above]{$r_0$} (A1)
			(A2) edge node[above]{$r_1$} (A3)
			(A4) edge node[above]{$r_1'$} (A3)
			(A4) edge node[above]{$r_2'$} (A5)
			(B2) edge node[left]{$\pi$} (A2)
			(B2) edge node[right]{$\ \pi'$} (A4);
			\draw (A2) edge[implies] node[above] {\scriptsize$\eta$} (A4);
		\end{tikzpicture}.
	\end{center}
	By the definition of the parallel section functor we find in our situation the following isomorphisms
	\begin{align} 
		\Par ((\Lambda',\lambda')\circ (\Lambda,\lambda)) &=  (r_2' \pi')_*({\pi'}^* \lambda')_*\varrho_1(\eta)_*(\pi^* \lambda)_*(r_0 \pi)^*\\ &\cong {r_2'}_*  \pi'_* ({\pi'}^* \lambda')_*\varrho_1(\eta)_*(\pi^* \lambda)_*\pi^* r_0^* \quad \left(\begin{array}{c}\text{Proposition~\ref{satzpulloperationsofk}, \ref{satzpulloperationsofkb} and} \\ \text{Proposition~\ref{satzmsceigpushmapsext}, \ref{satzmsceigpushmapsexta}} \end{array}\right) \\ &= {r_2'}_*  \pi'_* ({\pi'}^* \lambda')_*\varrho_1(\eta)_*\pi^* \lambda_* r_0^* \quad \text{(Proposition~\ref{satzpulloperationsofk}, \ref{satzpulloperationsofkd})} \\
		&\cong {r_2'}_*  \lambda'_* \pi'_* \varrho_1(\eta)_*\pi^* \lambda_* r_0^* \quad \text{(Proposition~\ref{satzmsceigpushmapsext}, \ref{satzmsceigpushmapsextb})} \\ &\cong {r_2'}_*  \lambda'_* {r_1'}^* {r_1}_* \lambda_* r_0^* \quad \text{(Proposition~\ref{mscextversionvonsatz440ausbscwoikeext})}\\ &= \Par (\Lambda',\lambda') \circ \Par (\Lambda,\lambda).
	\end{align} Their composition defines an isomorphism that we take as part of data of the 2-functor $\Par$. 
	
	\item Now we prove that the vertical composition of 2-morphisms is preserved. Again, it is obvious that identities are respected. For the proof of the composition law we take 2-morphisms
	\begin{center}
		\begin{tikzpicture}[scale=2,     implies/.style={double,double equal sign distance,-implies},
			dot/.style={shape=circle,fill=black,minimum size=2pt,
				inner sep=0pt,outer sep=2pt},]
			\node (A1) at (0,0) {$(\Gamma_0,\varrho_0)$};
			\node (A2) at (2,1) {$(\Lambda,\lambda)$};
			\node (A3) at (4,0) {$(\Gamma_1,\varrho_1)$};
			\node (B2) at (2,-1) {$(\Lambda',\lambda')$};
			\node (C) at (2,0) {$(\Omega,\omega)$};
			\node (B1) at (1,-0.5) {$$};
			\node (B3) at (2,-1) {$$};
			\path[->,font=\scriptsize]
			(A2) edge node[above]{$r_0$} (A1)
			(A2) edge node[above]{$r_1$} (A3)
			(B2) edge node[below]{$r_0'$} (A1)
			(B2) edge node[below]{$r_1'$} (A3)
			(C) edge node[right]{$t'$} (B2)
			(C) edge node[right]{$t$} (A2);
			\draw (1.5,0) edge[implies] node[above] {\scriptsize$\alpha_0$} (0.5,0);
			\draw (2.5,0) edge[implies] node[above] {\scriptsize$\alpha_1$} (3.5,0);
		\end{tikzpicture}
	\end{center} and 
	\begin{center}
		\begin{tikzpicture}[scale=2,     implies/.style={double,double equal sign distance,-implies},
			dot/.style={shape=circle,fill=black,minimum size=2pt,
				inner sep=0pt,outer sep=2pt},]
			\node (A1) at (0,0) {$(\Gamma_0,\varrho_0)$};
			\node (A2) at (2,1) {$(\Lambda',\lambda')$};
			\node (A3) at (4,0) {$(\Gamma_1,\varrho_1)$};
			\node (B2) at (2,-1) {$(\Lambda'',\lambda'')$};
			\node (C) at (2,0) {$(\widetilde{\Omega},\widetilde{\omega})$};
			\node (B1) at (1,-0.5) {$$};
			\node (B3) at (2,-1) {$$};
			\path[->,font=\scriptsize]
			(A2) edge node[above]{$r_0'$} (A1)
			(A2) edge node[above]{$r_1'$} (A3)
			(B2) edge node[below]{$r_0''$} (A1)
			(B2) edge node[below]{$r_1''$} (A3)
			(C) edge node[right]{$u''$} (B2)
			(C) edge node[right]{$u'$} (A2);
			\draw (1.5,0) edge[implies] node[above] {\scriptsize$\beta_0$} (0.5,0);
			\draw (2.5,0) edge[implies] node[above] {\scriptsize$\beta_1$} (3.5,0);
		\end{tikzpicture}
	\end{center} as well as the composition 
	\begin{center}
		\begin{tikzpicture}[scale=2,     implies/.style={double,double equal sign distance,-implies},
			dot/.style={shape=circle,fill=black,minimum size=2pt,
				inner sep=0pt,outer sep=2pt},]
			\node (A1) at (0,0) {$(\Gamma_0,\varrho_0)$};
			\node (A2) at (2,1) {$(\Lambda,\lambda)$};
			\node (A3) at (4,0) {$(\Gamma_1,\varrho_1)$};
			\node (B2) at (2,-1) {$(\Lambda'',\lambda'')$};
			\node (C) at (2,0) {$(\Omega \times_{\Lambda'} \widetilde \Omega,\omega \times_{\lambda'} \widetilde \omega)$};
			\node (B1) at (1,-0.5) {$$};
			\node (B3) at (2,-1) {$$};
			\path[->,font=\scriptsize]
			(A2) edge node[above]{$r_0$} (A1)
			(A2) edge node[above]{$r_1$} (A3)
			(B2) edge node[below]{$r_0''$} (A1)
			(B2) edge node[below]{$r_1''$} (A3)
			(C) edge node[right]{$v''$} (B2)
			(C) edge node[right]{$v$} (A2);
			\draw (1,0) edge[implies] node[above] {\scriptsize$\gamma_0$} (0.5,0);
			\draw (3,0) edge[implies] node[above] {\scriptsize$\gamma_1$} (3.5,0);
		\end{tikzpicture},
	\end{center}
	as given in Definition~\ref{defTRepGrpd}, \ref{defTRepGrpd2mor}. 
	According to Definition~\ref{defparext2} the natural transformation $\Par (\Omega \times_{\Lambda'} \widetilde \Omega,\omega \times_{\lambda'} \widetilde \omega)$ for the composition is given by
	\begin{center}
		\begin{tikzpicture}[scale=2, implies/.style={double,double equal sign distance,-implies},
			dot/.style={shape=circle,fill=black,minimum size=2pt,
				inner sep=0pt,outer sep=2pt},]
			\node (A0) at (0,0) {$\Par {r_0''}^* \varrho_0$};
			\node (A1) at (2,0) {$\Par {r_1''}^* \varrho_1$};
			\node (B0) at (0,1) {$\Par {v''}^* {r_0''}^* \varrho_0$};
			\node (B1) at (2,1) {$\Par {v''}^* {r_1''}^* \varrho_1$};
			\node (C0) at (0,2) {$\Par v^* r_0^* \varrho_0$};
			\node (C1) at (2,2) {$\Par v^* r_1^* \varrho_1$};
			\node (D0) at (0,3) {$\Par r_0^*\varrho_0$};
			\node (D1) at (2,3) {$\Par r_1^*\varrho_1$};
			\node (E0) at (-2,1.5) {$\Par \varrho_0$};
			\node (E1) at (4,1.5) {$\Par \varrho_1$};
			\node (F0) at (-0.7,2) {$$};
			\node (G0) at (-1.2,1.2) {$$};
			\node (F1) at (3.2,2) {$$};
			\node (G1) at (2.7,1.2) {$$};
			\path[->,font=\scriptsize]
			(D0) edge node[left]{$v^*$} (C0)
			(C0) edge node[left]{$\varrho_0(\gamma_0)_*$} (B0)
			(B0) edge node[left]{$v''_*$} (A0)
			(D1) edge node[right]{$v^*$} (C1)
			(C1) edge node[right]{$\varrho_1(\gamma_1)_*$} (B1)
			(B1) edge node[right]{$v''_*$} (A1)
			(D0) edge node[above]{$\lambda_*$} (D1)
			(C0) edge node[above]{$v^*\lambda_*$} (C1)
			(B0) edge node[above]{${v''}^*\lambda''_*$} (B1)
			(A0) edge node[above]{$\lambda_*''$} (A1)
			(E0) edge node[above]{$r_0^*$} (D0)
			(E0) edge node[below]{${r_0''}^* \ $} (A0)
			(D1) edge node[above]{$\ {r_1}_*$} (E1)
			(A1) edge node[below]{$\ {r_1''}_* \ $} (E1)
			;
			\draw (B1) edge[implies] node[above] {\scriptsize$\cong$} (A0);
			\draw (C1) edge[implies] node[above] {\scriptsize$(\omega\times_{\lambda'} \widetilde \omega)_* \quad \quad \quad$} (B0);
			\draw (F0) edge[implies] node[left] {\scriptsize${\gamma_0}_*$} (G0);
			\draw (F1) edge[implies] node[right] {\scriptsize$\gamma_1^*$} (G1);
		\end{tikzpicture}.
	\end{center} In a first step let us look at the inner ladder of this diagram. 
	The ladder is equal to \begin{center}
		\begin{tikzpicture}[scale=1.3, implies/.style={double,double equal sign distance,-implies},
			dot/.style={shape=circle,fill=black,minimum size=2pt,
				inner sep=0pt,outer sep=2pt},]
			\node (A0) at (0,0) {$\Par r_0^* \varrho_0$};
			\node (A1) at (3,0) {$\Par r_1^* \varrho_1$};
			\node (B0) at (0,-1) {$\Par t^* r_0^* \varrho_0 $};
			\node (B1) at (3,-1) {$\Par t^* r_1^* \varrho_1$};
			\node (C0) at (0,-2) {$\Par v^* r_0^* \varrho_0$};
			\node (C1) at (3,-2) {$\Par v^* r_1^* \varrho_1$};
			\node (D0) at (0,-3) {$\Par q^* {t'}^* {r_0'}^* \varrho_0$};
			\node (D1) at (3,-3) {$\Par q^* {t'}^* {r_1'}^* \varrho_1$};
			\node (E0) at (0,-4) {$\Par \widetilde q ^* {u'}^* {r_0'}^* \varrho_0$};
			\node (E1) at (3,-4) {$\Par \widetilde q^* {u'}^* {r_1'}^* \varrho_1$};
			\node (F0) at (0,-5) {$\Par \widetilde q^* {u''}^* {r_0''}^* \varrho_0$};
			\node (F1) at (3,-5) {$\Par \widetilde q^* {u''}^* {r_1''}^* \varrho_1$};
			\node (G0) at (0,-6) {$\Par {u''}^* {r_0''}^* \varrho_0$};
			\node (G1) at (3,-6) {$\Par {u''}^* {r_1''}^* \varrho_1$};
			\node (H0) at (0,-7) {$\Par {r_0''}^* \varrho_0$};
			\node (H1) at (3,-7) {$\Par {r_1''}^* \varrho_1$};
			\path[->,font=\scriptsize]
			(A0) edge node[left]{$t^*$} (B0)
			(B0) edge node[left]{$q^*$} (C0)
			(C0) edge node[left]{$\varrho_0(\alpha_0)_*$} (D0)
			(D0) edge node[left]{$\varrho_0(\eta)_*$} (E0)
			(E0) edge node[left]{$\varrho_0(\beta_0)_*$} (F0)
			(F0) edge node[left]{$\widetilde q_*$} (G0)
			(G0) edge node[left]{${u''_*}$} (H0)
			(A1) edge node[right]{$t^*$} (B1)
			(B1) edge node[right]{$q^*$} (C1)
			(C1) edge node[right]{$\varrho_1(\alpha_1)_*$} (D1)
			(D1) edge node[right]{$\varrho_1(\eta)_*$} (E1)
			(E1) edge node[right]{$\varrho_1(\beta_1)_*$} (F1)
			(F1) edge node[right]{$\widetilde q_*$} (G1)
			(G1) edge node[right]{${u''_*}$} (H1)
			(A0) edge node[above]{$\lambda_*$} (A1)
			(B0) edge node[above]{$t^* \lambda_*$} (B1)
			(C0) edge node[above]{$v^* \lambda_*$} (C1)
			(D0) edge node[above]{$q^* {t'}^* \lambda'_*$} (D1)
			(E0) edge node[above]{$\widetilde q^* {u'}^* \lambda'_*$} (E1)
			(F0) edge node[above]{$\widetilde q^* {u''}^* \lambda''_*$} (F1)
			(G0) edge node[above]{${u''}^* \lambda''_*$} (G1)
			(H0) edge node[above]{$\lambda''_*$} (H1);
			\draw (C1) edge[implies] node[above] {\scriptsize$(q^* \omega)_* \quad$} (D0);
			\draw (D1) edge[implies] node[above] {\scriptsize$\theta'_*$} (E0);
			\draw (E1) edge[implies] node[above] {\scriptsize$(\widetilde q^* \widetilde \omega)_* \ \ $} (F0);
			\draw (F1) edge[implies] node[above] {\scriptsize$\cong$} (G0);
			\draw (G1) edge[implies] node[above] {\scriptsize$\cong$} (H0);
		\end{tikzpicture}.
	\end{center}
	Here we have used the composition behaviour and naturality of the pull and push 1-morphisms (Proposition~\ref{satzpulloperationsofk} and Proposition~\ref{satzmsceigpushmapsext}), but we suppress the isomorphism $v_*'' \cong u_*'' \widetilde q_*$ for readability.  
	Additionally, we have unpacked the definition of $\omega\times_{\lambda'} \widetilde \omega$ (Definition~\ref{defTRepGrpd}, \ref{defTRepGrpd2mor}). Recall that the isomorphism $\theta'$ is the datum that $\lambda'$ comes equipped with.

	We investigate this ladder and obtain the equality \small   
	\begin{center}
		\begin{tikzpicture}[scale=1.3, implies/.style={double,double equal sign distance,-implies},
			dot/.style={shape=circle,fill=black,minimum size=2pt,
				inner sep=0pt,outer sep=2pt},]
			\node (A0) at (0,0) {$\Par r_0^* \varrho_0$};
			\node (A1) at (3,0) {$\Par r_1^* \varrho_1$};
			\node (B0) at (0,-1) {$\Par t^* r_0^* \varrho_0 $};
			\node (B1) at (3,-1) {$\Par t^* r_1^* \varrho_1$};
			\node (C0) at (0,-2) {$\Par v^* r_0^* \varrho_0$};
			\node (C1) at (3,-2) {$\Par v^* r_1^* \varrho_1$};
			\node (D0) at (0,-3) {$\Par q^* {t'}^* {r_0'}^* \varrho_0$};
			\node (D1) at (3,-3) {$\Par q^* {t'}^* {r_1'}^* \varrho_1$};
			\node (E0) at (0,-4) {$\Par \widetilde q ^* {u'}^* {r_0'}^* \varrho_0$};
			\node (E1) at (3,-4) {$\Par \widetilde q^* {u'}^* {r_1'}^* \varrho_1$};
			\node (F0) at (0,-5) {$\Par \widetilde q^* {u''}^* {r_0''}^* \varrho_0$};
			\node (F1) at (3,-5) {$\Par \widetilde q^* {u''}^* {r_1''}^* \varrho_1$};
			\node (G0) at (0,-6) {$\Par {u''}^* {r_0''}^* \varrho_0$};
			\node (G1) at (3,-6) {$\Par {u''}^* {r_1''}^* \varrho_1$};
			\node (H0) at (0,-7) {$\Par {r_0''}^* \varrho_0$};
			\node (H1) at (3,-7) {$\Par {r_1''}^* \varrho_1$};
			\path[->,font=\scriptsize]
			(A0) edge node[left]{$t^*$} (B0)
			(B0) edge node[left]{$q^*$} (C0)
			(C0) edge node[left]{$\varrho_0(\alpha_0)_*$} (D0)
			(D0) edge node[left]{$\varrho_0(\eta)_*$} (E0)
			(E0) edge node[left]{$\varrho_0(\beta_0)_*$} (F0)
			(F0) edge node[left]{$\widetilde q_*$} (G0)
			(G0) edge node[left]{${u''_*}$} (H0)
			(A1) edge node[right]{$t^*$} (B1)
			(B1) edge node[right]{$q^*$} (C1)
			(C1) edge node[right]{$\varrho_1(\alpha_1)_*$} (D1)
			(D1) edge node[right]{$\varrho_1(\eta)_*$} (E1)
			(E1) edge node[right]{$\varrho_1(\beta_1)_*$} (F1)
			(F1) edge node[right]{$\widetilde q_*$} (G1)
			(G1) edge node[right]{${u''_*}$} (H1)
			(A0) edge node[above]{$\lambda_*$} (A1)
			(B0) edge node[above]{$t^* \lambda_*$} (B1)
			(C0) edge node[above]{$v^* \lambda_*$} (C1)
			(D0) edge node[above]{$q^* {t'}^* \lambda'_*$} (D1)
			(E0) edge node[above]{$\widetilde q^* {u'}^* \lambda'_*$} (E1)
			(F0) edge node[above]{$\widetilde q^* {u''}^* \lambda''_*$} (F1)
			(G0) edge node[above]{${u''}^* \lambda''_*$} (G1)
			(H0) edge node[above]{$\lambda''_*$} (H1);
			\draw (C1) edge[implies] node[above] {\scriptsize$(q^* \omega)_* \quad$} (D0);
			\draw (D1) edge[implies] node[above] {\scriptsize$\theta'_*$} (E0);
			\draw (E1) edge[implies] node[above] {\scriptsize$(\widetilde q^* \widetilde \omega)_*\quad $} (F0);
			\draw (F1) edge[implies] node[above] {\scriptsize$\cong$} (G0);
			\draw (G1) edge[implies] node[above] {\scriptsize$\cong$} (H0);
		\end{tikzpicture}
		\begin{tikzpicture}[scale=1.3, implies/.style={double,double equal sign distance,-implies},
			dot/.style={shape=circle,fill=black,minimum size=2pt,
				inner sep=0pt,outer sep=2pt},]
			\node (A0) at (0,0) {$\Par r_0^* \varrho_0$};
			\node (A1) at (3,0) {$\Par r_1^* \varrho_1$};
			\node (B0) at (0,-1) {$\Par t^* r_0^* \varrho_0 $};
			\node (B1) at (3,-1) {$\Par t^* r_1^* \varrho_1$};
			\node (C0) at (0,-2) {$\Par {t'}^* {r_0'}^* \varrho_0$};
			\node (C1) at (3,-2) {$\Par {t'}^* {r_1'}^* \varrho_1$};
			\node (D0) at (0,-3) {$\Par q^* {t'}^* {r_0'}^* \varrho_0$};
			\node (D1) at (3,-3) {$\Par q^* {t'}^* {r_1'}^* \varrho_1$};
			\node (E0) at (0,-4) {$\Par \widetilde q ^* {u'}^* {r_0'}^* \varrho_0$};
			\node (E1) at (3,-4) {$\Par \widetilde q^* {u'}^* {r_1'}^* \varrho_1$};
			\node (F0) at (0,-5) {$\Par  {u'}^* {r_0'}^* \varrho_0$};
			\node (F1) at (3,-5) {$\Par  {u'}^* {r_1'}^* \varrho_1$};
			\node (G0) at (0,-6) {$\Par {u''}^* {r_0''}^* \varrho_0$};
			\node (G1) at (3,-6) {$\Par {u''}^* {r_1''}^* \varrho_1$};
			\node (H0) at (0,-7) {$\Par {r_0''}^* \varrho_0$};
			\node (H1) at (3,-7) {$\Par {r_1''}^* \varrho_1$};
			\node (EQ) at (-1.1,-3.5) {$=$};
			\path[->,font=\scriptsize]
			(A0) edge node[left]{$t^*$} (B0)
			(B0) edge node[left]{$\varrho_0(\alpha_0)_*$} (C0)
			(C0) edge node[left]{$q^*$} (D0)
			(D0) edge node[left]{$\varrho_0(\eta)$} (E0)
			(E0) edge node[left]{$\widetilde q_*$} (F0)
			(F0) edge node[left]{$\varrho_0(\beta_0)_*$} (G0)
			(G0) edge node[left]{${u''_*}$} (H0)
			(A1) edge node[right]{$t^*$} (B1)
			(B1) edge node[right]{$\varrho_1(\alpha_1)_*$} (C1)
			(C1) edge node[right]{$q^*$} (D1)
			(D1) edge node[right]{$\varrho_1(\eta)$} (E1)
			(E1) edge node[right]{$\widetilde q_*$} (F1)
			(F1) edge node[right]{$\varrho_1(\beta_1)_*$} (G1)
			(G1) edge node[right]{${u''_*}$} (H1)
			(A0) edge node[above]{$\lambda_*$} (A1)
			(B0) edge node[above]{$t^* \lambda_*$} (B1)
			(C0) edge node[above]{${t'}^* \lambda_*'$} (C1)
			(D0) edge node[above]{$q^* {t'}^* \lambda'_*$} (D1)
			(E0) edge node[above]{$\widetilde q^* {u'}^* \lambda'_*$} (E1)
			(F0) edge node[above]{${u'}^* \lambda'_*$} (F1)
			(G0) edge node[above]{${u''}^* \lambda''_*$} (G1)
			(H0) edge node[above]{$\lambda''_*$} (H1);
			\draw (B1) edge[implies] node[above] {\scriptsize$\omega_*$} (C0);
			\draw (D1) edge[implies] node[above] {\scriptsize$\theta'_*$} (E0);
			\draw (E1) edge[implies] node[above] {\scriptsize$\cong$} (F0);
			\draw (F1) edge[implies] node[above] {\scriptsize$\widetilde\omega_*$} (G0);
			\draw (G1) edge[implies] node[above] {\scriptsize$\cong$} (H0);
		\end{tikzpicture},
	\end{center} \normalsize where the changes only involve the second and the third as well as the fifth and the sixth square. By re-inserting the ladder we obtain\small 
	\begin{align}\begin{array}{c}
			\begin{tikzpicture}[scale=1.1, implies/.style={double,double equal sign distance,-implies},
				dot/.style={shape=circle,fill=black,minimum size=2pt,
					inner sep=0pt,outer sep=2pt},]
				\node (A0) at (0,0) {$\Par r_0^* \varrho_0$};
				\node (A1) at (3,0) {$\Par r_1^* \varrho_1$};
				\node (B0) at (0,-1) {$\Par t^* r_0^* \varrho_0 $};
				\node (B1) at (3,-1) {$\Par t^* r_1^* \varrho_1$};
				\node (C0) at (0,-2) {$\Par {t'}^* {r_0'}^* \varrho_0$};
				\node (C1) at (3,-2) {$\Par {t'}^* {r_1'}^* \varrho_1$};
				\node (D0) at (0,-3) {$\Par q^* {t'}^* {r_0'}^* \varrho_0$};
				\node (D1) at (3,-3) {$\Par q^* {t'}^* {r_1'}^* \varrho_1$};
				\node (E0) at (0,-4) {$\Par \widetilde q ^* {u'}^* {r_0'}^* \varrho_0$};
				\node (E1) at (3,-4) {$\Par \widetilde q^* {u'}^* {r_1'}^* \varrho_1$};
				\node (F0) at (0,-5) {$\Par  {u'}^* {r_0'}^* \varrho_0$};
				\node (F1) at (3,-5) {$\Par  {u'}^* {r_1'}^* \varrho_1$};
				\node (G0) at (0,-6) {$\Par {u''}^* {r_0''}^* \varrho_0$};
				\node (G1) at (3,-6) {$\Par {u''}^* {r_1''}^* \varrho_1$};
				\node (H0) at (0,-7) {$\Par {r_0''}^* \varrho_0$};
				\node (H1) at (3,-7) {$\Par {r_1''}^* \varrho_1$};
				\node (J0) at (-2,-3.5) {$\Par {r_0'}^* \varrho_0$};
				\node (J1) at (5,-3.5) {$\Par {r_1'}^* \varrho_1$};
				\node (K0) at (-4,-3.5) {$\Par  \varrho_0$};
				\node (K1) at (7,-3.5) {$\Par  \varrho_1$};
				\node (L0) at (-1,-3.5) {$\stackrel{\cong}{\Leftarrow}$};
				\node (L1) at (4,-3.5) {$\stackrel{\cong}{\Leftarrow}$};
				\node (M0) at (-2,-2.5) {$$};
				\node (M1) at (-1,-1.5) {$$};
				\node (N0) at (5.2,-4.3) {$$};
				\node (N1) at (4.2,-5.3) {$$};
				\node (O0) at (-1.2,-4.2) {$$};
				\node (O1) at (-2,-5) {$$};
				\node (P0) at (4.8,-1.7) {$$};
				\node (P1) at (4,-2.5) {$$};
				\path[->,font=\scriptsize]
				(A0) edge node[left]{$t^*$} (B0)
				(B0) edge node[left]{$\varrho_0(\alpha_0)_*$} (C0)
				(C0) edge node[left]{$q^*$} (D0)
				(D0) edge node[left]{$\varrho_0(\eta)_*$} (E0)
				(E0) edge node[left]{$\widetilde q_*$} (F0)
				(F0) edge node[left]{$\varrho_0(\beta_0)_*$} (G0)
				(G0) edge node[left]{${u''_*}$} (H0)
				(A1) edge node[right]{$t^*$} (B1)
				(B1) edge node[right]{$\varrho_1(\alpha_1)_*$} (C1)
				(C1) edge node[right]{$q^*$} (D1)
				(D1) edge node[right]{$\varrho_1(\eta)_*$} (E1)
				(E1) edge node[right]{$\widetilde q_*$} (F1)
				(F1) edge node[right]{$\varrho_1(\beta_1)_*$} (G1)
				(G1) edge node[right]{${u''_*}$} (H1)
				(A0) edge node[above]{$\lambda_*$} (A1)
				(B0) edge node[above]{$t^* \lambda_*$} (B1)
				(C0) edge node[above]{${t'}^* \lambda_*'$} (C1)
				(D0) edge node[above]{$q^* {t'}^* \lambda'_*$} (D1)
				(E0) edge node[above]{$\widetilde q^* {u'}^* \lambda'_*$} (E1)
				(F0) edge node[above]{$ {u'}^* \lambda'_*$} (F1)
				(G0) edge node[above]{${u''}^* \lambda''_*$} (G1)
				(H0) edge node[above]{$\lambda''_*$} (H1)
				(K0) edge node[above]{$r_0^*$} (A0)
				(K0) edge node[below]{${r_0''}^*\ $} (H0)
				(K0) edge node[above]{${r_0'}^*$} (J0)
				(C0) edge node[above]{$t_*'$} (J0)
				(J1) edge node[above]{${r_1'}_*$} (K1)
				(C1) edge node[above]{$t_*'$} (J1)
				(J0) edge node[below]{${u'}^*\ $} (F0)
				(J1) edge node[below]{${u'}^*\ $} (F1)
				(A1) edge node[above]{$\ {r_1}_*$} (K1)
				(H1) edge node[below]{$\ {r_1''}_*$} (K1);
				\draw (B1) edge[implies] node[above] {\scriptsize$\omega_*$} (C0);
				\draw (D1) edge[implies] node[above] {\scriptsize$\theta'_*$} (E0);
				\draw (E1) edge[implies] node[above] {\scriptsize$\cong$} (F0);
				\draw (F1) edge[implies] node[above] {\scriptsize$\widetilde\omega_*$} (G0);
				\draw (G1) edge[implies] node[above] {\scriptsize$\cong$} (H0);
				\draw (M1) edge[implies] node[left] {\scriptsize${\alpha_0}_*\ $} (M0);
				\draw (N0) edge[implies] node[right] {\scriptsize${\ \beta_1^*} $} (N1);
				\draw (O0) edge[implies] node[above] {\scriptsize${{\beta_0}_* \quad} $} (O1);
				\draw (P0) edge[implies] node[right] {\scriptsize${\ \alpha_1^* } $} (P1);
			\end{tikzpicture}.\end{array} \label{2parproof1}  
	\end{align}\normalsize
	Here we have also decomposed the triangles containing ${\gamma_0}_*$ and $\gamma_1^*$ into three smaller triangles each, which we will justify in step~\ref{proof2pardreiecke}. If we accept this for a moment, we can observe that the natural isomorphisms in the inner hexagon yield a natural isomorphism
	\begin{align} {u'}^* t_*' ({t'}^* \lambda'_*) \to ({u'}^* \lambda'_*){ u'}^* t_*' 
	\end{align} between 2-linear maps $\Par {t'}^* {r_0'}^* \varrho_0 \to \Par {u'}^* {r_1'}^* \varrho_1$. Evaluated on $s \in \Par {t'}^* {r_0'}^* \varrho_0$ and $\widetilde z \in \widetilde \Omega$ it consists of the isomorphism from 
	\begin{align}
		(({u'}^* t_*' ({t'}^* \lambda'_*))s)(\widetilde z)=\lim_{(z,g) \in {t'}^{-1}[u'(\widetilde z)]} \varrho_1(r_1'(g)) \lambda'_{t'(z)} s(z) \end{align} to \begin{align} (({u'}^* \lambda'_*) { u'}^* t_*')s)(\widetilde z)=\lim_{(z,g) \in {t'}^{-1}[u'(\widetilde z)]} \lambda'_{u'(\widetilde z)}  \varrho_0(r_0'(g)) s(z)
	\end{align}\normalsize described as follows: 
	We have to pull back the diagram underlying the first limit along the equivalence ${t'}^{-1}[u'(\widetilde z)]\cong \widetilde q^{-1}[\widetilde z]$ coming from the fiberwise characterization of the homotopy pullback to a diagram over $q^{-1}[\widetilde z]$, which amounts just to a change of variables. Then we apply the (pulled back version of) $\theta'$. Finally, we pull back the diagram to ${t'}^{-1}[u'(\widetilde z)]$ using (the inverse of) ${t'}^{-1}[u'(\widetilde z)]\cong \widetilde q^{-1}[\widetilde z]$. But this isomorphism is obviously equal to the isomorphism 
	\begin{align}
		\lim_{(z,g) \in {t'}^{-1}[u'(\widetilde z)]} \varrho_1(r_1'(g)) \lambda'_{t'(z)} s(z) \to \lim_{(z,g) \in {t'}^{-1}[u'(\widetilde z)]} \lambda'_{u'(\widetilde z)}  \varrho_0(r_0'(g)) s(z)
	\end{align} just coming from $\theta'$. 
	
	This allows us to simplify the inner hexagon in \eqref{2parproof1} and gives us
	\begin{center}
		\begin{tikzpicture}[scale=1.2, implies/.style={double,double equal sign distance,-implies},
			dot/.style={shape=circle,fill=black,minimum size=2pt,
				inner sep=0pt,outer sep=2pt},]
			\node (A0) at (0,0) {$\Par r_0^* \varrho_0$};
			\node (A1) at (3,0) {$\Par r_1^* \varrho_1$};
			\node (B0) at (0,-1) {$\Par t^* r_0^* \varrho_0 $};
			\node (B1) at (3,-1) {$\Par t^* r_1^* \varrho_1$};
			\node (C0) at (0,-2) {$\Par {t'}^* {r_0'}^* \varrho_0$};
			\node (C1) at (3,-2) {$\Par {t'}^* {r_1'}^* \varrho_1$};
			\node (D0) at (0,-3) {$\Par  {r_0'}^* \varrho_0$};
			\node (D1) at (3,-3) {$\Par  {r_1'}^* \varrho_1$};
			\node (E0) at (0,-4) {$\Par  {u'}^* {r_0'}^* \varrho_0$};
			\node (E1) at (3,-4) {$\Par {u'}^* {r_1'}^* \varrho_1$};
			\node (F0) at (0,-5) {$\Par  {u''}^* {r_0''}^* \varrho_0$};
			\node (F1) at (3,-5) {$\Par  {u''}^* {r_1''}^* \varrho_1$};
			\node (G0) at (0,-6) {$\Par  {r_0''}^* \varrho_0$};
			\node (G1) at (3,-6) {$\Par  {r_1''}^* \varrho_1$};
			\node (K0) at (-4,-3) {$\Par  \varrho_0$};
			\node (K1) at (7,-3) {$\Par  \varrho_1$};
			\node (M0) at (-2,-2.5) {$$};
			\node (M1) at (-1,-1.5) {$$};
			\node (N0) at (5.2,-3.3) {$$};
			\node (N1) at (4.2,-4.3) {$$};
			\node (O0) at (-1,-3.5) {$$};
			\node (O1) at (-2,-4.5) {$$};
			\node (P0) at (5,-1.7) {$$};
			\node (P1) at (4,-2.7) {$$};
			\path[->,font=\scriptsize]
			(A0) edge node[left]{$t^*$} (B0)
			(B0) edge node[left]{$\varrho_0(\alpha_0)_*$} (C0)
			(C0) edge node[left]{$t_*'$} (D0)
			(D0) edge node[left]{${u'}^*$} (E0)
			(E0) edge node[left]{$\varrho_0(\beta_0)_*$} (F0)
			(F0) edge node[left]{$u_*''$} (G0)
			(A1) edge node[right]{$t^*$} (B1)
			(B1) edge node[right]{$\varrho_1(\alpha_1)_*$} (C1)
			(C1) edge node[right]{$t_*'$} (D1)
			(D1) edge node[right]{${u'}^*$} (E1)
			(E1) edge node[right]{$\varrho_1(\beta_1)_*$} (F1)
			(F1) edge node[right]{$u_*''$} (G1)
			(A0) edge node[above]{$\lambda_*$} (A1)
			(B0) edge node[above]{$t^* \lambda_*$} (B1)
			(C0) edge node[above]{${t'}^* \lambda_*'$} (C1)
			(D0) edge node[above]{$\lambda'_*$} (D1)
			(E0) edge node[above]{${u'}^* \lambda'_*$} (E1)
			(F0) edge node[above]{$ {u''}^* \lambda''_*$} (F1)
			(G0) edge node[above]{${u''}^* \lambda''_*$} (G1)
			(K0) edge node[above]{$r_0^*$} (A0)
			(K0) edge node[above]{${r_0'}^*$} (D0)
			(D1) edge node[above]{${r_1'}_*$} (K1)
			(A1) edge node[above]{$\ {r_1}_*$} (K1)
			(G1) edge node[right]{${r_1''}^*$} (K1)
			(K0) edge node[left]{${r_0''}^*\ $} (G0);
			\draw (B1) edge[implies] node[above] {\scriptsize$\omega_*$} (C0);
			\draw (D1) edge[implies] node[above] {\scriptsize$\cong$} (E0);
			\draw (E1) edge[implies] node[above] {\scriptsize$\widetilde\omega_*$} (F0);
			\draw (F1) edge[implies] node[above] {\scriptsize$\cong$} (G0);
			\draw (M1) edge[implies] node[left] {\scriptsize${\alpha_0}_*\ $} (M0);
			\draw (N0) edge[implies] node[right] {\scriptsize${\ \beta_1^*} $} (N1);
			\draw (O0) edge[implies] node[above] {\scriptsize${{\beta_0}_* }\ \ \ $} (O1);
			\draw (P0) edge[implies] node[right] {\scriptsize${\ \alpha_1^* } $} (P1);
		\end{tikzpicture}.
	\end{center} Here we have replaced the inner hexagon by two squares. One of them commutes strictly (Proposition~\ref{satzpulloperationsofk}, \ref{satzpulloperationsofkd}), the other up to a natural isomorphism coming from $\theta'$ (Proposition~\ref{satzmsceigpushmapsext}, \ref{satzmsceigpushmapsextb}). This proves the preservation of the vertical composition. 
	
	\item \label{proof2pardreiecke}We still have to justify the decomposition of ${\gamma_0}_*$ and $\gamma_1^*$ that we have used to obtain \eqref{2parproof1}. 
	We only do this for ${\gamma_0}_*$ because it is the more difficult case (involving pushforward maps instead of only pullback maps).
	First note that the small inner triangles (the ones being part of the inner hexagon) come from a homotopy pullback, so the corresponding natural transformations are actually isomorphisms by Proposition~\ref{mscextversionvonsatz440ausbscwoikeext}. To prove that ${\gamma_0}_*$ is equal to the composition of the transformation living on the three triangles on the left side of \eqref{2parproof1}, we choose $s \in \Par \varrho_0$ and $y'' \in \Lambda$. Now both transformations in question correspond to maps
	\begin{align} 	(u_*'' \varrho_0(\beta_0)_* \widetilde q_* \varrho_0(\eta)_* q^* \varrho(\alpha_0)_* t^* r_0^* s)(y'') \to s(r_0''(y'')).\label{proofverticalcompeqn}\end{align}
	Using the definition of $\gamma_0$ in Definition~\ref{defTRepGrpd}, \ref{defTRepGrpd2mor} we can identify \begin{align}	(u_*'' \varrho_0(\beta_0)_* \widetilde q_* \varrho_0(\eta)_* q^* \varrho(\alpha_0)_* t^* r_0^* s)(y'')\end{align} with $(v_*'' \varrho_0(\gamma_0)_* v^* r_0^* s)(y'')$. Hence, we will see the maps \eqref{proofverticalcompeqn} as maps
	\begin{align}
		(v_*'' \varrho_0(\gamma_0)_* v^* r_0^* s)(y'')\to s(r_0''(y'')).
	\end{align} Now the composition of the three triangles amounts to the composition 
	\begin{align} &(v_*'' \varrho_0(\gamma_0)_* v^* r_0^* s)(y'')\\  \xrightarrow{\eta_*}& (u_*'' \varrho_0(\beta_0)_* {u'}^* t'_* \varrho(\alpha_0)_* t^* r_0^* s)(y'')\\ \xrightarrow{{\alpha_0}_*} &(u_*'' \varrho_0(\beta_0)_* {u'}^* {r_0'}^* s)(y'')\\ \xrightarrow{{\beta_0}_*} &s(r_0''(y'')),\label{footnotesizeeqn}
	\end{align} \normalsize and we have to show that it is given by ${\gamma_0}_*$. To see this, observe that the object $(v_*'' \varrho_0(\gamma_0)_* v^* r_0^* s)(y'')$ is a limit over the groupoid ${v''}^{-1} [y'']$, whereas \begin{align} (u_*'' \varrho_0(\beta_0)_* {u'}^* t'_* \varrho(\alpha_0)_* t^* r_0^* s)(y'')\end{align} is a limit over ${u''}^{-1}[y''] \times_{\Lambda'} \Omega$. The first map $\eta_*$ is the pushforward along the equivalence
	\begin{align}
		{v''}^{-1}[y''] = {(u''\circ \widetilde q)}^{-1}[y'']\cong {u''}^{-1}[y''] \times_{\widetilde \Omega} (\Omega \times_{\Lambda'} \widetilde \Omega) \cong {u''}^{-1}[y''] \times_{\Lambda'} \Omega.
	\end{align} Here, by pushforward we always mean pushforward of limits, i.e.\ pushforward of sections of ordinary vector bundles in the sense of Section~\ref{secrecallpullpush}. As a next step, $(u_*'' \varrho_0(\beta_0)_* {u'}^* {r_0'}^* s)(y'')$ is a limit over ${u''}^{-1}[y'']$ and ${\alpha_0}_*$ is the pushforward along the projection ${u''}^{-1}[y''] \times_{\Lambda'} \Omega \to {u''}^{-1}[y'']$. Finally, $s(r_0''(y''))$ is a limit over the terminal groupoid $\star$ and ${\beta_0}_*$ is the pushforward along the functor ${u''}^{-1}[y'']\to \star$. 
	By Proposition~\ref{satzcompositionslawspullpush}, \ref{satzcompositionslawspullpushb} we conclude that the composition \eqref{footnotesizeeqn} is the pushforward along the composition
	\begin{align}
		{v''}^{-1}[y''] \cong {u''}^{-1}[y''] \times_{\Lambda'} \Omega \to {u''}^{-1}[y''] \to \star
	\end{align} of functors, i.e.\ it integrates over the homotopy fiber ${v''}^{-1}[y'']$ with respect to groupoid cardinality. Hence, by Corollary~\ref{korbeckchevgrpdcard} it is equal to ${\gamma_0}_*$. 
	This gives us the missing step in the derivation of \eqref{2parproof1}.

	\item Next we prove that the horizontal composition of 2-morphisms is respected up to the isomorphisms specified for the composition of 1-morphisms. To this end, we take 2-morphisms
	\begin{center}
		\begin{tikzpicture}[scale=2,     implies/.style={double,double equal sign distance,-implies},
			dot/.style={shape=circle,fill=black,minimum size=2pt,
				inner sep=0pt,outer sep=2pt},]
			\node (A1) at (0,0) {$(\Gamma_0,\varrho_0)$};
			\node (A2) at (2,1) {$(\Lambda,\lambda)$};
			\node (A3) at (4,0) {$(\Gamma_1,\varrho_1)$};
			\node (B2) at (2,-1) {$(\Lambda',\lambda')$};
			\node (C) at (2,0) {$(\Omega,\omega)$};
			\node (B1) at (1,-0.5) {$$};
			\node (B3) at (2,-1) {$$};
			\path[->,font=\scriptsize]
			(A2) edge node[above]{$r_0$} (A1)
			(A2) edge node[above]{$r_1$} (A3)
			(B2) edge node[below]{$r_0'$} (A1)
			(B2) edge node[below]{$r_1'$} (A3)
			(C) edge node[right]{$t'$} (B2)
			(C) edge node[right]{$t$} (A2);
			\draw (1.5,0) edge[implies] node[above] {\scriptsize$\alpha_0$} (0.5,0);
			\draw (2.5,0) edge[implies] node[above] {\scriptsize$\alpha_1$} (3.5,0);
		\end{tikzpicture}
	\end{center} and 
	\begin{center}
		\begin{tikzpicture}[scale=2,     implies/.style={double,double equal sign distance,-implies},
			dot/.style={shape=circle,fill=black,minimum size=2pt,
				inner sep=0pt,outer sep=2pt},]
			\node (A1) at (0,0) {$(\Gamma_1,\varrho_1)$};
			\node (A2) at (2,1) {$(\Pi,\pi)$};
			\node (A3) at (4,0) {$(\Gamma_2,\varrho_2)$};
			\node (B2) at (2,-1) {$(\Pi',\pi')$};
			\node (C) at (2,0) {$(\widetilde{\Omega},\widetilde{\omega})$};
			\node (B1) at (1,-0.5) {$$};
			\node (B3) at (2,-1) {$$};
			\path[->,font=\scriptsize]
			(A2) edge node[above]{$v_1$} (A1)
			(A2) edge node[above]{$v_2$} (A3)
			(B2) edge node[below]{$v_1'$} (A1)
			(B2) edge node[below]{$v_2'$} (A3)
			(C) edge node[right]{$u'$} (B2)
			(C) edge node[right]{$u$} (A2);
			\draw (1.5,0) edge[implies] node[above] {\scriptsize$\beta_0$} (0.5,0);
			\draw (2.5,0) edge[implies] node[above] {\scriptsize$\beta_1$} (3.5,0);
		\end{tikzpicture}
	\end{center} and their horizontal composition
	\begin{center}
		\begin{tikzpicture}[scale=2,     implies/.style={double,double equal sign distance,-implies},
			dot/.style={shape=circle,fill=black,minimum size=2pt,
				inner sep=0pt,outer sep=2pt},]
			\node (A1) at (0,0) {$(\Gamma_0,\varrho_0)$};
			\node (A2) at (2,1) {$(\Lambda \times_{\Gamma_1} \Pi,\lambda \times_{\varrho_1} \pi)$};
			\node (A3) at (4,0) {$(\Gamma_2,\varrho_2)$};
			\node (B2) at (2,-1) {$(\Lambda' \times_{\Gamma_1} \Pi',\lambda' \times_{\varrho_1} \pi')$};
			\node (C) at (2,0) {$(\Omega \times_{\Gamma_1} \widetilde \Omega,\omega \times_{\varrho_1} \widetilde \omega)$};
			\node (B1) at (1,-0.5) {$$};
			\node (B3) at (2,-1) {$$};
			\path[->,font=\scriptsize]
			(A2) edge node[above]{$r_0p$} (A1)
			(A2) edge node[above]{$v_2q$} (A3)
			(B2) edge node[below]{$r_0'p'$} (A1)
			(B2) edge node[below]{$v_2'q'$} (A3)
			(C) edge node[right]{$d$} (B2)
			(C) edge node[right]{$c$} (A2);
			\draw (1,0) edge[implies] node[above] {\scriptsize$\delta_0$} (0.5,0);
			\draw (3,0) edge[implies] node[above] {\scriptsize$\delta_1$} (3.5,0);
		\end{tikzpicture}
	\end{center} as given in Definition~\ref{defTRepGrpd}, \ref{defTRepGrpd2morh}. We have to show the equality of natural transformations
	\begin{align}\begin{array}{c}
			\begin{tikzpicture}[scale=2, implies/.style={double,double equal sign distance,-implies},
				dot/.style={shape=circle,fill=black,minimum size=2pt,
					inner sep=0pt,outer sep=2pt},]
				\node (A0) at (-2,0) {$\Par (\Omega \times_{\Gamma_1} \widetilde \Omega,\omega \times_{\Gamma_1} \widetilde \omega)$};
				\node (EQ) at (-0.7,0) {$=$};
				\node (A1) at (0,0) {$\Par \varrho_0$};
				\node (A2) at (2,0) {$\Par \varrho_1$};
				\node (A3) at (4,0) {$\Par \varrho_2$};
				\node (C) at (1,0.5) {$$};
				\node (D) at (1,-0.5) {$$};
				\node (E) at (3,0.5) {$$};
				\node (F) at (3,-0.5) {$$};
				\node (G) at (2,1.2) {$$};
				\node (H) at (2,0.5) {$$};
				\node (I) at (2,-1.2) {$$};
				\node (J) at (2,-0.5) {$$};
				\draw[line width=0.5pt,->]
				(A1)   \dir{90}{90} node[above] {\scriptsize$\Par (\Lambda,\lambda)$} (A2) ;
				\draw[line width=0.5pt,->]
				(A1)   \dir{270}{270} node[below] {\scriptsize$\Par (\Lambda',\lambda')$} (A2) ;
				\draw[line width=0.5pt,->]
				(A2)   \dir{90}{90} node[above] {\scriptsize$\Par (\Pi,\pi)$} (A3) ;
				\draw[line width=0.5pt,->]
				(A2)   \dir{270}{270} node[below] {\scriptsize$\Par  (\Pi',\pi')$} (A3) ;
				\draw (C) edge[implies] node[right] {\scriptsize$\Par (\Omega,\omega)$} (D);
				\draw (E) edge[implies] node[right] {\scriptsize$\Par (\widetilde \Omega,\widetilde \omega)$} (F);
				\draw (G) edge[implies] node[right] {\scriptsize$\cong$} (H);
				\draw (I) edge[implies] node[right] {\scriptsize$\cong$} (J);
				\draw[line width=0.5pt,->]
				(A1)   \dir{90}{90} node[above] {\scriptsize$\Par (\Lambda\times_{\Gamma_1}\Pi,\lambda\times_{\Gamma_1} \pi)$} (A3) ;
				\draw[line width=0.5pt,->]
				(A1)   \dir{270}{270} node[below] {\scriptsize$\Par (\Lambda'\times_{\Gamma_1}\Pi',\lambda'\times_{\Gamma_1} \pi')$} (A3) ;
			\end{tikzpicture}.
		\end{array}  \label{2parproof2}  \end{align}
	We abbreviate the left hand side by $L$ and the right hand side by $R$.
	Using Definition~\ref{defTRepGrpd}, \ref{defTRepGrpd2morh} and the labels therein we find for the left hand side
	\scriptsize
	\begin{center}
		\begin{tikzpicture}[scale=1.5,     implies/.style={double,double equal sign distance,-implies},
			dot/.style={shape=circle,fill=black,minimum size=2pt,
				inner sep=0pt,outer sep=2pt},]
			\node (A0) at (-2.5,-0.5) {$L$};
			\node (EQ) at (-2.2,-0.5) {$=$};
			\node (A1) at (0,1) {$\Par (r_0p)^* \varrho_0$};
			\node (A2) at (6,1) {$\Par (v_2q)^* \varrho_2$};
			\node (A11) at (2,1) {$\Par (r_1p)^*\varrho_1$};
			\node (A12) at (4,1) {$\Par (r_1q)^* \varrho_1$};
			\node (B1) at (0,0) {$\Par (r_0pc)^* \varrho_0$};
			\node (B2) at (6,0) {$\Par (v_2qc)^* \varrho_2$};
			\node (B11) at (2,0) {$\Par (r_1pc)^*\varrho_1$};
			\node (B12) at (4,0) {$\Par (v_1pc)^* \varrho_1$};
			\node (C1) at (0,-1) {$\Par (r_0'pd)^* \varrho_0$};
			\node (C2) at (6,-1) {$\Par (v_2'q'd)^* \varrho_1$};
			\node (C11) at (2,-1) {$\Par (r_1'p'd)^*\varrho_1$};
			\node (C12) at (4,-1) {$\Par (v_1'q'd)^* \varrho_1$};
			\node (D1) at (0,-2) {$\Par (r_0'p')^* \varrho_0$};
			\node (D2) at (6,-2) {$\Par (v_2'q')^* \varrho_1$};
			\node (D11) at (2,-2) {$\Par (r_1'p')^*\varrho_1$};
			\node (D12) at (4,-2) {$\Par (v_1'q')^* \varrho_1$};
			\node (E0) at (-1.7,-0.5) {$\Par \varrho_0$};
			\node (E1) at (7.7,-0.5) {$\Par \varrho_2$};
			\node (F0) at (-1.1,-0.8) {$$};
			\node (G0) at (-0.6,-0.1) {$$};
			\node (G1) at (7.2,-0.1) {$$};
			\node (F1) at (6.7,-0.8) {$$};
			\path[->,font=\tiny]
			(A1) edge node[left]{$c^*$} (B1)
			(B1) edge node[left]{$\varrho_0(\ell^* \alpha_0)_*$} (C1)
			(C1) edge node[left]{$d_*$} (D1)
			(A11) edge node[left]{$c^*$} (B11)
			(B11) edge node[right]{$\varrho_1(\ell^* \alpha_1)_*$} (C11)
			(C11) edge node[left]{$d_*$} (D11)
			(A12) edge node[left]{$c^*$} (B12)
			(B12) edge node[left]{$\varrho_1(\widetilde \ell^* \beta_0)_*$} (C12)
			(C12) edge node[left]{$d_*$} (D12)
			(A2) edge node[left]{$c^*$} (B2)
			(B2) edge node[left]{$\varrho_2(\widetilde \ell^* \beta_1)_*$} (C2)
			(C2) edge node[left]{$d_*$} (D2)
			(A1) edge node[above]{$p^*\lambda_*$} (A11)
			(A11) edge node[above]{$\varrho_1(\eta)$} (A12)
			(A12) edge node[above]{$q^* \pi_*$} (A2)
			(B1) edge node[above]{$(pc)^*\lambda_*$} (B11)
			(B11) edge node[above]{$\varrho_1(\varepsilon)$} (B12)
			(B12) edge node[above]{$(qc)^* \pi$} (B2)
			(C1) edge node[above]{$(p'd)^*\lambda_*'$} (C11)
			(C11) edge node[above]{$\varrho_1(d^* \eta')$} (C12)
			(C12) edge node[above]{$(q'd)^* \pi'_*$} (C2)
			(D1) edge node[above]{${p'}^*\lambda_*'$} (D11)
			(D11) edge node[above]{$\varrho_1(\eta')$} (D12)
			(D12) edge node[above]{${q'}^* \pi'_*$} (D2)
			(E0) edge node[above]{$(r_0p)^*\ \ $} (A1)
			(E0) edge node[below]{$(r_0'p')^* \qquad$} (D1)
			(A2) edge node[above]{$\qquad ({v_2q})_*$} (E1)
			(D2) edge node[below]{$\qquad (v_2'q')_* $} (E1);
			\draw (G0) edge[implies] node[left] {\scriptsize${\delta_0}_*$} (F0);
			\draw (G1) edge[implies] node[left] {\scriptsize${\delta_1^*}$} (F1);
			\draw (C11) edge[implies] node[above] {\scriptsize$\cong$} (D1);
			\draw (C2) edge[implies] node[above] {\scriptsize$\cong$} (D12);
			\draw (B11) edge[implies] node[above] {\scriptsize$\ell^*\omega_*\ \ $} (C1);
			\draw (B2) edge[implies] node[above] {\scriptsize$\widetilde\ell^*\widetilde\omega_*\ \ $} (C12);
		\end{tikzpicture}, 
	\end{center} 
	\normalsize
	while the right hand side of \eqref{2parproof2} is given by
	\scriptsize
	\begin{center}
		\begin{tikzpicture}[scale=1.5,     implies/.style={double,double equal sign distance,-implies},
			dot/.style={shape=circle,fill=black,minimum size=2pt,
				inner sep=0pt,outer sep=2pt},]
			\node (R) at (-2.4,-0.5) {$R$};
			\node (EQ) at (-2.1,-0.5) {$=$};
			\node (A1) at (0,1) {$\Par r_0^* \varrho_0$};
			\node (A2) at (6,1) {$\Par v_2^* \varrho_2$};
			\node (A11) at (2,1) {$\Par r_1^*\varrho_1$};
			\node (A12) at (4,1) {$\Par v_1^* \varrho_1$};
			\node (O1) at (0,2) {$\Par (r_0p)^* \varrho_0$};
			\node (O11) at (2,2) {$\Par (r_1p)^* \varrho_2$};
			\node (O12) at (4,2) {$\Par (v_1q)^*\varrho_1$};
			\node (O2) at (6,2) {$\Par (v_2q)^* \varrho_1$};
			\node (P1) at (0,-3) {$\Par (r_0'p')^* \varrho_0$};
			\node (P11) at (2,-3) {$\Par (r_1'p')^* \varrho_2$};
			\node (P12) at (4,-3) {$\Par (v_1'q')^*\varrho_1$};
			\node (P2) at (6,-3) {$\Par (v_2'q')^* \varrho_1$};
			\node (P) at (3,-0.5) {$\Par \varrho_1$};
			\node (B1) at (0,0) {$\Par (r_0t)^* \varrho_0$};
			\node (B2) at (6,0) {$\Par (v_2u)^* \varrho_2$};
			\node (B11) at (2,0) {$\Par (r_1t)^*\varrho_1$};
			\node (B12) at (4,0) {$\Par (v_1u)^* \varrho_1$};
			\node (C1) at (0,-1) {$\Par (r_0't')^* \varrho_0$};
			\node (C2) at (6,-1) {$\Par (v_2'u')^* \varrho_2$};
			\node (C11) at (2,-1) {$\Par (r_1't')^*\varrho_1$};
			\node (C12) at (4,-1) {$\Par (v_1'u')^* \varrho_1$};
			\node (D1) at (0,-2) {$\Par {r_0'}^* \varrho_0$};
			\node (D2) at (6,-2) {$\Par {v_2'}^* \varrho_2$};
			\node (D11) at (2,-2) {$\Par {r_1'}^*\varrho_1$};
			\node (D12) at (4,-2) {$\Par {v_1'}^* \varrho_1$};
			\node (E0) at (-1.7,-0.5) {$\Par \varrho_0$};
			\node (E1) at (7.7,-0.5) {$\Par \varrho_2$};
			\node (F0) at (-1,-0.8) {$$};
			\node (G0) at (-0.5,-0.1) {$$};
			\node (G1) at (7,-0.1) {$$};
			\node (F1) at (6.5,-0.8) {$$};
			\node (X1) at (2.5,-0.1) {$$};
			\node (X2) at (2.5,-0.9) {$$};
			\node (Y1) at (3.5,-0.1) {$$};
			\node (Y2) at (3.5,-0.9) {$$};
			\node (Z1) at (3,1.5) {$$};
			\node (Z2) at (3,0.5) {$$};
			\node (ZA) at (3,-1) {$$};
			\node (ZB) at (3,-2.5) {$$};
			\path[->,font=\tiny]
			(A1) edge node[left]{$p^*$} (O1)
			(A11) edge node[left]{$p^*$} (O11)
			(O12) edge node[left]{$q_*$} (A12)
			(O2) edge node[left]{$q_*$} (A2)
			(D1) edge node[right]{${p'}^*$} (P1)
			(D11) edge node[left]{${p'}^*$} (P11)
			(P12) edge node[left]{$q'_*$} (D12)
			(P2) edge node[left]{$q'_*$} (D2)
			(A1) edge node[left]{$t^*$} (B1)
			(B1) edge node[left]{$\varrho_0(\alpha_0)_*$} (C1)
			(C1) edge node[left]{$t'_*$} (D1)
			(A11) edge node[left]{$t^*$} (B11)
			(B11) edge node[left]{$\varrho_1(\alpha_1)_*$} (C11)
			(C11) edge node[left]{$t'_*$} (D11)
			(A12) edge node[left]{$t^*$} (B12)
			(B12) edge node[right]{$\varrho_1(\beta_0)_*$} (C12)
			(C12) edge node[left]{$t'_*$} (D12)
			(A2) edge node[left]{$t^*$} (B2)
			(B2) edge node[left]{$\varrho_2( \beta_1)_*$} (C2)
			(C2) edge node[left]{$t'_*$} (D2)
			(A1) edge node[above]{$\lambda_*$} (A11)
			(O11) edge node[above]{$\varrho_1(\eta)$} (O12)
			(O1) edge node[above]{$p^* \lambda_*$} (O11)
			(O12) edge node[above]{$q^* \pi_*$} (O2)
			(A12) edge node[above]{$\pi_*$} (A2)
			(B1) edge node[above]{$t^*\lambda_*$} (B11)
			(B12) edge node[above]{$u^* \pi_*$} (B2)
			(C1) edge node[above]{${t'}^*\lambda_*'$} (C11)
			(C12) edge node[above]{${u'}^* \pi'_*$} (C2)
			(D1) edge node[above]{$\lambda_*'$} (D11)
			(P11) edge node[above]{$\varrho_1(\eta')$} (P12)
			(P1) edge node[above]{${p'}^* \lambda'_*$} (P11)
			(P12) edge node[above]{${q'}^* \pi'_*$} (P2)
			(D12) edge node[above]{$\pi'_*$} (D2)
			
			(A11) edge node[right]{${r_1}_*$} (P)
			(P) edge node[left]{$v_1^*$} (A12)
			(D11) edge node[right]{${r_1'}_*$} (P)
			(P) edge node[left]{${v_1'}^*$} (D12)

			(E0) edge node[above]{$r_0^* $} (A1)
			(E0) edge node[left]{$(r_0p)^* $} (O1)
			(E0) edge node[left]{$(r_0'p')^* $} (P1)
			(E0) edge node[above]{$\ \ {r_0'}^*$} (D1)
			(A2) edge node[above]{$\ {v_2}_*$} (E1)
			(O2) edge node[right]{$(v_2q)_*$} (E1)
			(P2) edge node[right]{$(v_2' q')_*$} (E1)
			(D2) edge node[above]{${v_2'}_*\ \ $} (E1);
			\draw (G0) edge[implies] node[left] {\scriptsize${\alpha_0}_*$} (F0);
			\draw (Z1) edge[implies] node[left] {\scriptsize${\eta}_*$} (Z2);
			\draw (ZA) edge[implies] node[left] {\scriptsize$\eta^*$} (ZB);
			\draw (G1) edge[implies] node[left] {\scriptsize${\beta_1^*}$} (F1);
			\draw (X1) edge[implies] node[left] {\scriptsize${\alpha_1^*}$} (X2);
			\draw (Y1) edge[implies] node[right] {\scriptsize${\beta_0}_*$} (Y2);
			\draw (C11) edge[implies] node[above] {\scriptsize$\cong$} (D1);
			\draw (O2) edge[implies] node[above] {\scriptsize$\cong$} (A12);
			\draw (D2) edge[implies] node[above] {\scriptsize$\cong$} (P12);
			\draw (C2) edge[implies] node[above] {\scriptsize$\cong$} (D12);
			\draw (B11) edge[implies] node[above] {\scriptsize$\omega_*\ \ $} (C1);
			\draw (B2) edge[implies] node[above] {\scriptsize$\widetilde\omega_*\ \ $} (C12);
		\end{tikzpicture}.\end{center}
	\normalsize
	We now describe the 2-morphisms \begin{align} L, \quad  R: \quad (v_2q)_* (q^* \pi_*)\varrho_1(\eta) (p^* \lambda_*) (r_0p)^* \to (v_2'q')_* ({q'}^* \pi'_*) \varrho_1(\eta') ({p'}^* \lambda'_*) (r_0'p')^* \end{align} explicitly by chasing through the diagrams. We will look at the component for a parallel section $s \in \Par \varrho_0$. The image of $s$ under the transformations will be evaluated at $x_2 \in \Gamma_2$. In the following step-by-step description of both transformations some of the obvious isomorphisms will not be mentioned explicitly in order to not obscure the main ideas: 
	
	\underline{Description of $L$:} Using the canonical equivalences
	\begin{align} (v_2 q)^{-1}[x_2] \cong v_2^{-1}[x_2] \times_\Pi (\Lambda \times_{\Gamma_1} \Pi) \cong v_2^{-1}[x_2] \times_{\Gamma_1} \Lambda \end{align} we obtain\small  
	\begin{align}
		((v_2q)_* (q^* \pi_*)\varrho_1(\eta) (p^* \lambda_*) (r_0p)^* s)(x_2) = \lim_{\substack{v_2^{-1}[x_2] \times_{\Gamma_1} \Lambda: \\ \bar y \in \Pi , v_2(\bar y) \stackrel{\xi}{\cong} x_2\\ y \in \Lambda ,  r_1(y) \stackrel{\nu}{\cong} v_1(\bar y)}} \varrho_2(\xi) \pi_{\bar y} \varrho_1(\nu) \lambda_y s(r_0(y)).\label{parfunctorproof3}  
	\end{align}\normalsize The groupoid $v_2^{-1}[x_2] \times_{\Gamma_1} \Lambda$ is the index groupoid for the diagram that we need to compute the limit of. We have written the index groupoid below the limit symbol. After a double point we also listed all the dummy variables. We will use this notation in the sequel.
	
	The first transformation we have to apply is $\delta_1^*$. Since
	\begin{align} (v_2' q'd)^{-1}[x_2] \cong {v_2'}^{-1}[x_2] \times_{\Pi'} (\Omega \times_{\Gamma_1} \widetilde \Omega) \end{align} its target is
	\begin{align}
		\begin{array}{rl}
			&( ( v_2'q'd)_* \varrho_2(\widetilde \ell^* \beta_1) c^* (q^* \pi_*) \varrho_1(\eta) (p^* \lambda_*) (r_0p)^* s )(x_2) \\ = &  \displaystyle \lim_{\substack{  {v_2'}^{-1}[x_2] \times_{\Pi'} (\Omega \times_{\Gamma_1} \widetilde \Omega): \\ \bar y' \in \Pi' , v_2'(\bar y') \stackrel{\xi'}{\cong} x_2 \\ z \in \Omega , \widetilde z \in \widetilde \Omega , r_1t(z) \stackrel{\mu}{\cong} v_1u(\widetilde z) \\ \bar y' \stackrel{\kappa}{\cong} u'(\widetilde z)  }} \varrho_2(\xi'v_2'(\kappa)\beta_{1,\widetilde z}) \pi_{u(\widetilde z)} \varrho_1(\mu) \lambda_{t(z)} s(r_0t(z)).\end{array}\label{parfunctorproof4}  
	\end{align}\normalsize
	The needed map from \eqref{parfunctorproof3} to \eqref{parfunctorproof4} is the pullback along the functor
	\begin{align} {v_2'}^{-1}[x_2] \times_{\Pi'} (\Omega \times_{\Gamma_1} \widetilde \Omega) \to v_2^{-1}[x_2] \times_{\Gamma_1} \Lambda\end{align} which, on the level of dummy variables as established in \eqref{parfunctorproof3} and \eqref{parfunctorproof4} sends
	\begin{align}
		\left(\bar y',v_2'(\bar y') \stackrel{\xi'}{\cong} x_2,z,\widetilde z, r_1t(z) \stackrel{\mu}{\cong} v_1 u(\widetilde z), \bar y' \stackrel{\kappa}{\cong} u'(\widetilde z) \right) \in {v_2'}^{-1}[x_2] \times_{\Pi'} (\Omega \times_{\Gamma_1} \widetilde \Omega) \end{align} to \begin{align} \left(u(\widetilde z), v_2u(\widetilde z) \stackrel{\beta_1}{\cong} v_2'u'(\widetilde z) \stackrel{\kappa}{\cong} v_2'(\bar y')\stackrel{\xi'}{\cong} x_2 , t(z), r_1t(z) \stackrel{\mu}{\cong} v_1 u(\widetilde z)\right) \in v_2^{-1}[x_2] \times_{\Gamma_1} \Lambda.\end{align} The next transformation does not change the index groupoids, but is the vertex-wise transformation
	\begin{align}
		\varrho_2(\xi'v_2'(\kappa))\varrho_2(\beta_{1,\widetilde z}) \pi_{u(\widetilde z)} \varrho_1(\mu) \lambda_{t(z)} s(r_0t(z)) \xrightarrow{\widetilde \omega_{\widetilde z}} \varrho_2(\xi'v_2'(\kappa) )\pi'_{u'(\widetilde z)} \varrho_1(\beta_{0,z}\mu) \lambda_{t(z)} s(r_0t(z)) 
	\end{align} applied to the diagram on the right hand side of \eqref{parfunctorproof4}.
	In the next step we have to replace the groupoid $\Omega \times_{\Gamma_1} \widetilde \Omega$ as a part of the index groupoid in \eqref{parfunctorproof4} by the canonically equivalent groupoid $(\Omega \times_{\Gamma_1} \widetilde \Omega)'$, see
	Definition~\ref{defTRepGrpd}, \ref{defTRepGrpd2morh} for the notation. More concretely, we replace $z \in \Omega$ and $\widetilde z \in \widetilde \Omega$ together with $r_1t(z) \stackrel{\mu}{\cong} v_1u(\widetilde z)$ by the same pair $(z,z') \in \Omega \times \widetilde \Omega$, but now with $r_1't'(z) \stackrel{\mu'}{\cong} r_1' t'(z)$, where $\mu' \alpha_{1,z} = \beta_{0,\widetilde z} \mu$. This leaves us with vertices
	\begin{align}\varrho_2(\xi'v_2'(\kappa) )  \pi'_{u'(\widetilde z)} \varrho_1(\mu' \alpha_{1,z}) \lambda_{t(z)} s(r_0t(z)) \end{align} such that we can apply the vertex-wise transformation\small  
	\begin{align}
		\varrho_2(\xi'v_2'(\kappa) ) \pi'_{u'(\widetilde z)} \varrho_1(\mu' \alpha_{1,z}) \lambda_{t(z)} s(r_0t(z)) &\xrightarrow{\omega_z} \varrho_2(\xi'v_2'(\kappa) ) \pi'_{u'(\widetilde z)} \varrho_1(\mu') \lambda_{t'(z)} \varrho_0(\alpha_{\alpha_0,z}) s(r_0t(z)) \\& \cong \varrho_2(\xi'v_2'(\kappa) ) \pi'_{u'(\widetilde z)} \varrho_1(\mu') \lambda_{t'(z)}  s(r_0't'(z)),
	\end{align} \normalsize where the last isomorphism comes from parallelity of $s$. We need to perform one last step, namely the application of ${\delta_0}_*$: Our current index groupoid is ${v_2'}^{-1}[x_2] \times_{\Pi'} (\Omega \times_{\Gamma_1} \widetilde \Omega)'$, whereas the index groupoid for the final result \begin{align} ((v_2'q')_* ({q'}^* \pi'_*) \varrho_1(\eta') ({p'}^* \lambda'_*) (r_0'p')^* s)(x_2)\end{align} is ${v_2'}^{-1}[x_2] \times_{\Gamma_1} \Lambda'$. The needed map comes from pushing along the functor
	\begin{align}
		{v_2'}^{-1}[x_2] \times_{\Pi'} (\Omega \times_{\Gamma_1} \widetilde \Omega)' \to {v_2'}^{-1}[x_2] \times_{\Gamma_1} \Lambda',
	\end{align} which is induced from projection to $\Omega$ and $t' : \Omega \to \Lambda'$. In summary, the natural transformation $L=\Par (\Omega \times_{\Gamma_1} \widetilde \Omega,\omega \times_{\Gamma_1} \widetilde \omega)$ consists of the maps
	\begin{align}
		((v_2q)_* (q^* \pi_*)\varrho_1(\eta) (p^* \lambda_*) (r_0p)^* s)(x_2) \to ((v_2'q')_* ({q'}^* \pi'_*) \varrho_1(\eta') ({p'}^* \lambda'_*) (r_0'p')^* s)(x_2)
	\end{align} obtained by performing two operations:
	\begin{itemize}
		\item Apply $\omega$ and $\widetilde \omega$ vertex-wise to the diagrams involved.
		
		\item Compute on the level of the index groupoids the pull-push map along the span\small 
		\begin{align}
			v_2^{-1}[x_2] \times_{\Gamma_1} \Lambda\longleftarrow \left\{\begin{array}{c} {v_2'}^{-1}[x_2] \times_{\Pi'} (\Omega \times_{\Gamma_1} \widetilde \Omega)\\ \cong {v_2'}^{-1}[x_2] \times_{\Pi'} (\Omega \times_{\Gamma_1} \widetilde \Omega)'\end{array}\right\} \to {v_2'}^{-1}[x_2] \times_{\Gamma_1} \Lambda'. \label{parfunctorproof5}  
		\end{align}\normalsize
		
	\end{itemize}
	These two operations obviously commute.

	\underline{Description of $R$:} The maps
	\begin{align}
		((v_2q)_* (q^* \pi_*)\varrho_1(\eta) (p^* \lambda_*) (r_0p)^* s)(x_2) \to ((v_2'q')_* ({q'}^* \pi'_*) \varrho_1(\eta') ({p'}^* \lambda'_*) (r_0' p')^* s)(x_2)
	\end{align} that $R$ consists of can be described similarly as for $L$. Since no new ideas enter, we just give the result. Again, we have to perform two commuting operations:
	\begin{itemize}
		\item Apply $\omega$ and $\widetilde \omega$ vertex-wise to the diagrams involved.
		
		\item Compute on the level of the index groupoids the pull-push map along the two composable spans
		\begin{align}
			{v_2}^{-1}[x_2] \times_{\Gamma_1} \Lambda \longleftarrow  {v_2'}^{-1}[x_2] \times_{\Pi'} (\widetilde \Omega \times_{\Gamma_1} \Lambda) \to {v_2'}^{-1}[x_2] \times_{\Gamma_1} \Lambda  
		\end{align}and
		\begin{align}
			{v_2'}^{-1}[x_2] \times_{\Gamma_1} \Lambda \longleftarrow  {v_2'}^{-1}[x_2] \times_{\Gamma_1} \Omega  \longleftarrow   {v_2'}^{-1}[x_2] \times_{\Gamma_1} \Lambda'.
		\end{align}
	\end{itemize}
	The composition of these two spans is (equivalent to) 
	the span in \eqref{parfunctorproof5}. 
	Indeed, we find the canonical equivalences
	\begin{align}
		&	\left({v_2'}^{-1} [x_2] \times_{\Pi'} (\widetilde \Omega \times_{\Gamma_1} \Lambda)\right) \times _{{v_2'}^{-1} [x_2] \times _{\Gamma_1} \Lambda} \left({v_2'}^{-1} [x_2] \times_{\Gamma_1} \Omega\right) \\ \cong &\left( \left(  {v_2'}^{-1} [x_2] \times_{\Gamma_1} \Lambda \right)  \times_{\Pi'} \widetilde \Omega \right)    \times _{{v_2'}^{-1} [x_2] \times _{\Gamma_1} \Lambda}   \left({v_2'}^{-1} [x_2] \times_{\Gamma_1} \Omega\right)\\
		\cong & \widetilde \Omega \times_{\Pi'}  \left({v_2'}^{-1} [x_2] \times_{\Gamma_1} \Omega\right) \\ \cong & {v_2'}^{-1} [x_2] \times_{\Pi'} (\Omega \times_{\Gamma_1} \widetilde \Omega)
	\end{align}
	Applying Proposition~\ref{satzparsecbeckchevalley} (equivariant Beck-Chevalley condition) now finishes the proof of \eqref{2parproof2}.

	\item We endow $\Par$ with a monoidal structure. For this we use for any 2-vector bundle $\varrho$ over $\Gamma$ and $\xi$ over $\Omega$ the obvious 2-linear maps
	\begin{align}
		\Phi : \Par \varrho \boxtimes \Par \xi \to \Par (\varrho \boxtimes \xi)
	\end{align} defined using the universal property of the Deligne product. Note that we suppress the groupoids in the notation, i.e.\ we use the shorthand $\Par \varrho = \Par (\Gamma,\varrho)$ etc. It remains to show that these 2-linear maps are equivalences. For the proof we can assume without loss of generality that $\Gamma$ and $\Omega$ are finite groups $G$ and $H$, in which case $\varrho$ and $\xi$ send the one object to a 2-vector space $\cat{V}$ and $\cat{W}$, respectively. 
	Now we use Proposition~\ref{satzthm35kir} and the notation used therein to write the relevant 2-vector spaces of parallel sections as
	\begin{align} \Par \varrho &\cong \bigoplus_{\mathcal{O} \in \mathscr{S}/G} \mathfrak{A}_{\alpha}(G,\mathcal{O})\text{-}\Mod,\\
		\Par \xi &\cong \bigoplus_{\mathcal{P} \in \mathscr{T}/H} \mathfrak{A}_\beta(H,\mathcal{P})\text{-}\Mod, \end{align} where $\mathscr{S}$ and $\mathscr{T}$ is the set of isomorphism classes of simple objects in $\cat{V}$ and $\cat{W}$ with representing systems $(X_s)_{s\in \mathscr{S}}$ and $(Y_t)_{t\in \mathscr{T}}$, respectively, and $\alpha \in H^2(G;\Map(\mathscr{S},\mathbb{C}^\times))$ and $\beta \in H^2(H;\Map(\mathscr{T},\mathbb{C}^\times))$ are the corresponding gerbes, see Section~\ref{parsecsecfor2vecbund}. 
	The set of isomorphism classes of $\cat{V}\boxtimes \cat{W}$ is given by $\mathscr{S} \times \mathscr{T}$ with a representing system $(X_s \boxtimes Y_t)_{(s,t) \in \mathscr{S} \times \mathscr{T}}$. Denote the corresponding gerbe by $\gamma \in H^2(G\times H; \Map(\mathscr{S}\times \mathscr{T},\mathbb{C}^\times))$. Using $(\mathscr{S}\times \mathscr{T}) // (G\times H)\cong \mathscr{S} //G \times \mathscr{T} //H$ 
	we obtain
	\begin{align} \Par (\varrho \boxtimes \xi) \cong \bigoplus_{(\mathcal{O},\mathcal{P}) \in \mathscr{S} /G \times \mathscr{T} /H} \mathfrak{A}_\gamma(G\times H; (\mathcal{O},\mathcal{P}))\text{-}\Mod .  \end{align} 
	This yields the weakly commutative diagram \footnotesize
	\begin{center}
		\begin{tikzpicture}[scale=2]
			\node (A1) at (0,1) {$\Par \varrho \boxtimes \Par \xi$};
			\node (A2) at (4,1) {$\Par (\varrho \boxtimes \xi)$};
			\node (B1) at (0,0) {$\left(  \bigoplus_{\mathcal{O} \in \mathscr{S}} \mathfrak{A}_\alpha(G,\mathcal{O})\text{-}\Mod   \right) \boxtimes \left(  \bigoplus_{\mathcal{P} \in \mathscr{T}} \mathfrak{A}_\beta(H,\mathcal{P})\text{-}\Mod \right)$};
			\node (B2) at (4,0) {$\bigoplus_{(\mathcal{O},\mathcal{P}) \in \mathscr{S} /G \times \mathscr{T} /H} \mathfrak{A}_\gamma(G\times H; (\mathcal{O},\mathcal{P}))\text{-}\Mod$};
			\path[->,font=\scriptsize]
			(A1) edge node[above]{$\Phi$} (A2)
			(A1) edge node[left]{$\cong$} (B1)
			(B1) edge node[above]{$\Psi$} (B2)
			(A2) edge node[right]{$\cong$} (B2);
		\end{tikzpicture},
	\end{center}\normalsize in which the vertical equivalences are the ones just discussed. The functor $\Psi$ admits the following description: The projections $\mathscr{S} \times \mathscr{T} \to \mathscr{S}$ and $\mathscr{S} \times \mathscr{T} \to \mathscr{T}$ yield maps
	\begin{align} 
		\Map (\mathscr{S} ,\mathbb{C}^\times) &\to \Map (\mathscr{S}\times \mathscr{T},\mathbb{C}^\times),\\
		\Map (\mathscr{T},\mathbb{C}^\times) &\to \Map (\mathscr{S}\times \mathscr{T},\mathbb{C}^\times) .\end{align} Together with the projections $G\times H \to G$ and $G\times H \to H$ they induce a map
	\begin{align}
		&	H^2(G;\Map (\mathscr{S} ,\mathbb{C}^\times)) \oplus H^2(H;\Map (\mathscr{T} ,\mathbb{C}^\times))\\ \to& H^2(G\times H;\Map (\mathscr{S}\times \mathscr{T} ,\mathbb{C}^\times)) \oplus H^2(G\times H;\Map (\mathscr{S} \times \mathscr{T} ,\mathbb{C}^\times)).
	\end{align}\normalsize
	Using the group operation in $H^2(G\times H;\Map (\mathscr{S}\times \mathscr{T} ,\mathbb{C}^\times))$ we obtain a map
	\begin{align}
		H^2(G;\Map (\mathscr{S} ,\mathbb{C}^\times)) \oplus H^2(H;\Map (\mathscr{T} ,\mathbb{C}^\times)) \to H^2(G\times H;\Map (\mathscr{S}\times \mathscr{T} ,\mathbb{C}^\times)) 
	\end{align} sending $(\alpha,\beta)$ to $\gamma$. With this observation in mind, $\Psi$ is the equivalence\footnotesize
	\begin{align}
		\left(  \bigoplus_{\mathcal{O} \in \mathscr{S}} \mathfrak{A}_\alpha(G,\mathcal{O})\text{-}\Mod   \right) \boxtimes \left(  \bigoplus_{\mathcal{P} \in \mathscr{T}} \mathfrak{A}_\beta(H,\mathcal{P})\text{-}\Mod \right) &\cong \bigoplus_{(\mathcal{O},\mathcal{P}) \in \mathscr{S} /G \times \mathscr{T} /H} \mathfrak{A}_\alpha(G,\mathcal{O})\text{-}\Mod \boxtimes \mathfrak{A}_\beta(H,\mathcal{P})\text{-}\Mod \\&\cong \bigoplus_{(\mathcal{O},\mathcal{P}) \in \mathscr{S} /G \times \mathscr{T} /H} \left( \mathfrak{A}_\alpha(G,\mathcal{O}) \otimes \mathfrak{A}_\beta(H,\mathcal{P}) \right)\text{-}\Mod \\ & \cong \bigoplus_{(\mathcal{O},\mathcal{P}) \in \mathscr{S} /G \times \mathscr{T} /H} \mathfrak{A}_\gamma(G\times H; (\mathcal{O},\mathcal{P}))\text{-}\Mod.
	\end{align} \normalsize Since $\Psi$ is an equivalence, so is $\Phi$. 
	
	To fully specify the monoidal structure we need to exhibit for 1-morphisms \begin{align}
		(\Gamma_0,\varrho_0) \stackrel{r_0}{\longleftarrow} (\Lambda,\lambda) \stackrel{r_1}{\to} (\Gamma_1,\varrho_1) 
	\end{align} and
	\begin{align} (\Gamma_0',\varrho_0') \stackrel{r_0'}{\longleftarrow} (\Lambda',\lambda') \stackrel{r_1'}{\to} (\Gamma_1',\varrho_1') \end{align} in $\TRepGrpd$
	natural 2-isomorphisms \begin{center}
		\begin{tikzpicture}[scale=2, implies/.style={double,double equal sign distance,-implies},
			dot/.style={shape=circle,fill=black,minimum size=2pt,
				inner sep=0pt,outer sep=2pt},]
			\node (A1) at (0,1) {$\Par (\Gamma_0,\varrho_0) \boxtimes \Par (\Gamma_0',\varrho_0') $};
			\node (A2) at (4,1) {$\Par (\Gamma_1,\varrho_1) \boxtimes \Par (\Gamma_1',\varrho_1') $};
			\node (B1) at (0,0) {$\Par ( (\Gamma_0,\varrho_0) \boxtimes (\Gamma_0',\varrho_0'))$};
			\node (B2) at (4,0) {$\Par ( (\Gamma_1',\varrho_1') \boxtimes (\Gamma_1',\varrho_1'))$};
			\path[->,font=\scriptsize]
			(A1) edge node[above]{$\Par (\Lambda,\lambda) \boxtimes \Par (\Lambda',\lambda')$} (A2)
			(A1) edge node[left]{$\Phi$} (B1)
			(A2) edge node[right]{$\Phi$} (B2)
			(B1) edge node[below]{$\Par ( (\Lambda,\lambda) \boxtimes (\Lambda,\lambda) )$} (B2);
			\draw (A2) edge[implies] node[above] {\scriptsize$\cong $} (B1);
		\end{tikzpicture}.
	\end{center} Evaluated on $s \in \Par (\Gamma_0,\varrho_0)$, $s' \in \Par (\Gamma_0,\varrho_0)$ and $(y,y') \in \Gamma_1 \times \Gamma_1'$ they are given by\footnotesize
	\begin{align}
		\left(	\left(\Phi \circ \left(\Par (\Lambda,\lambda) \boxtimes \Par (\Lambda',\lambda') \right) \right) (s\boxtimes s') \right)(y,y') &= (\Par (\Lambda,\lambda) s)(y) \boxtimes (\Par (\Lambda',\lambda') s')(y') \\ &= \lim_{(x,g) \in r_1^{-1}[y]} g. \lambda_x s(r_0(x)) \boxtimes \lim_{(x',g') \in r_1'^{-1}[y']} g'. \lambda_{x'} s'(r_0'(x')) \\ &\cong \lim_{ ((x,g), (x',g')) \in r_1^{-1}[y] \times r_1'^{-1}[y']} g. \lambda_x s(r_0(x)) \boxtimes g'. \lambda_{x'} s'(r_0'(x')) \\ &\cong \lim_{ ((x,g), (x',g')) \in r_1^{-1}[y] \times r_1'^{-1}[y']} (g.\lambda_x \boxtimes g'.\lambda_{x'}) (s(r_0(x)) \boxtimes s'(r_0'(x'))) \\ &\cong \lim_{ ((x,x'),(g,g')) \in (r_1 \times r_1')^{-1}[y,y']} (g.\lambda_x \boxtimes g'.\lambda_{x'}) (s(r_0(x)) \boxtimes s'(r_0'(x'))) \\ &= \left(((\Par ( (\Lambda,\lambda) \boxtimes (\Lambda,\lambda) ) \circ \Phi))(s\boxtimes s')\right)(y,y').
	\end{align} \normalsize Here we used that the Deligne product preserves limits (because $\boxtimes$-tensoring with an object is exact). This concludes the definition of the monoidal structure. This monoidal structure is also symmetric: For a monoidal functor between symmetric monoidal bicategories the symmetry is \emph{structure} and is given by natural isomorphisms  
	\begin{center}
		\begin{tikzpicture}[scale=2, implies/.style={double,double equal sign distance,-implies},
			dot/.style={shape=circle,fill=black,minimum size=2pt,
				inner sep=0pt,outer sep=2pt},]
			\node (A1) at (0,1) {$\Par \varrho \boxtimes \Par \xi $};
			\node (A2) at (4,1) {$\Par (\varrho \boxtimes \xi)$};
			\node (B1) at (0,0) {$\Par \xi \boxtimes\Par \varrho $};
			\node (B2) at (4,0) {$\Par (\xi \boxtimes \xi)$};
			\path[->,font=\scriptsize]
			(A1) edge node[above]{$\Phi$} (A2)
			(A1) edge node[left]{$c_{\Par \varrho , \Par \xi  }$} (B1)
			(A2) edge node[right]{$\Par c_{\varrho,\xi}$} (B2)
			(B1) edge node[below]{$\Phi$} (B2);
			\draw (A2) edge[implies] node[above] {\scriptsize$\cong $} (B1);
		\end{tikzpicture}
	\end{center} for all  2-vector bundles $\varrho$ over $\Gamma$ and $\xi$ over $\Omega$, where the horizontal maps are the monoidal structure of $\Par$ and the vertical maps are given by the braiding and the image thereof under the parallel section functor. 
	Indeed, we can easily exhibit the needed diagonal isomorphisms: For $(x,y) \in \Gamma \times \Omega$, $s \in \Par$ and $s' \in \Par \xi$ we recall that $	(   \Par c_{\varrho, \xi} \circ \Phi (s \boxtimes s'))(y,x)$ is computed as a limit over the homotopy fiber of the flip $\Gamma \times \Omega \to \Omega \times \Gamma$ over $(y,x)$. But this homotopy fiber is equivalent to the discrete groupoid with one point $(x,y)$, so we obtain a canonical isomorphism
	\begin{align}
		(   \Par c_{\varrho, \xi} \circ \Phi (s \boxtimes s'))(y,x) \cong s'(y) \boxtimes s(x) = (\Phi \circ c_{\Par \varrho , \Par \xi  } (s \boxtimes s')) (y,x) .
	\end{align}
	This isomorphism is natural in $x,y,s$ and $s'$ and gives us the symmetric structure.	
\end{enumerate}  
\end{proof}

\noindent The parallel section functors for 2-vector bundles generalizes the parallel section functor from \cite{trova} used in \cite{schweigertwoikeofk}:

\begin{proposition}\label{satzresofparfunctor}
Restriction of the parallel section functor $\Par : \TRepGrpd \to \TwoVect$ to the endomorphisms of the respective monoidal units yields the parallel section functor $\RepGrpd{} \to \FinVect$ from \cite[Theorem~3.17]{schweigertwoikeofk}, i.e.\ the square 
\begin{center}
	\begin{tikzpicture}[scale=2, implies/.style={double,double equal sign distance,-implies},
		dot/.style={shape=circle,fill=black,minimum size=2pt,
			inner sep=0pt,outer sep=2pt},]
		\node (A1) at (0,1) {$\End _{\TRepGrpd} (\star) $};
		\node (A2) at (3,1) {$\End_{\TwoVect}(\FinVect) $};
		\node (B1) at (0,0) {$\RepGrpd{}$};
		\node (B2) at (3,0) {$\FinVect$};
		\path[->,font=\scriptsize]
		(A1) edge node[above]{$\Par$} (A2)
		(A1) edge node[left]{$\Phi$} (B1)
		(A2) edge node[right]{$\Psi$} (B2)
		(B1) edge node[below]{$\Par$} (B2);
	\end{tikzpicture}
\end{center} featuring the equivalences from Proposition~\ref{satzendtrepgrpd} and Example~\ref{extwovectandalg2} is weakly commutative.
\end{proposition}

\begin{proof} The proof proceeds very much like the proof of Proposition~\ref{satzendtrepgrpd} and also uses the notation established therein. An object in $\End _{\TRepGrpd} (\star)$ is a span $\star \stackrel{t}{\longleftarrow} \Gamma \stackrel{t}{\to} \star$ together with an intertwiner $\lambda : t^* \tau \to t^* \tau$, where $\tau$ is the trivial representation of ther terminal groupoid $\star$ on $\FinVect$. The restriction
\begin{align}\Par :  \End _{\TRepGrpd} (\star) \to \End_{\TwoVect}(\FinVect) \end{align} of the parallel section functor sends this object to a 2-linear map $\FinVect \to \FinVect$. Under $\Psi$ this 2-linear map is identified with the vector space
\begin{align}
	(	t_* \lambda_* t^* s)(\star),
\end{align} where $s\in \Par(\star,\tau)$ is the parallel section sending $\star$ to $\mathbb{C}$.  But by definition
\begin{align}
	(	t_* \lambda_* t^* s)(\star) = \lim_{x\in \Gamma} \varrho_\lambda (x),\end{align} where the representation $\varrho_\lambda$ of $\Gamma$ is the image of $\Gamma$ and $\lambda$ under $\Phi$ (Proposition~\ref{satzendtrepgrpd}), and the limit of $\varrho_\lambda$ is just the space of parallel sections of $\varrho_\lambda$. 

Consider now a morphism
\begin{align}
	\begin{tikzpicture}[scale=2,     implies/.style={double,double equal sign distance,-implies},
		dot/.style={shape=circle,fill=black,minimum size=2pt,
			inner sep=0pt,outer sep=2pt},]
		\node (A1) at (0,0) {$(\star,\tau)$};
		\node (A2) at (2,1) {$(\Gamma_0,\lambda_0)$};
		\node (A3) at (4,0) {$(\star,\tau)$};
		\node (B2) at (2,-1) {$(\Gamma_1,\lambda_1)$};
		\node (C) at (2,0) {$(\Omega,\omega)$};
		\node (B1) at (1,-0.5) {$$};
		\node (B3) at (2,-1) {$$};
		\path[->,font=\scriptsize]
		(A2) edge node[above]{$t$} (A1)
		(A2) edge node[above]{$t$} (A3)
		(B2) edge node[below]{$t$} (A1)
		(B2) edge node[below]{$t$} (A3)
		(C) edge node[right]{$r_0$} (B2)
		(C) edge node[right]{$r_1$} (A2);
	\end{tikzpicture}
\end{align} in $\End _{\TRepGrpd} (\star)$. By what we have already seen in this proof the parallel section functor assigns to this morphism the transformation
\begin{align}\begin{array}{c}
		\begin{tikzpicture}[scale=1, implies/.style={double,double equal sign distance,-implies},
			dot/.style={shape=circle,fill=black,minimum size=2pt,
				inner sep=0pt,outer sep=2pt},]
			\draw[line width=0.5pt,->]
			(0,2) node {$\FinVect$}
			(4,2) node {$\FinVect$}
			(2,4.2) node {${\displaystyle \mathbb{C} \mapsto \lim_{x_0 \in \Gamma_0} \varrho_{\lambda_0}(x_0)}$}
			(2,-0.2) node {${\displaystyle\mathbb{C} \mapsto \lim_{x_1 \in \Gamma_1} \varrho_{\lambda_1}(x_1)} $}
			(0,2.6)   \dir{90}{90} (4,2.6) ;
			\draw[line width=0.5pt,->]
			(0,1.4)   \dir{270}{270} (4,1.5) 
			;
			\draw (2,3.2) edge[implies] node[left] {$$} (2,0.8);
		\end{tikzpicture} 
	\end{array},
\end{align} whose image under $\Psi$ is the linear map
\begin{align}
	\Par \varrho_{\lambda_0} = \lim_{x_0 \in \Gamma_0} \varrho_{\lambda_0}(x_0) \to  \Par \varrho_{\lambda_1} = \lim_{x_1 \in \Gamma_1} \varrho_{\lambda_1}(x_1)
\end{align} being by Definition~\ref{defparext2} the composition
\begin{align}
	\Par \varrho_{\lambda_0} \stackrel{r_0^* }{\to} \Par r_0^* \varrho_{\lambda_0}  \stackrel{\omega_* }{\to} \Par r_1^* \varrho_{\lambda_1} \stackrel{{r_1}_*}{\to} \Par \varrho_{\lambda_1} .
\end{align} The names chosen for these three maps suggestively already coincide with the corresponding maps used for the parallel section functor in \cite{schweigertwoikeofk}. For the first map (the pullback map) this is indeed obvious. For the second map this follows from the fact that $\omega$ can be seen as an intertwiner $r_0^* \varrho_{\lambda_0} \to r_1^* \varrho_{\lambda_1}$ as observed in the proof of Proposition~\ref{satzendtrepgrpd}. Finally, the fact that ${r_1}_*$ is really given by integral over homotopy fibers of $r_1$ with respect to groupoid cardinality (as the pushforward maps in \cite{schweigertwoikeofk}) follows from the application of Corollary~\ref{korbeckchevgrpdcard} to the square
\begin{align}
	\begin{tikzpicture}[scale=2, implies/.style={double,double equal sign distance,-implies},
		dot/.style={shape=circle,fill=black,minimum size=2pt,
			inner sep=0pt,outer sep=2pt},]
		\node (A1) at (0,1) {$\Omega$};
		\node (A2) at (1,1) {$\star$};
		\node (B1) at (0,0) {$\Gamma_1$};
		\node (B2) at (1,0) {$\star$};
		\path[->,font=\scriptsize]
		(A1) edge node[above]{$$} (A2)
		(A1) edge node[left]{$r_1$} (B1)
		(A2) edge node[right]{$$} (B2)
		(B1) edge node[below]{$$} (B2);
	\end{tikzpicture}.     
\end{align}
\end{proof}

\end{document}